\documentclass{amsart}
\usepackage{xspace,amssymb,amsfonts,amsthm,ifpdf,fge, rotating}
\usepackage{epsfig,euscript,enumerate,cancel}
\usepackage[enableskew]{youngtab}
\usepackage{palatino}
\usepackage[all]{xy}
\usepackage{graphicx}

\newcommand{\PL}{\operatorname{Pl}}

\newcommand{\cA}{{\mathcal A}}

\newcommand{\cF}{{\mathcal F}}

\newcommand{\cH}{{\mathcal H}}

\newcommand{\cN}{{\mathcal N}}

\newcommand{\cS}{{\mathcal S}}
\newcommand{\cT}{{\mathcal T}}
\newcommand{\cU}{{\mathcal U}}

\newcommand{\mg}{\mathfrak{g}}
\newcommand{\mh}{\mathfrak{h}}

\newcommand{\mC}{\mathbb{C}}

\newcommand{\mZ}{\mathbb{Z}}

\newcommand{\la}{\lambda}

\newcommand{\Pfin}{P_{\op{fin}}}
\newcommand{\ten}{10}
\newcommand{\cansix}{\cancel6}
\newcommand{\cannine}{\cancel9}

\newcommand{\canfive}{\cancel5}
\newcommand{\canthree}{\cancel3}

\newcommand{\END}{\operatorname{End}}

\newcommand{\op}{\operatorname}

\newcommand{\inj}{\hookrightarrow}

\newcommand{\surj}{\mbox{$\rightarrow\!\!\!\!\!\rightarrow$}}

\newtheorem{theorem}{Theorem}[section]
\theoremstyle{plain}

\newtheorem{corollary}[theorem]{Corollary}

\newtheorem{definition}[theorem]{Definition}
\newtheorem{example}[theorem]{Example}

\newtheorem{lemma}[theorem]{Lemma}

\newtheorem{proposition}[theorem]{Proposition}
\newtheorem{remark}[theorem]{Remark}

\numberwithin{equation}{section}
\numberwithin{figure}{section}
\input{tcilatex}

\begin{document}
\title[A combinatorial approach to fusion rings]{The $\widehat{\mathfrak{sl}}(n)_{k}$-WZNW Fusion Ring: a combinatorial construction and a realisation as quotient of quantum cohomology}
\author{Christian Korff}
\address{Department of Mathematics, University of Glasgow, 15 University Gardens, Glasgow G12 8QW, Scotland, UK}
\email{c.korff@maths.gla.ac.uk}
\urladdr{http://www.maths.gla.ac.uk/\symbol{126}ck/}
\thanks{C.K. is financially supported by a University Research Fellowship of the Royal Society}
\author{Catharina Stroppel}
\address{Mathematik Zentrum Bonn, Endenicher Allee 60, 53115 Bonn, Germany}
\email{stroppel@uni-bonn.de}
\urladdr{http://www.math.uni-bonn.de/people/stroppel/}

\begin{abstract}
A simple, combinatorial construction of the $\widehat{\mathfrak{sl}}(n)_{k}$-WZNW fusion ring, also known as Verlinde algebra, is given. As a byproduct of the construction
one obtains an isomorphism between the fusion ring and a particular quotient
of the small quantum cohomology ring of the Grassmannian $Gr_{k,k+n}$. We explain how our approach naturally fits into known combinatorial descriptions of the quantum cohomology ring, by establishing what one could call a `Boson-Fermion-correspondence' between the two rings. We also present new recursion formulae for the structure constants of both rings, the fusion coefficients and the Gromov-Witten invariants.
\end{abstract}

\maketitle
\tableofcontents

\section{Introduction}
Given $k, n\in\mathbb{N}$ natural numbers, we can associate the following two rings:
\subsubsection*{The quantum cohomology ring $qH^\bullet(\op{Gr}_{k,n+k})$}
Let $\op{Gr}_{k,n+k}$ be the Grassmannian of $k$-planes in $\mC^{k+n}$, considered as an algebraic variety. Then the {\it quantum cohomology ring} $qH^\bullet(\op{Gr}_{k,n+k})$ is a deformation of the ordinary (integral) cohomology ring $H^\bullet(\op{Gr}_{k,n+k})$ of $\op{Gr}_{k,n+k}$. The latter has a basis given by the Schubert classes $[\Omega_\la]$, indexed by partitions $\la$ whose Young diagram fits into a bounding box of size $k\times n$. The structure constants are just intersection numbers, the so-called Littlewood-Richardson coefficients. The quantum cohomology can be viewed as the free $\mZ[q]$-module with the same basis, but the structure constants $\sum_d C_{\la,\mu}^{\nu^*, d}q^d$ are now the so-called 3-point genus 0 Gromov-Witten invariants, which count the number of rational curves of degree $d$ passing through generic translates of the involved three Schubert cells. By a theorem of Siebert and Tian \cite{ST}, there is an isomorphism of rings $qH^{\bullet }(\op{Gr}_{k,k+n})\cong \mathbb{Z}[q]\otimes_\mZ\mZ[e_1,e_2,\ldots, e_k]/
\left\langle h_{n+1},\ldots,h_{n+k-1},h_{n+k}+(-1)^{k}q\right\rangle$, where $e_i$, $1\leq i\leq k$ denotes the ordinary $i^{\text{th}}$ elementary symmetric polynomial and $h_{r}$, $r\geq 1$, denotes the $r^{\text{th}}$ complete symmetric function in $k$ variables.
\subsubsection*{The fusion ring $\cF(\widehat{\mathfrak{sl}}(n))_k $}
The {\it fusion ring} $\cF=\cF(\widehat{\mathfrak{sl}}(n))_k$ of integrable highest weight representations of $\widehat{\mathfrak{sl}}(n)$ at level $k$ has a natural basis indexed by integrable highest weights $\hat{\la}\in P_k^+$ of level $k$ which can be identified with partitions $\la$ whose Young diagram fits into a bounding box of size $(n-1)\times k$. The structure constants in this ring can be described in various ways, for instance via characters of the Lie algebra ${\mathfrak{sl}}(n)$ (Verlinde formula), in terms of combinatorics of the affine Weyl group (Kac-Walton formula), geometrically as dimensions of certain moduli spaces of generalised theta functions etc.

To motivate our discussion we recall the following theorem due to Witten \cite{Witten}, based on earlier work of Gepner \cite{Gepner} and Vafa \cite{Vafa1} and Intrilligator \cite{Intrilligator1} which states that there is an isomorphism of rings
$$\cF(\widehat{\mathfrak{gl}}(n))_k\cong qH^\bullet(\op{Gr}_{n,k+n})_{q=1}$$
between the level $k$ fusion ring of $\widehat{\mathfrak{gl}}(n)$ (resp. $\widehat{\mathfrak{u}}(n)$) and the specialisation at $q=1$ of the quantum cohomology ring. (A mathematical proof seems to be contained in the unpublished,
  unfortunately not anymore  available, work \cite{Agnihotri}.)

The main goal of this paper is to give a (mathematically rigorous) realisation of  $\cF(\widehat{\mathfrak{sl}}(n))_k$ as a quotient of $qH^\bullet(\op{Gr}_{k,k+n})$ with the defining relations explained via Bethe Ansatz equations of a quantum integrable system.

Both rings will also be described combinatorially in terms of symmetric polynomials (Schur polynomials) in pairwise non-commuting variables. These variables will be interpreted in either case in two ways: as particle hopping operators of a certain quantum integrable system and as generators of an algebra appearing naturally in representation theory. A crucial observation here is that the generating function of the elementary symmetric polynomials  turns out to be the transfer matrix of the quantum integrable system. This is the starting point of our analysis of the two rings.

As a consequence we obtain a simple particle formulation of both rings leading to novel recursion formulae for the structure constants of the fusion ring as well as for the Gromov-Witten invariants. Certain symmetries in these constants follow easily from our quite elementary description of both rings.

The Verlinde formula emerges naturally from our approach and in particular equips the combinatorial fusion algebra with a natural $PSL(2,\mZ)$-action (see Remark \ref{modinv}), a structural data of the abstract Verlinde algebras as introduced in \cite[0.4.1]{Cherednik}.

An abstract isomorphism of the two rings as mentioned above could be of course also obtained much more directly, since both rings are (after complexification) semisimple, and hence it is enough to compare their spectra (see for instance \cite{Beauville} for a description). The interesting new aspect here is our connection with Bethe Ansatz techniques of quantum integrable systems `in the crystal limit'.

\subsubsection*{Outline of the results and methods of the paper}

Let $\cH_k$ denote the underlying vector space of the fusion ring $\cF$. Any basis vector of $\cH_k$ can be realised as a partition, or alternatively as an $n$-tuple $(m_0,m_1,\ldots, m_{n-1})$ of non-negative integers which sum up to $k$. We therefore can consider it as a particle configuration on the circular lattice with $n$ sites and study this underlying integrable model by viewing $\cH_k$ as the space of $k$-particle states. (Details are explained in Section 2.)

\begin{figure}
\includegraphics[scale=0.4]{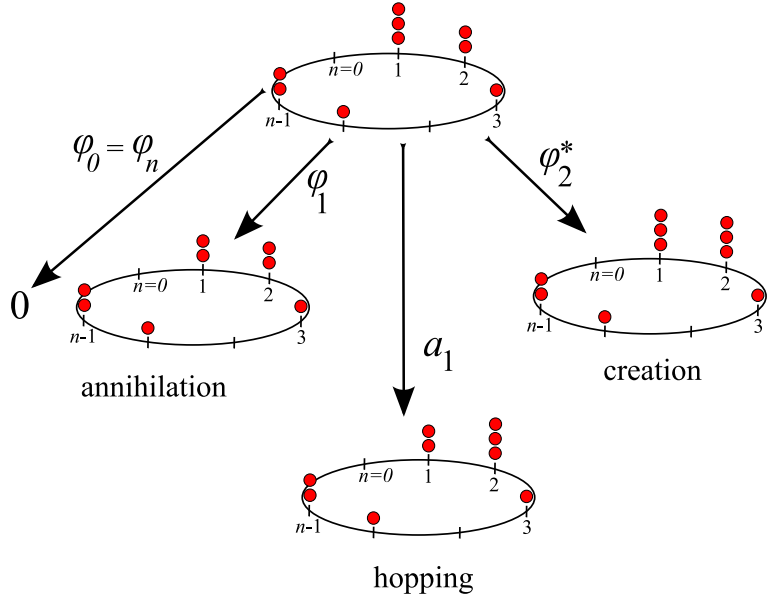}
\caption{The dominant integral weight $(m_0,m_1,\ldots, m_{n-1})=(0,3,2,1,0,1,2)$ of level $k=9$ as particle configuration on a circle; the processes of creating and annihilating, and hopping applied to it.}
\label{fig:bosons}
\end{figure}

We have the obvious operations of particle creation and annihilation $\varphi_i^*$, $\varphi_i$ at site $i$ and {\it particle hopping} $a_i=\varphi_{i+1}^*\varphi_i$; see Figure \ref{fig:bosons}. The endomorphisms $a_i$ of $\cH$ form an interesting algebra which we call the {\it local affine plactic algebra} or the {\it generic affine Hall algebra}, since it generalises the local plactic algebra, a certain quotient (considered by Fomin and Greene [15]) of an algebra originally introduced by Lascoux and Sch\"{u}tzenberger [31]. Its geometric incarnation as generic Hall algebra was introduced by Reineke \cite{Reineke}. These algebras can be viewed as the {\it $q=0$ crystal limit} of the positive part $U_q(\mathfrak{n^+})$ of the quantised universal enveloping algebra of  $\widehat{\mathfrak{sl}}(n)$. The action of the $a_i$'s sends a basis vector to a basis vector or zero and defines a crystal isomorphic to the crystal of the $k$-symmetric power of the natural representation of $U_q(\widehat{\mathfrak{sl}}(n))$. (Details can be found in Section 3 and 4).

We then define elementary symmetric functions $e_r=e_r(a_0,a_1,\ldots, a_{n-1})$ in these non-commuting generators and get the following first result:

\begin{theorem}
The space $\cH_k$ of states can be equipped with the structure of an associative unital algebra $\cF_{\op{comb}}$ by defining
$$\la \circledast \mu=s_\la(a_0,a_1,\ldots a_{n-1})\mu,$$
where $s_\la(a_0,a_1,\ldots ,a_{n-1})$ denotes the (noncommutative) Schur polynomials defined via the usual {\em Jacobi-Trudi formula} $s_{\lambda }=\det \left(e_{\lambda _{i}^{t}-i+j}\right)$. This combinatorially defined fusion product coincides with the fusion product in the Verlinde algebra, hence $\cF_{\op{comb}}\cong\cF(\widehat{\mathfrak{sl}}(n))_k$.
\end{theorem}
 To ensure that the definition of the noncommutative Schur polynomials makes sense, we have to show that the $e$'s pairwise commute. Instead of giving a purely combinatorial proof of this fact \cite{Rohde} we construct a simultaneous eigenbasis for the action of the $s_\la(a_0,a_1,\ldots, a_{n-1})$'s using the algebraic Bethe Ansatz and prove:

\begin{theorem}
\begin{enumerate}
\item The generating function for the noncommutative elementary symmetric functions $e_r(a_0,a_1,\ldots a_{n-1})$ is given by the transfer matrix $T(u)$ associated with a quantum integrable system, known as the {\em phase model}.
\item This system is integrable in the sense that $[T(u),T(v)]=0$ for any $u,v\in\mC$, in particular, the $e_r(a_0,a_1,\ldots a_{n-1})$'s commute, and the $s_\la(a_0,a_1,\ldots a_{n-1})$'s are well-defined.
\item
The Bethe vectors \eqref{bethevec} form an orthogonal eigenbasis for the action of the noncommutative Schur-functions $s_\la(a_0,a_1,\ldots a_{n-1})$.
\item The transformation matrix between the standard basis and the basis of Bethe vectors can be expressed in terms of Weyl characters evaluated at certain $(n+k)^{\text{th}}$ roots of unity. The combinatorial fusion product satisfies the Verlinde formula \eqref{Verlinde} expressed in the entries of the modular S-matrix.
\end{enumerate}
\end{theorem}
Let $\Lambda^{(k)}=\mathbb{Z}[e_{1},\ldots,e_{k}]$ be the ring of symmetric polynomials in $k$ variables. Our main theorem, connecting the fusion algebra with the quantum cohomology ring can be stated as follows:
\begin{theorem}
The assignment $P^+_k\ni\hat{\la}\mapsto s_{\la^t}$ defines an isomorphism of rings
\begin{equation*}
\cF_{\op{comb}}=\cF(\widehat{\mathfrak{sl}}(n),\mathbb{Z})\cong \mathbb{Z}[e_{1},\ldots,e_{k}]/\langle h_n-1,h_{n+1},\ldots,h_{n+k-1},h_{n+k}+(-1)^{k}e_{k}\rangle
,
\end{equation*}
hence, the fusion ring is isomorphic to the quotient of $qH^\bullet(\op{Gr}_{k,n+k})$ obtained by imposing the additional relations $q=e_k$ and $h_n=1$.
\end{theorem}
Note that the algebra structure on the left hand side of the isomorphism is defined in terms of noncommutative Schur polynomials whereas the right hand side is in terms of ordinary commutative Schur polynomials. By Lemma \ref{baeequiv}, the defining relations of the quotient ring coincide exactly with the Bethe Ansatz equations \eqref{bae}. The case $n=3$, $k=1$ is worked out in detail in Example~\ref{example3and1}.

In the second part of the paper we develop the integrable system for the quantum cohomology ring with the resulting combinatorial description of the ring using  noncommutative Schur functions in the generators of the affine nil-Temperley-Lieb algebra, reproving a result of Postnikov \cite{Postnikov} with our methods. The defining relations are then the Bethe Ansatz equations \eqref{baeequiv2} of a free fermion system. The resulting basis of Bethe vectors agrees with a basis introduced by Rietsch \cite{Rietsch} in her comparison of the quantum cohomology ring with the coordinate ring of Peterson varieties. By applying our method of deducing the Verlinde formula on the fusion ring side once more on the quantum cohomology side, gives the celebrated Bertram-Vafa Intrilligator formula \eqref{BVI} expressed in terms of Schur polynomials evaluated at roots of unity (as first developed by Rietsch). Employing our free fermion formulation we present new identities relating Gromov-Witten invariants at different dimensions $k$ leading to an inductive algorithm which starting from $k=0$ allows one to compute and relate the entire hierarchy of structure constants for all $k>0$; see Remark \ref{GWalgo}. We will provide an explicit example. This algorithm differs from the known rim-hook algorithms in \cite{BCF} and \cite{Sottile} (see also \cite{Buch} and references therein) which compute Gromov-Witten invariants for fixed $k$ only.

The paper finishes with Part III, where we summarize the parallel construction in what we shall call {\em Boson-Fermion correspondence} in the present context because of its close analogy with the well-known case (see e.g. \cite{MJD} and \cite{KacRaina}). A detailed study of this correspondence from a representation theory point of view will appear in a forthcoming paper.\cite{KS}\\

{\bf Acknowledgement}:
The first author thanks Michio Jimbo, Ulrich Kr\"ahmer, Eugene Mukhin, Jonathan Nimmo, Arun Ram and Simon Ruijsenaars for helpful discussions. He would like to express his special gratitude to Alastair Craw for his insightful comments, support and advice. The second author would like to thank Evgeny Feigin, Sergey Fomin, Konstanze Rietsch, Christoph Schweigert and Ivan Cherednik for helpful discussions. Both authors are grateful to Ken Brown, Rinat Kedem and in particular Anne Schilling for sharing knowledge and many ideas. We would like to thank the referee for a very thorough reading of the article and helpful comments. Most of the collaboration took place during stays by both authors at
the Isaac Newton Institute in Spring 2009.
We thank the INI staff and both, the Discrete Integrable Systems and the Algebraic Lie Theory programme
organisers for the invitation and for providing the possibility to carry out this research.

\subsubsection*{General conventions}
In the following vector spaces are always defined over the complex numbers $\mC$, by an {\it algebra} we mean an associative unital algebra over $\mC$, and we will abbreviate $\otimes=\otimes_\mC$, $\END=\END_\mC$. To avoid confusion with indices, we denote the imaginary complex number $(0,1)$ by $\iota$.

\section*{Part I: The fusion ring: bosons on a circle}
\section{The affine Dynkin diagram and the vector space of states}
\subsection{The Lie algebra $\widehat{\mathfrak{sl}}(n)$, weights and partitions}
\label{sl}
We start by setting up the notation and basics needed from the (untwisted) affine Lie algebra
$\widehat{\mathfrak{sl}}(n)$ (see e.g. \cite{Kac}, \cite{WZX} for the general theory).

Let $\mg={\mathfrak{sl}}(n)$ be the Lie algebra of complex traceless $n\times n$ matrices with its standard Cartan
subalgebra $\mh$ given by the diagonal matrices. If $\mC[t,t^{-1}]$ denotes the ring of formal Laurent series in $t$ we
have the loop algebra $\mg\otimes\mC[t,t^{-1}]$ with Lie bracket $[x\otimes P, y\otimes Q]=[x,y]\otimes PQ$ for $x,y\in
\mg, P,Q\in \mC[t,t^{-1}]$. This Lie algebra has a unique central extension $\mg\otimes\mC[t,t^{-1}]\otimes\mC{\bf C}$,
and $\hat{\mg}:={\widehat{\mathfrak{sl}}}(n)$ is then obtained by adding an exterior derivative ${\bf D}=t
\frac{d}{dt}$, in formulae
\begin{equation*}
\widehat{\mathfrak{sl}}(n)=\mathfrak{sl}(n)\otimes\mathbb{C}[t,t^{-1}]\oplus \mathbb{C}\boldsymbol{C}\oplus
\mathbb{C}\boldsymbol{D}\;.
\end{equation*}
with Lie bracket $[{\bf D},{\bf C}]=0$,  $[{\bf D},x\otimes P]=x\otimes \frac{d}{dt}P$ where $x\in\mg, P\in
\mC[t,t^{-1}]$. Consider the extended Cartan subalgebra $\hat{\mh}:=\mh\oplus\mC{\bf C}\oplus \mC{\bf D}$ and let
$\delta, \hat{\omega}_0\in\hat{\mh}^*$ such that $\delta({\bf D})=1$,  $\delta({\bf C})=\delta(h)=0$ and
$\hat{\omega}_0({\bf C})=1$, $\hat{\omega}_0({\bf D})=\hat{\omega}_0(h)=0$, for $h\in\mh$.

Let $\{\varepsilon _{i}\}_{i=1}^{n}\in\mh^*\subset\hat{\mh}^*$ be the linear function which picks out the $i^{\text{th}}$
diagonal entry of an element $h\in\mh$. The elements $\alpha _{i}=\varepsilon_{i}-\varepsilon _{i+1}$ for $i=1,\ldots,n-1$
form a basis of $\mh^*$ and can be identified with the simple roots of $\mathfrak{sl(n)}$. Viewed as elements of
$\hat{\mh}^*$ (by extending trivially by zero), together with
$\alpha _{0}=\varepsilon_{n}-\varepsilon _{1}+{\delta }\in\widehat{\mh}^*$ they form exactly the simple roots of
$\widehat{\mg}$. The $\mathbb{Z}$-span of $\{\alpha _{i}\}_{i=0}^{n-1}$ is the affine root lattice $Q$. Let $\alpha_i^\vee$, $0\leq i\leq n$, be the dual roots $\alpha_i^\vee\in\mh$ defined by $\langle \alpha_i, \alpha_j^\vee\rangle=a_{i,j}$, the $(i,j)^{\text{th}}$ entry in the Cartan matrix $A$.

We denote by
$\Gamma$ the Dynkin diagram of $\widehat{\mathfrak{sl}}(n)$, which we view as a circle with $n$ equidistant marked
points. We name these points $0,1,\ldots, n-1$ in clockwise direction (see Figure \ref{fig:bosons}).

There is a unique non-degenerate bilinear pairing on $\mh^*$ given by $(\varepsilon_i,\varepsilon_j)=\delta_{i,j}$ and
we have the fundamental weights
\begin{equation}
\omega _{i}=\varepsilon _{1}+\varepsilon _{2}+\cdots +\varepsilon _{i}-\frac{%
i}{n}\tsum_{j=1}^{n}\varepsilon _{j},\quad i=1,\ldots,n-1  \label{fund weights}
\end{equation}%
for $\mg$ characterized by $\frac{2(\omega_i, \alpha_j)}{(\alpha_j, \alpha_j)}=\langle \omega_i,\alpha_j^\vee\rangle=\delta_{i,j}$ for $i=1,\ldots,n-1$. Let $\Pfin$ denote the corresponding integral weight
lattice, that is the $\mZ$-span of the fundamental weights, with its positive cone $\Pfin^+$ given by the $\mZ_{\geq
0}$-span. Set $\hat{\omega}_i=\omega_i+\hat{\omega}_0$ for $1\leq i\leq {n-1}$.

Then we have the (affine) integral weight lattice $P=\Pfin\oplus\mZ \hat{\omega}_0\oplus\mC\delta$. Now fix
$k\in\mZ_{>0}$ and $\gamma\in\mC$ and consider
\begin{eqnarray}
\label{levelk}
P^+_k=\left\{\hat{\la}=\sum_{i=0}^{n-1}m_i\hat{\omega}_i+\gamma\delta\;\left |\;\sum_{i=0}^{n-1}m_i=k\right.\!,\;m_i\in\mathbb{Z}_{\geq 0}\right\},
\end{eqnarray}
the set of {\it integral dominant weights of level $k$}.  Since (up to a
grading induced by the action of $\delta$) the integrable highest weight modules are independent of the choice of
$\gamma$, the particular choice will not be important for us and we therefore assume from now on $\gamma=0$. The $m_{i}$ appearing here are often called Dynkin labels and will also be denoted $m_i(\hat{\la})$ in the following. For
$\hat{\la}=\la+k\hat{\omega}_0\in P^+_k$ we call $\la$ the {\it finite part} of $\hat{\la}$. Note that $\la=\sum_{i=1}^{n-1}m_i\omega_i$.

It is convenient to encode affine and non-affine weights in terms of
partitions. For any two integers $r,c\in \mathbb{Z}_{\geq 0}$ let
\begin{equation*}
\mathfrak{P}_{\leq r,c}=\{\mu =(\mu _{1},\ldots,\mu_{r},0,0,\ldots)~|~
\mu _{i}\in \mathbb{Z}_{\geq 0},\;\mu _{i}\geq
\mu _{i+1},\;\mu _{1}\leq c\}
\end{equation*}%
which is the set of all partitions whose associated Young diagram has at most $r$ rows and $c$ columns, i.e. they fit
into a boundary box of height $r$ and width $c$. Hereby we use the English notation to denote Young diagrams, we will often just
write $(\mu _{1},\ldots,\mu_{r})$ instead of $(\mu _{1},\ldots,\mu
_{r},0,0,\ldots)$ and use the symbol $\emptyset$ to denote the (unique) partition of zero; so for instance
\begin{eqnarray}
\label{PYoung}
\mathfrak{P}_{\leq 2,3}=\left\{\emptyset, (1), (2), (3),  (1,1), (2,1), (2,2), (3,1), (3,2), (3,3)\right\}
\end{eqnarray}
corresponds to
\begin{eqnarray}
\label{Young}
{\tiny\Yvcentermath1\left\{\emptyset, \yng(1), \yng(2), \yng(3),  \yng(1,1), \yng(2,1), \yng(2,2), \yng(3,1), \yng(3,2),
\yng(3,3)\right\}}
\end{eqnarray}
The following Lemma passes between weights and partitions. For instance, it identifies the Young diagrams in \eqref{Young} via the map $\op{P}^{-1}$ in \eqref{weight2part} with the weights
$$3\hat{\omega}_0, 2\hat{\omega}_0+\hat{\omega}_1,
\hat{\omega}_0+2\hat{\omega}_1,3\hat{\omega}_1,
2\hat{\omega}_0+\hat{\omega}_2,
\hat{\omega}_0+\hat{\omega}_1+\hat{\omega}_2,
\hat{\omega}_0+2\hat{\omega}_2,
2\hat{\omega}_1+\hat{\omega}_2, \hat{\omega}_1+2\hat{\omega}_2,3\hat{\omega}_2.$$

\begin{lemma}
\label{weights}
\begin{enumerate}
\item With the notation from \eqref{levelk}, there is a bijection of sets
\begin{eqnarray}
\op{P}:\quad P_{k}^{+}&\longrightarrow&\mathfrak{P}_{\leq n-1,k}\nonumber\\
\hat{\lambda}&\longmapsto&(\mu_{1},\ldots,\mu
_{n-1},0,\ldots) \text{ with }\mu_{i}-\mu_{i+1}=m_i.\label{weight2part}
\end{eqnarray}
(The associated Young diagram has then exactly
$m_{i}$ columns of length $i$.)
\item There is an injective map
\begin{eqnarray}
\hat{\op{P}}:\quad P_{k}^{+}&\longrightarrow&\mathfrak{P}_{\leq n,k}\nonumber\\
\hat{\lambda}&\longmapsto&(\mu_{1}+m_{0},\ldots,\mu
_{n-1}+m_{0},m_{0},0,\ldots),\text{ with }\mu_{i}-\mu_{i+1}=m_i.\label{affine partition}
\end{eqnarray}
(The Young diagram $\hat{\op{{P}}}(\hat{\lambda})$ is obtained from $\op{P}(\hat\la)$ by adding
$m_{0}=k-\sum_{i>0}m_{i}$ columns with $n$ boxes.)
\end{enumerate}
\end{lemma}
Given $\hat{\la}\in P_{k}^{+}$ we also denote by $\hat{\la}$ the associated Young diagram as in \eqref{affine partition}
and with $\la$ the Young diagram associated to the finite part $\la$ of $\hat{\la}$ as in \eqref{weight2part}.
Given $\mu\in\mathfrak{P}_{\leq n,k}$ one can remove all columns of length $n$ and obtain a partition
$\tilde{\mu}\in\mathfrak{P}_{\leq {n-1},k}$. We denote by $\mu'=\op{P}^{-1}(\tilde{\mu})$, the preimage under $\op{P}$. In the following we identify the elements from $\{0,1,\ldots, n-1\}$ with the elements from $\mZ_n:=\mZ/n\mZ$.
\subsection{Diagram automorphisms and the vector space of states}
\label{autos}
The group of Dynkin diagram automorphisms of $\Gamma$ is generated by the rotation $\op{rot}$ of order $n$ which sends
vertex $i$ to vertex $i+1$ modulo $n$ and the diagram automorphism $\op{flip}$ coming from the non-affine Dynkin diagram
which fixes the $0^{\text{th}}$ vertex and otherwise maps vertex $i$ to vertex $n-i$. There are corresponding automorphisms
$P^+_{k}\rightarrow P^+_{k}$:
\begin{eqnarray}
\hat{\lambda}=\tsum_{i\in \mathbb{Z}_{n}}m_{i}\hat{\omega}_{i}&\longmapsto&
\op{rot}(\hat{\lambda}):=\tsum_{i\in \mathbb{Z}_{n}}m_{i}\hat{\omega}%
_{i-1}\label{Dynkin auto}\\
\hat{\lambda}=k\hat{\omega}_{0}+\tsum_{i=1}^{n-1}m_{i}\omega _{i}&\longmapsto&
\op{flip}(\hat{\lambda}):=k\hat{\omega}_{0}+\tsum_{i=1}^{n-1}m_{i}%
\omega _{n-i}  \label{Dynkin auto 2}
\end{eqnarray}%
In terms of partitions (via \eqref{weight2part}) the rotation corresponds to adding a
row of $k$ boxes and then removing all columns which contain $%
n$ boxes. The action of $\op{flip}$ is best described in terms of taking complements of
partitions: for a partition $\mu$ fitting into a box of height $n$ and width $k$ we denote by $\mu^\vee$ the complementary partition of $\mu$ in this box. Then $\op{flip}(\hat\la)$ gets mapped to the elements
\begin{equation}\label{complement}
\hat\la^\ast:=\op{\hat P}(\op{flip}(\hat\la))=(\op{P}(\hat\la))^\vee\quad\text{and}\quad
\la^\ast:=\op{P}(\op{flip}(\hat\la)).
\end{equation}
For a partition $\la$ we denote by $\la^{t}$ its transpose partition (in terms of Young diagrams it just means we
reflect the diagram in the diagonal).

\subsection{The vector space of states and the phase algebra}
We now introduce what we call the vector space of states. Instead of considering only a fixed
level $k$ it turns out to be more convenient to allow $k$ ranging over the positive integers.
That is, we consider the infinite-dimensional vector space
\begin{equation}
\cH=\tbigoplus_{k\in \mathbb{Z}_{\geq 0}}\cH_{k},\qquad
\cH_{k}=\mathbb{C}P^+_{k}\;,
\end{equation}%
with $\cH_{0}=\mathbb{C\{\emptyset \}}=\mathbb{C}$.
For each $0\leq i\leq n-1$ we define linear endomorphisms of $\cH$ by the following assignment on basis vectors
\begin{eqnarray}
\varphi _{i}^{\ast }:\quad P^+_{k}\rightarrow P^+_{k+1}&&\varphi _{i}^{\ast }\hat{\lambda}=\hat{\lambda}+\hat{\omega%
}_{i}\\
\varphi _{i}:\quad P^+_{k}\rightarrow P^+_{k-1}\cup\{0\}&&\varphi _{i}\hat{\lambda}=
\begin{cases}
\hat{\lambda}-\hat{\omega}_{i}, & \text{if $\hat{\lambda}-\hat{\omega}_{i}\in P^+_{k-1}$,}
\\
0, & \text{otherwise}%
\end{cases}
\label{phase1}
\end{eqnarray}%

In terms of partitions the map $\varphi _{i}^{\ast }$ acts on the Young diagram
associated with $\lambda $ by adding a column with $i$ boxes in the
appropriate place. In contrast, $\varphi _{i}$ is the map which deletes from $%
\lambda $ a column with $i$ boxes or, if it has none, sends $\la$ to zero. For $i=0$ the corresponding maps simply increase or decrease the width of the bounding box by one, provided this is allowed by the shape of the diagram.

With the notation from \eqref{levelk} and $0\leq i\leq n-1$ let $N_{i}$ be the linear endomorphism of $\cH$ which
multiplies every basis vector by $m_i$, in formulae: $N_{i}\hat{\lambda}=m_{i}(\hat{\lambda})\hat{\lambda}$.
The subalgebra of $\END\cH$ generated by $\{\varphi
_{i},\varphi _{i}^{\ast },N_{i}\}$ has been introduced previously in the
physics literature and is called the {\it phase algebra}; compare with \cite{Bogoliubovetal}. The operators $\varphi
_{i}$ and $\varphi _{i}^{\ast }$ can be interpreted as
respectively particle {\it annihilation} and {\it creation operators} at site $i$ of a
circular lattice which coincides with the Dynkin diagram of $\widehat{%
\mathfrak{sl}}(n)$. The Dynkin label $m_{i}(\hat{\lambda})$ is the occupation number at site $i$ and thus
$N=\sum_{i}N_{i}$ is the {\it total particle number operator}, and the subspace $\cH_{k}$ corresponds to the physical
states which contain
$k$ particles. The map $\pi _{i}=1-\varphi _{i}^{\ast }\varphi _{i}$ projects onto
the subspace where no particle is sitting at site $i$. In the following we will denote $\varphi_0$ also  by
$\varphi_n$, and similarly $\varphi_0^*$ also  by
$\varphi_n^*$. We refer to Section~\ref{phasemodel} for the so-called phase model and to Figure \ref{fig:bosons} for an illustration.

\section{The phase algebra and the quantum Yang-Baxter-algebra and their connection to integrable systems}
\subsection{The phase algebra (creation and annihilation of particles)}
\begin{proposition}[Phase algebra]
\label{phasealgebra}
The $\varphi_i,\varphi^\ast_i$ and $N_i$ generate a subalgebra $\hat{\Phi}$ of $\END(\cH)$ which can be realized as the algebra $\Phi$ with the following generators and relations for
$0\leq i,j\leq n-1$:
\begin{gather}
\varphi _{i}\varphi _{j}=\varphi _{j}\varphi _{i},\quad \varphi _{i}^{\ast }\varphi _{j}^{\ast
}=\varphi _{j}^{\ast }\varphi _{i}^{\ast },\quad N_{i}N_{j}=N_{j}N_{i}\label{comm}\\
N_{i}\varphi _{j}-\varphi _{j}N_{i}=-\delta _{ij}\varphi _{i},\quad N_{i}\varphi
_{j}^{\ast }-\varphi _{j}^{\ast }N_{i}=\delta _{ij}\varphi _{i}^{\ast },\label{comm2}\\
\varphi
_{i}\varphi _{i}^{\ast }=1,\quad\varphi
_{i}\varphi _{j}^{\ast }=\varphi _{j}^{\ast }\varphi_i\;\text{ if }\; i\neq j,\label{comm3}\\
N_{i}(1-\varphi _{i}^{\ast }\varphi _{i})=0=(1-\varphi _{i}^{\ast }\varphi _{i})N_{i}. \label{Npi}
\end{gather}%

The algebra $\Phi$ has a basis $B$ of the form
\begin{eqnarray}
\label{PBW}
\{B_{\bf{ b, a, c}}:=
{\varphi_0^*}^{b_0}{\varphi^*_1}^{b_1}\cdots
{\varphi_{n-1}^*}^{b_{n-1}}\varphi_0^{a_0}\varphi_1^{a_1}\cdots\varphi_{n-1}^{a_{n-1}}N_0^{c_0}N_1^{c_1}\cdots N_{n-1}^{c_{n-1}}\},
\end{eqnarray}
where $a_i, b_i, c_i\in\mZ_{\geq 0}, a_ib_ic_i=0$ for $0\leq i\leq n-1$. If we introduce the scalar product on the vector space $\cH$ by
$$\langle \alpha\hat{\lambda},\beta\hat{\mu}\rangle =\overline{\alpha}\beta\tprod_{i=0}^{n-1}\delta _{m_{i}(%
\hat{\lambda}),m_{i}(\hat{\mu})},$$
for $\alpha,\beta\in\mC$, then
$\langle \varphi _{i}^{\ast }\hat{\lambda},\hat{\mu}\rangle =\langle \hat{\lambda%
},\varphi _{i}\hat{\mu}\rangle \;.$
\end{proposition}
\begin{remark}{\rm
One can easily check that with $\pi _{i}=1-\varphi _{i}^{\ast }\varphi _{i}$ the following (probably better known)
relations hold: $\pi _{i}\varphi _{i}^{\ast }=\varphi _{i}\pi _{i}=0, \pi _{i}^{2}=\pi _{i}$.}
\end{remark}

\begin{proof}
Straightforward calculations show that the asserted relations hold in the algebra $\hat{\Phi}$. Hence, $\Phi$ surjects canonically onto $\hat{\Phi}$. The commutator relations, the
relation $\varphi_i\varphi_i^*=1$ and \eqref{Npi} imply that the elements in the proposed basis $B$ at least span $\Phi$. To
see that they are linearly independent, we look at their action on $\cH$.
For the following argument we write a basis vector $\hat{\la}$ of $\cH$ (see \eqref{levelk}) as a tuple $\hat{\la}=(m_0,m_1,\ldots, m_{n-1})={\bf m}$.
Assume we have a finite linear combination
\begin{eqnarray}
\label{lincombX}
Z:=\sum_{{\bf b}, {\bf a}, {\bf c}}\gamma_{{\bf b}, {\bf a}, {\bf c}}B_{{\bf b}, {\bf a}, {\bf c}}=0
\end{eqnarray}
with all $\gamma_{{\bf b}, {\bf a}, {\bf c}}\in\mC$ non-zero. We have to produce a contradiction.
Obviously, $Z\hat{\la}=0$ for any $\hat{\la}={\bf m}$. If we choose the $m_i$'s big enough then
$$Z\hat{\la}=\sum_{{\bf b}, {\bf a}, {\bf c}}\gamma_{{\bf b}, {\bf a}, {\bf c}}m_0^{c_0}m_1^{c_1}\ldots m_{n-1}^{c_{n-1}}(m_0-a_0+b_0,m_1-a_1+b_1,\ldots, m_{n-1}-a_{n-1}+b_{n-1}).$$
The coefficient of $(m_0-x_0,m_1-x_1,\ldots, m_{n-1}-x_{n-1})$ for fixed ${\bf x}=(x_0,x_1,\ldots, x_{n-1})$ is equal to
\begin{eqnarray}
\label{coeffs}
\sum_{({\bf b}, {\bf a}, {\bf c})\in X}\gamma_{{\bf b}, {\bf a}, {\bf c}}m_0^{c_0}m_1^{c_1}\ldots m_{n-1}^{c_{n-1}}
\end{eqnarray}
where the sum runs over the set $X$ of all triples with ${\bf a}-{\bf b}={\bf x}$. All these coefficients have to be zero. In other words, the polynomial
\begin{eqnarray}
\label{poly}
P(t_0,t_1,\ldots, t_{n-1})&=&\sum_{({\bf b}, {\bf a}, {\bf c})\in X}\gamma_{{\bf b}, {\bf a}, {\bf c}}t_0^{c_0}t_1^{c_1}\ldots t_{n-1}^{c_{n-1}}
\end{eqnarray}
satisfies $P(m_0,m_1,\ldots, m_{n-1})=0$ whenever the $m_i$'s are big enough. Now consider $P$ as a polynomial $\tilde{P}$ in $t_0$ by putting fixed values $t_i=m_i$ for $i>0$. Then $\tilde{P}$ has infinitely many zeroes (all the big enough $m_0$), so it is constant zero, and each coefficient of $t_0^{c_0}$ for fixed $c_0$ has to be zero. Repeating this argument finally implies that $P$ is the zero polynomial. (Alternatively, $P$ vanishes on the set of all the ${\bf m}$'s for big enough $m_i$'s. They form a Zariski dense set in the affine space $\mathbb{A}^n$, hence $P$ has to be the zero polynomial).  In particular we can fix ${\bf c}$, consider the set $X({\bf c})$ of all triples $({\bf b}, {\bf a}, {\bf c})$ in $X$ with ${\bf c}$ our fixed choice and get
\begin{eqnarray}
\label{big}
\sum_{({\bf b}, {\bf a}, {\bf c})\in X({\bf c})}\gamma_{{\bf b}, {\bf a}, {\bf c}}=0.
\end{eqnarray}
In case the components $c_i$ of ${\bf c}$ are all non-zero, then the condition $a_ib_ic_i=0$ implies that for any $i$, $a_i=0$ or $b_i=0$, and so there is a unique element in $X({\bf c})$ (namely $x_i=b_i$, $a_i=0$ if $x_i\geq 0$ and $x_i=a_i$, $b_i=0$ if $x_i\leq 0$), because $a_i, b_i\geq0$. Therefore, $\gamma_{{\bf b}, {\bf a}, {\bf c}}=0$. Hence all the in $Z$ occurring ${\bf c}$'s contain at least one zero. For $\gamma_{{\bf b}, {\bf a}, {\bf c}}$ occurring in $Z$ we let $l(\gamma_{{\bf b}, {\bf a}, {\bf c}})$ be the number of zeroes in ${\bf c}$. Let $l$ be the minimum of all these $l(\gamma_{{\bf b}, {\bf a}, {\bf c}})$.  From above we know $l>0$.

Now we choose some ${\bf c}$ occurring in $Z$ with exactly $l$ zeroes. Call it $\tilde{\bf c}$. Without restriction we may assume $\tilde{c}_0=\tilde{c}_1=\ldots =\tilde{c}_{l-1}=0$. Consider the $\gamma_{{\bf b}, {\bf a}, \tilde{\bf c}}$'s appearing in \eqref{lincombX}. Amongst these pick the ones with $a_0$, call it $\tilde{a}_0$, minimal. Amongst these $\gamma_{{\bf b}, {\bf a}, \tilde{\bf c}}$ with $a_0=\tilde{a}_0$ choose $a_1$ minimal, call it $\tilde{a}_1$, etc. Carrying on like this defines  $\tilde{a}_i$ for $0\leq i\leq l-1$. Set ${\bf m}=(\tilde{a}_0, \tilde{a}_1, \ldots, \tilde{a}_{l-1}, m_l,\ldots, m_{n-1})$ and abbreviate $$M({\bf a}, {\bf b},{\bf c})=(\tilde{a}_0-a_0+b_0,\dots\tilde{a}_{l-1}-a_{l-1}+b_{l-1}, m_l-a_l+b_l,\ldots, m_{n-1}-a_{n-1}+b_{n-1}).$$
If we choose the $m_i$'s big enough then
$$0=Z{\bf m}=\sum_{{\bf b}, {\bf a}, {\bf c}}\gamma_{{\bf b}, {\bf a}, {\bf c}}\tilde{a}_0^{c_0}\ldots\tilde{a}_{l-1}^{c_{l-1}}m_l^{c_l}\ldots m_{n-1}^{c_{n-1}} M({\bf a}, {\bf b},{\bf c}),$$
where the sum runs over all triples with the extra condition $a_i\leq\tilde{a}_i$ for $0\leq i\leq l-1$.
As in \eqref{coeffs} we consider the coefficient of a fixed basis vector $\hat{\la}$  and use the polynomial argument from \eqref{poly} to deduce that
\begin{eqnarray}
\label{sumX}
\sum_{{\bf b}, {\bf a}, {\bf c}\in X(\tilde{c},a_i\leq\tilde{a}_i )}\tilde{a}_0^{c_0}\tilde{a}_1^{c_1}\ldots\tilde{a}_{l-1}^{c_{l-1}}\gamma_{{\bf b}, {\bf a}, {\bf c}}=0,
\end{eqnarray}
where the set $X(\tilde{c},a_i\leq\tilde{a}_i)$ consists of all triples with $a_i\leq\tilde{a}_i$ for $0\leq i\leq  l-1$ and $c_i=\tilde{c}_i$ for $l\leq i\leq n-1$ as well as ${\bf a}-{\bf b}=\bf{x}$ fixed. The minimality of $l$ implies however that $\gamma_{{\bf b}, {\bf a}, {\bf c}}\not=0$ and $c_i=\tilde{c}_i$ for $l\leq i\leq n-1$ forces $\bf{c}=\tilde{\bf c}$, and so we must have ${\bf c}=\tilde{\bf c}$. Now \eqref{sumX} refines \eqref{big} in the following sense
\begin{eqnarray}
\label{sumXX}
\sum_{({\bf b}, {\bf a}, \tilde{\bf c})\in X(\tilde{\bf c})}\gamma_{{\bf b}, {\bf a}, {\bf c}}=0,
\end{eqnarray}
with $a_i\leq\tilde{a}_i$ for $0\leq i\leq l-1$. The minimality property of the $\tilde{a}_i$ and the choice of ${\bf c}$ allows to replace the above condition $a_i\leq\tilde{a}_i$ by $a_i=\tilde{a}_i$ for $0\leq i\leq  l-1$. The conditions $c_i\not=0$ for $l\leq i\leq n-1$ and $a_ib_ic_i=0$ imply that for any such $i$, $a_i=0$ or $b_i=0$. Therefore, the sum \eqref{sumXX} has only one summand $\gamma_{{\bf b}, {\bf a}, {\bf c}}$ and we get $\gamma_{{\bf b}, {\bf a}, {\bf c}}=0$. This is a contradiction. Therefore the elements of $B$ are linearly independent and \eqref{PBW} is in fact a basis of $\Phi$.
The action of $\Phi$ on $\cH$ is therefore faithful and factors through $\hat\Phi$. Hence the composition $\Phi\stackrel{\op{can}}\surj\hat{\Phi}\inj\END(\cH)$ is injective, hence $\op{can}$ is an isomorphism. Finally the last claim of the proposition follows directly from the definitions.
\end{proof}

\subsection{The quantum Yang-Baxter algebra and Pieri rules}

Using the phase algebra we want to describe a solution to the quantum
Yang-Baxter equation arising in \cite{Bogoliubovetal}. To calculate with endomorphisms of $\mC^2\otimes\cH$ we use the
abbreviation
\begin{equation*}
\begin{pmatrix}
f_1&f_2\\
f_3&f_4
\end{pmatrix}
:=
\left(
\begin{array}{cc}
1 & 0 \\
0 & 0%
\end{array}%
\right) \otimes f_1+\left(
\begin{array}{cc}
0 & 1 \\
0 & 0%
\end{array}%
\right) \otimes f_2+\left(
\begin{array}{cc}
0 & 0 \\
1 & 0%
\end{array}%
\right) \otimes f_3+
u\left(
\begin{array}{cc}
0 & 0 \\
0 & 1%
\end{array}%
\right) \otimes f_4,
\end{equation*}
where the left hand side is a $2\times 2$ matrix with entries in $\END(\cH)$.
It is straightforward to check that the composition of endomorphisms corresponds to the usual matrix multiplication.

For $i\in\{1,2,\ldots,n\}$, the $i^{\text{th}}$ {\it Lax matrix} or $L$-operator $L_{i}=L_{i}(u)$ is the following one-parameter
family of
endomorphisms of $\mathbb{C}^{2}\otimes \cH$
\begin{equation}
L_{i}(u)=\left(
\begin{array}{cc}
1 & u\varphi _{i}^{\ast } \\
\varphi _{i} & u1%
\end{array}%
\right) \in \limfunc{End}(\mathbb{C}^{2}\otimes \cH)\;.  \label{Lax}
\end{equation}%
The complex variable $u\in \mathbb{C}$ is called the {\it spectral parameter}.
The {\it monodromy
matrix} is defined as
\begin{equation}
M(u)=L_{n}(u)\cdots L_{1}(u)=\left(
\begin{array}{cc}
A(u) & B(u) \\
C(u) & D(u)%
\end{array}%
\right) \in \limfunc{End}(\mathbb{C}^{2}\otimes \cH).  \label{mom}
\end{equation}%
 If we identify $\mC^4$ with $\mC^2\otimes \mC^2$ by mapping the standard basis $e_1,e_2,e_3,e_4$ to $e_1\otimes e_1$,
 $e_1\otimes e_2$, $e_2\otimes e_1$, $e_2\otimes e_2$, then we get
\begin{lemma}[cf. \cite{Bogoliubovetal}]
The monodromy matrix is a solution  to the following $RTT$-relation
in $\limfunc{End}(\mC^2\otimes \mC^2\otimes\cH)$
\begin{equation}
R_{12}(u/v)M_{1}(u)M_{2}(v)=M_{2}(v)M_{1}(u)R_{12}(u/v),\qquad u,v,x\in
\mathbb{C}, u\not=v\not=0, \label{qYBE}
\end{equation}%
with%
\begin{equation}
R(x)=\left(
\begin{array}{cccc}
\frac{x}{x-1} & 0 & 0 & 0 \\
0 & 0 & \frac{x}{x-1} & 0 \\
0 & \frac{1}{x-1} & 1 & 0 \\
0 & 0 & 0 & \frac{x}{x-1}%
\end{array}%
\right) \in \limfunc{End}\mathbb{C}^{4}\cong \limfunc{End}(\mathbb{C}%
^{2}\otimes \mathbb{C}^{2})\;.  \label{Rmatrix}
\end{equation}
The lower indices indicate in which of the two $\mC^2$-spaces of the tensor product the respective operators act.
\end{lemma}
\begin{proof}
It is easy to check that the equation holds when we replace $M$ by any $(L_i)_1$ or $(L_i)_2$. Then the definition of $M$ and the
fact that the $(L_i)_1$ and $(L_j)_2$ pairwise commute for $i\not=j$ imply the claim.
\end{proof}

For $F$ from $A,B,C,D$ as in \eqref{mom} we introduce the power series decomposition $F(u)=\sum_{r\geq 0}F_{r}u^{r}$ with respect to the
spectral parameter $u$. Note that it follows from (\ref{Lax}) and (\ref{mom}) that
$A_{r},B_{r},C_r, D_r=0$ for $r>n$.

\begin{lemma}
\label{ABCD}
The monodromy matrix elements are explicitly given by
\begin{equation*}
\begin{array}[t]{lcllcl}
A(u)&=&\sum_{0\leq r\leq n/2}\sum u^{(j_{1}-i_{1})+\cdots +(j_{r}-i_{r})}a_{i_{1},j_{1}}\cdots
a_{i_{r},j_{r}},&C(u)&=&\varphi_{n}A(u),\\
D(u)&=&\sum_{0\leq r\leq n/2}\sum u^{n-j_{1}+i_{1}-\cdots -j_{r}+i_{r}}a_{j_{1},i_{1}}\cdots
a_{j_{r},i_{r}},&B(u)&=&D(u)\varphi_{n}^\ast,
\end{array}
\end{equation*}
where $a_{i,j}=\varphi_i\varphi_j^\ast$, and the not specified sums run through the set of tuples $0\leq i_{1}<j_{1}<\ldots<i_{r}<j_{r}\leq
n-1$.
\end{lemma}
\begin{proof}
For $n=2$ (or $n=1$) this formula is clear by \eqref{mom}. In general it follows from an easy induction on $n$.
\end{proof}

\begin{definition}{\rm
The {\it quantum Yang-Baxter algebra} is the algebra $\mathfrak{B}$ generated by the $
\{A_{r},B_{r},C_{r},D_{r}\}_{r\geq 0}$ subject to the commutation relations \eqref{qYBE} via the decompositions
$F(u)=\sum_{r\geq 0}F_{r}u^{r}$.
}
\end{definition}

\begin{remark}\label{coprod}{\rm
The Yang-Baxter algebra can be equipped with a bialgebra structure with comultiplication  $\Delta (A) =A\otimes A+C\otimes B$, $\Delta (B) =B\otimes A+D\otimes B$, $\Delta (C) =A\otimes C+C\otimes D$,  $\Delta (D) =B\otimes C+D\otimes D$ and co-unit $\varepsilon (A)=\varepsilon (D)=1$, $\varepsilon
(C)=\varepsilon (B)=0$. The above construction resembles the RTT-construction of Yangians. However, in our case the $L_i(0)$ are not invertible, so that for instance the usual construction \cite[(1.27)]{Molev} of the antipode is not applicable. One can in fact  show that  the bialgebra structure does not extend to a Hopf algebra structure.}
\end{remark}

\begin{remark}{\rm
The term `quantum' Yang-Baxter algebra has its origin in the physical interpretation of the integrable model which is underlying our construction. By forming the state space $\mathcal H$ of particle configurations on a circle we allow for their complex linear superpositions which is the hallmark of a quantum mechanical system in physics. The adjective `quantum' sets our construction apart from algebraic constructions connected with the so-called classical Yang-Baxter equation which differs from relation \eqref{qYBE}.}
\end{remark}

Expanding the quantum Yang--Baxter equation (or RTT-equation) leads for instance to the identities
\begin{eqnarray}
(u-v)A(u)B(v) &=&vB(u)A(v)-vB(v)A(u),  \label{aba1} \\
\left( u-v\right) D(u)B(v) &=&uB(v)D(u)-vB(u)D(v)  \notag\\
C(u)B(v)&=&\frac{v}{u-v}[A(v)D(u)-A(u)D(v)]\ .  \label{aba2}
\end{eqnarray}

The action of the phase algebra on the state space $\cH$ induces an action of the Yang-Baxter algebra. We describe this
action again combinatorially. Let $\la$ and $\mu$ be partitions which we identify with their Young diagrams. Assume that the diagram $\la$ contains
the diagram $\mu$. Then the {\it skew diagram} $\la/\mu$  is obtained by removing $\mu$ from $\la$. It is a {\it
vertical strip} if it contains at most 1 box in each row, or equivalently if $0\leq\lambda _{i}-\mu _{i}\leq 1$. We call
it a  {\it vertical $r$-strip}, denoted $(1^r)$, if it is a vertical strip containing exactly $r$ boxes. Horizontal strips are defined analogously.

\begin{proposition}[Pieri-type formulae]
\label{Pieri}
The space of states $\cH$ can be be turned into a $\mathfrak{B}$-module such that the action of the generators on the
basis is given as follows:
\begin{eqnarray*}
A_{r}\hat{\mu}=\tsum_{\substack{ \hat{\lambda}-\hat{\mu}=(1^{r})  \\ \hat{%
\lambda}\in \cH_{k}}}\hat{\lambda},\quad B_{r}\hat{\mu}=\tsum
_{\substack{ \hat{\lambda}-\hat{\mu}=(1^{r})  \\ \hat{\lambda}\in \cH_{k+1}}}\hat{\lambda}~,\quad
C_{r}\hat{\mu}=\tsum_{\substack{ \hat{\mu}-\hat{\lambda}=(1^{r})  \\ \hat{%
\lambda}\in \cH_{k-1}}}\hat{\lambda},\quad D_{r}\hat{\mu}=\tsum
_{\substack{ \hat{\mu}-\hat{\lambda}=(1^{r})  \\ \hat{\lambda}\in \cH_{k}}}\hat{\lambda}\;.
\end{eqnarray*}
In particular, $B$ increases, whereas $C$ decreases the level.
\end{proposition}

\begin{proof}
The claim can be checked easily by induction using the comultiplication in Remark \ref{coprod}.  Alternatively one may use the explicit formulae from Lemma~\ref{ABCD}.
\end{proof}

The above action of the Yang-Baxter algebra on $\cH$ can be described in terms of skew
Schur functions. We first recall the necessary notions to explain the result and refer to \cite{FultonYT} or
\cite{MacDonald} for more details.
Given a Young diagram $\la$ a {\it semi-standard tableau} (or just {\it tableau}) is a filling of the boxes of $\la$ with
the numbers from $\{1,2, \ldots, n\}$ such that the entries are strictly increasing downwards along the columns and
weakly increasing to the right along the rows. Given a tableau $T$ we have the associated monomial
$x^T:=x_1^{q_1}x_2^{q_2}x_3^{q_3}\ldots x_n^{q_n}$ where $q_i$ denotes the number of boxes filled  with $i$.
The Schur polynomial $s_\la$ is then the sum $\sum x^T$, where $T$ runs through all semi-standard tableaux of shape
$\la$.
\begin{corollary}
Let $\hat{\mu}\in \mathfrak{P}_{\leq n,k}$ be the partition associated with an affine weight of level $k$ via
(\ref{affine partition}). For any $x_{1},\ldots,x_{\ell }$ complex numbers or invertible formal variables we have%
\begin{eqnarray*}
A(x_{1})\cdots A(x_{\ell })\hat{\mu} &=&\sum_{\hat{\lambda}\in \op{P}^+_{k}}s_{%
\hat{\lambda}^{t}/\hat{\mu}^{t}}(x_{1},\ldots,x_{\ell })\hat{\lambda}, \\
B(x_{1})\cdots B(x_{\ell })\hat{\mu} &=&\sum_{\hat{\lambda}\in \op{P}^+_{k+\ell
}}s_{\hat{\lambda}^{t}/\hat{\mu}^{t}}(x_{1},\ldots,x_{\ell })\hat{\lambda}, \\
C(x_{1})\cdots C(x_{\ell })\hat{\mu} &=&(x_{1}\cdots x_{\ell })^{n}\sum_{%
\hat{\lambda}\in \op{P}^+_{k-\ell }}^{n}s_{\hat{\mu}^{t}/\hat{\lambda}%
^{t}}(x_{1}^{-1},\ldots,x_{\ell }^{-1})\hat{\lambda}, \\
D(x_{1})\cdots D(x_{\ell })\hat{\mu} &=&(x_{1}\cdots x_{\ell })^{n}\sum_{%
\hat{\lambda}\in \op{P}^+_{k}}s_{\hat{\mu}^{t}/\hat{\lambda}%
^{t}}(x_{1}^{-1},\ldots,x_{\ell }^{-1})\hat{\lambda}\ .
\end{eqnarray*}
Note that for $\hat{\mu}=\emptyset $ we get the usual Schur polynomials
\begin{equation}
B(x_{1})\cdots B(x_{\ell })\hat{\mu}=\sum_{\hat{\lambda}\in \op{P}^+_{\ell }}s_{%
\hat{\lambda}^{t}}(x_{1},\ldots,x_{\ell })\hat{\lambda}\ .  \label{BetheSchur}
\end{equation}
\end{corollary}

\begin{proof}
Each semi-standard tableau $T$ determines uniquely a sequence of partitions
$\mu =\lambda ^{(1)},\lambda ^{(2)},\ldots,\lambda ^{(\ell )}=\lambda $ where $\mu_i$ is obtained from $\la$ by removing all
boxes in $T$ filled with a number greater than $i$. Semi-standardness implies that $\lambda ^{(r+1)}/\lambda ^{(r)}$ is
a horizontal strip. Conversely, every sequence of partitions differing by horizontal strips arise from a semi-standard tableau in this way. Applying
Proposition \ref{Pieri} yields the desired result.
\end{proof}

\section{Transfer matrix and the phase model}
\label{phasemodel}

We now employ the quantum Yang-Baxter algebra to define a discrete quantum
integrable system, called the {\it phase model} in \cite{Bogoliubovetal} due to its
similarity to constructions in quantum optics. Our approach is motivated by this physical model with the
so-called transfer matrix playing a central role in our combinatorial construction: as we
will see below it is the generating function of the cyclic noncommutative
elementary symmetric polynomials (Proposition \ref{genfunction}).\smallskip

First we extend scalars of the vector space of states from $\mC$ to $\mC[z]$, the ring of polynomials in one variable
and denote it $\cH[z]:=\mathbb{C}[z]\otimes _{\mathbb{C}}\cH$.
In Section~\ref{sec:qcoho}, when the quantum cohomology comes into the picture, $z$ will play a role analogous to the
deformation parameter $q$ in the small quantum cohomology ring. In physical applications it is a magnetic flux parameter (or number) related to
quasi-periodic boundary conditions. It enters the following definition of
the one-parameter family of row-to-row transfer matrices,
\begin{equation}
T(u)=\sum_{r\geq 0}T_{r}u^{r}=A(u)+zD(u)\;.  \label{tcoeff}
\end{equation}

The term `transfer matrix' has its origin in yet another physical interpretation of the phase model. Namely, consider an $m\times n$ square lattice with periodic boundary conditions in both directions, i.e. a toroidal lattice. Fix a particle configuration of our circular lattice with $n$ sites by choosing a partition. Then the operator T maps, `transfers', this configuration into a linear combination of other configurations weighed with its matrix elements which for a special value of the spectral parameter are interpreted as Boltzmann weights, i.e. statistical probabilities that a particular configuration occurs. Taking the $m^{th}$ power of the transfer matrix we end up with a cylinder of height $m$ and by taking its trace we compute the so-called partition function,
$$ Z(u)={\limfunc{Tr}}_{\cH[z]}T(u)^{m}\;,$$
which is the sum over all allowed configuration on the torus and which is the central physical object when interpreting the phase model as a system in statistical mechanics rather than quantum mechanics.

Employing the Yang-Baxter equation one easily shows (the original idea of the proof goes back to Baxter, but the argument can also be found in e.g. \cite[\S3, Lemma 2]{FrenkelResh}) that
\begin{equation}
\lbrack T(u),T(v)]=0  \label{int}
\end{equation}%
for any pair of spectral parameters $u,v\in \mathbb{C}$. Hence the model is {\it integrable}, as it possesses an infinite
number of conserved
quantities.

Alternatively, one can define a quantum system in terms of the Hamiltonian %
\begin{equation}
H=-\frac{T_{1}+T_{n-1}}{2}=-\frac{1}{2}\sum_{i=1}^{n}\left( \varphi _{i}\varphi
_{i+1}^{\ast }+z\varphi _{i}^{\ast }\varphi _{i+1}\right) ,
\end{equation}%
where $\varphi _{n+1}=z^{-1}\varphi _{1}$ and $\varphi _{n+1}^{\ast }=z\varphi _{1}^{\ast }$%
, i.e. one considers particles moving on a circle of $n$ sites. More
generally, we can define higher Hamiltonians, conserved charges by setting,%
\begin{equation}
H_{r}^{\pm }=-\frac{T_{r}\pm T_{n-r}}{2}\;.  \label{imotion}
\end{equation}%
The latter are in involution, $[H_{r}^{\pm },H_{r^{\prime }}^{\pm }]=0$ for
any $1\leq r,r^{\prime }\leq (n-1)/2$, and again we have an integrable
model. However, it is not the existence of these integrals of motions, but
rather the quantum Yang-Baxter algebra which allows one to solve the model
explicitly, i.e. to compute the eigenstates and eigenvalues of the
Hamiltonian and the conserved charges.

\section{Plactic algebras, universal enveloping algebras and crystal graphs}

In the previous section we have largely reviewed well-known algebraic
structures from the physics literature on quantum integrable systems. In
this section we introduce tools from algebraic combinatorics and Lie theory and show
their relations with the phase model. These connections are novel results and are crucial for our setup.

\subsection{The (local affine) plactic algebra}

We start by defining an algebra $\PL(\mathcal{A})$, motivated by the plactic
algebra introduced by Lascoux and Sch\"{u}tzenberger in \cite{LasSchutz} with the following natural quotient studied for
instance in \cite{FG}:

\begin{definition}
The \emph{local plactic algebra} $\PL_{\text{fin}}$ is the free algebra generated by the elements of $\{a_1,a_2,\ldots, a_{n-1}\}$ modulo the relations
\begin{eqnarray}
a_{i}a_{j}=a_{j}a_{i}, &&\text{ if  }|i-j|>1,\label{PLfin1}\\
a_{i+1}a_{i}^{2}=a_{i}a_{i+1}a_{i},&& a_{i+1}^{2}a_{i}=a_{i+1}a_{i}a_{i+1},\label{PLfin2}
\end{eqnarray}
whenever the expressions are defined.
\end{definition}

Recall from \cite{LasSchutz} the plactic monoid $M$ defined by the Knuth relations. The Robinson-Schenstedt algorithm gives a bijection between the equivalence
classes of the monoid $M$ and the set $\mathcal{T}$ of tableaux with filling from $\{1,2,\ldots, n-1\}$ \cite{FultonYT}. Given a tableau $\cT$ we obtain the corresponding word by reading the columns from left to right and bottom to top, replacing each number $i$ by the generator $a_i$. These words form then in particular a basis of the monoid algebra $\mC[M]$. The {\it local} plactic algebra is the quotient of $\mC[M]$ with the additional local relations \eqref{PLfin1}. In the following we summarize a few properties of the local plactic algebra:

\begin{proposition}[Local plactic algebra]\hfill
\label{finiteplactic}
\begin{enumerate}
\item \label{1}
Let $N=\frac{(n-1)n}{2}$, then the words of the form\hfill\\
$a_1^{\alpha_1}(a_2a_1)^{\alpha_2}(a_2)^{\alpha_3}(a_3a_2a_1)^{\alpha_4}(a_3a_2)^{\alpha_5}(a_3)^{\alpha_6}\ldots
(a_{n-1}a_{n-2})^{\alpha_{N-1}}(a_{n-1})^{\alpha_N}$\hfill\\
where the $\alpha$'s run through
the nonnegative integers form a basis $B$ of $\PL_{\text{fin}}$.
\item \label{2}
The Robinson-Schenstedt correspondence defines a bijection between $B$ and the set
$\mathcal{T}_{\text str}$ of tableaux in $\mathcal{T}$, where the number of $j$'s appearing in row $i$ is smaller or
equal the number of $(j-1)$'s appearing in row $i-1$ for any $i>1$, $j\in\{1,2,\ldots, n-1\}$.
\item \label{3}
There is a 2-parameter deformation $U^+_{r,s}$ (with generic $r$ and $s$) of the universal enveloping algebra $U^+$ of the strictly upper
triangular matrices of $\mathfrak{sl}(n)$ such that the specialisation $U^+_{q,q^{-1}}$ is
isomorphic to the usual Drinfeld-Jimbo quantisation $U_q^+$ of $U^+$ (with generic $q$) and $U_{0,1}$ is isomorphic to
$\PL_{\text{fin}}$. The basis above is a specialisation of both, a canonical basis and a PBW-basis.
\end{enumerate}
\end{proposition}
\begin{proof} For \eqref{2} we indicate the bijection for the case $n=4$. The general case is completely
analogous. Note that the Robinson-Schenstedt algorithm transfers a word $a_ia_{i-1}a_{i_3}\ldots a_j$ with $i\geq j$  into a  tableau consisting of a single column with entries $j,j+1,\ldots,i-1,i$ (starting from the top). Then more generally, the word
$a_1^{\alpha_1}(a_2a_1)^{\alpha_2}(a_2)^{\alpha_3}(a_3a_2a_1)^{\alpha_4}(a_3a_2)^{\alpha_5}(a_3)^{\alpha_6}$ is mapped
under the Robinson-Schensted algorithm to the tableau where the first row contains $\alpha_1+\alpha_2+\alpha_4$ ones,
$\alpha_3+\alpha_5$ twos and $\alpha_6$ threes, the second row contains $\alpha_2+\alpha_4$ twos and $\alpha_5$ threes,
the last row contains $\alpha_4$ threes. In particular, the tableau is in $\mathcal{T}_{\text str}$ and the map is injective on $B$.
On the other hand varying the $\alpha's$ provides all possible tableaux in
$\mathcal{T}_{\text str}$. Hence \eqref{2} follows.
Part \eqref{1} follows from part \eqref{3} or from the explicit algorithm in the next subsection.
Part \eqref{3} goes back to \cite{Takeuchi} where Takeuchi defined a $\mathbb{C}[r,s]$-algebra $U^+_{r,s}$ with generators $E_i$, $1\leq i\leq n$ and relations $E_iE_j$ if $|i-j|>1$ and
\small
\begin{eqnarray*}
E_{i+1}^2E_i-(r+s)E_{i+1}E_iE_{i+1}+rs E_iE_{i+1}^2,&&E_{i+1}E_i^2-(r+s)E_{i}E_{i+1}E_{i}+rs E_i^2E_{i+1}.
\end{eqnarray*}
\normalsize
Setting $r=0$ and $s=1$ gives exactly the relations \eqref{PLfin1} and \eqref{PLfin2}. Setting $r=q$ and $s=q^{-1}$ gives the well-known quantum Serre relations (see e.g.
\cite[4.3,4.6,4.12b]{Jantzen}).

Using standard arguments (see e.g. \cite[\S8]{Jantzen}) one can show that this algebra is a free  $\mathbb{C}[r,s]$-module,
with a basis given by the elements
$$E_1^{\alpha_1}(E_2E_1)^{\alpha_2}(E_2)^{\alpha_3}(E_3E_2E_1)^{\alpha_4}(E_3E_2)^{\alpha_5}\ldots
(E_{n-1}E_{n-2})^{\alpha_{\frac{(n-1)n-2}{2}}}(E_n)^{\alpha_\frac{(n-1)n}{2}},$$ that is the PBW basis associated with
the special reduced expression $$w_0=s_1(s_2s_1)s_2(s_3s_2s_1)(s_3s_2)\ldots$$ of the longest element of Weyl group
(i.e. the symmetric group) $S_n$ of $\mathfrak{sl}(n)$ (here $s_i$ denotes the elementary transposition $(i,i+1)$).
The existence of a canonical basis for $U^+_{r,s}$ is proved in \cite{Reineke} with the property that it specialises to the basis in \eqref{1}.
\end{proof}
\begin{remark}{\rm
\begin{itemize}
\item The quantised universal enveloping algebra $U_q^+$ is isomorphic to (a twisted version of) {\it Ringel's Hall algebra} \cite{Ringel}. In this context, the plactic algebra appears as the Hall algebra with multiplication defined using generic extensions, hence a {\it generic Hall algebra} \cite{Reineke}. Reineke also defines the Hall algebra version of $U^+_{r,s}$.
\item In Section \ref{Appendix} we give an explicit algorithm which turns an arbitrary tableau into a tableau in $\mathcal{T}_{\text str}$ using the relations of the local plactic algebra.
\end{itemize}
}
\end{remark}

\begin{definition}{\rm
Let $\cA=\{a_0,a_1,a_2,\ldots a_{n-1}\}$. The \emph{\ affine local plactic algebra} $\PL=\PL(
\mathcal{A})$ is the free algebra generated by the elements of $\cA$ modulo the relations
\begin{eqnarray}
a_{i}a_{j}-a_{j}a_{i}=0,&&\text{ if $|i-j|\neq 1\mod n$},\label{PL1}\\
a_{i+1}a_{i}^{2}=a_{i}a_{i+1}a_{i},&& a_{i+1}^{2}a_{i}=a_{i+1}a_{i}a_{i+1},\label{PL2}
\end{eqnarray}
where in \eqref{PL2} all variables are understood as elements in $\cA$ by taking indices modulo $n$.}
\end{definition}
\begin{example}
 {\rm If $n=3$ then the defining relations are $a_2a_1^2=a_1a_2a_1$, $a_0a_2^2=a_2a_0a_2$,
$a_2^2a_1=a_2a_1a_2$, $a_0^2a_2=a_0a_2a_0$. (Note that $a_2$ and $a_0$ do not commute.)}
\end{example}

\subsection{The generalised Robinson-Schensted algorithm}
An {\it $n$-multi-partition} is an $n$-tuple $\pi=(\pi^{(0)},\ldots, \pi^{(n-1)})$ of partitions (resp. Young diagrams). A
multi-partition is {\it  aperiodic} if for any positive number $l$ there is at least one Young diagram which does not
have a column of height $l$. An {\it aperiodic multi-tableau} or just a multi-tableau is given by taking an aperiodic
multi-partition $\pi$ and putting the number $i+r-1$ (modulo $n$) into the boxes in the $r^{\text{th}}$ row of $\pi^{(i)}$.

\begin{figure}
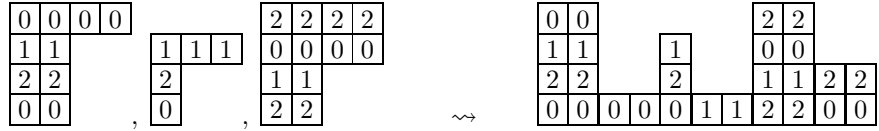

\caption{An example of an aperiodic multipartition, $((4,2,2,2),(3,1,1,0),(4,4,2,2))$. Displayed are its multi-tableau and the multi-tableau with all columns pushed down which then  corresponds to the word $a_0^2a_2^4a_1^5a_0^6a_2^3a_1^2a_0^3a_2^2$.}
\begin{equation*}
{{\;\young(0000,11,22,00),\;\young(111,2,0),\;\young(2222,0000,11,22)\;}}\quad\quad\rightsquigarrow\quad\quad
{{\young(00,11,22,00)\young(00)\young(1,2,0)\young(11)\young(22,00,11,22)\young(2,0)\young(2,0)}}
\end{equation*}
\label{fig:towers}
\end{figure}

Given a word in the local affine monoid, Deng and Du \cite{DengDu} assign (following ideas of Lusztig and Ringel) a multi-partition by the following
algorithm: given a multi-partition $\pi=(\pi^{(0)},\ldots, \pi^{(n-1)})$ and $i\in\{0,\ldots, n-1\}$ then $a_i\pi$ is the
multi-partition obtained from $\pi$ by adding to $\pi^{(i)}$ an extra box in the first row if $\pi^{(i+1)}=\emptyset$, and otherwise remove the first column, $c$, of $\pi^{(i+1)}$ and place a column of length one box longer than $c$ in the partition
$\pi^{(i)}$. (The indices are again taken modulo $n$). The multi-partition $\pi$ associated to a word $w\in\PL$ is
then defined as $\pi=w(\emptyset, \ldots, \emptyset)$. One can easily see that this multi-partition is aperiodic.
There is also an explicit way to read off a word from a multi-partition as follows: first we convert the multi-partition into a multi-tableau (there is a unique way to do this), then we draw the multi-partitions as a diagram as indicated in
Figure \ref{fig:towers} by pushing down all the columns to a common baseline. Considering only the numbers on top of a column, remove all those $0$'s appearing higher than all $1$'s. Considering again only the numbers on top of a column, remove all those $1$'s appearing higher than all $2$'s etc. Repeat this procedure (with the cyclic ordering) as long as possible. The aperiodicity guaranties that this process can be carried on until there are no boxes left. The order in which the numbers got removed defines a word $w(\pi)$ in the local affine plactic algebra. We call these
words {\it standard words}.

\begin{proposition}\cite[Theorem 4.1 and its proof]{DengDu}
\label{Standardbasis}
\begin{enumerate}
\item The two algorithms are inverse to each other, i.e. $\pi(w(\pi))=\pi$ and $w(\pi(w))$ is equivalent to $w$ in
    the local plactic monoid.
\item The standard words associated to aperiodic multi-partitions (resp. multi-tableaux) form a basis of $\PL$.
\item If we replace the letters in the standard words by the usual Chevalley generators $E_i$, $0\leq i\leq n-1$ of
    the positive part $U^+$ of the universal enveloping algebra $U(\widehat{\mathfrak{sl}}(n))$ then we obtain a monomial
    basis of $U^+$.
\end{enumerate}
\end{proposition}

\begin{remark}{\rm
For a 2-parameter deformation of $U^+$ and a monomial basis in terms of Lyndon words we refer to \cite{HuRossoZhang}.}
\end{remark}

\begin{proposition}[Faithfulness]
\label{faithfulness}
There is a homomorphism of algebras $\Psi_{\text{fin}}:\PL_{\op{fin}}\rightarrow\Phi$ such that
\begin{equation*}
a_{j}\mapsto \varphi _{j+1}^{\ast }\varphi _{j},\quad j=1,...,n-1
\end{equation*}
In particular, the representation \eqref{phase1} of the phase algebra $\Phi$ lifts to a representation of the local
plactic algebra $\PL_{\text fin}$. This representation is faithful. Moreover, it lifts to a representation of $\PL$ on $\cH[z]$ by mapping $a_{0}=a_n$ to $z\varphi
_{1}^{\ast }\varphi _{0}$ and the $a_{j}$ as above. This representation is again faithful.
\end{proposition}

Explicitly, the action on $\cH_k$ in terms of the basis vectors given by affine weights reads for $i>0$%
\begin{equation*}
a_{i}:P^+_{k}\rightarrow P^+_{k}\cup\{0\},\qquad \hat{\lambda}
\mapsto
\begin{cases}
\hat{\lambda}+\hat{\omega}_{i+1}-\hat{\omega}_{i} &
\text{if $\hat{\lambda}+\hat{\omega}_{i+1}-\hat{\omega}_{i}\in P_k^+$}\\
\\
0 & \text{otherwise.}
\end{cases}
\end{equation*}%
In terms of Young diagrams, the endomorphism $a_i$ adds a box in the $(i+1)^{\text{th}}$ row of the Young diagram $\op{\hat{P}}(\hat\la)$ (see Lemma \ref{weights}) provided the result is again a Young diagram. For $i=n$ it removes an $n$-column (if there is none, the result is zero), adds a box in the first row and multiplies with $z$. Note that these actions physically correspond to moving single particles by one site in clockwise direction on the affine Dynkin diagram (see Figure ~\ref{fig:bosons}).

\begin{remark}
\label{crystal}
{\rm
The set $P_k^+$ with the action of the local (affine) plactic generators $a_i$ (and the obvious weight function) has the structure of an abstract crystal. This crystal is known (see \cite{JMMO}) to be isomorphic to the {\it crystal of the $k^{\text{th}}$ symmetric tensor representation of the vector representation} of the corresponding quantized universal enveloping algebra of adjoint type (in the sense of \cite[4.5]{Jantzen}).
}
\end{remark}
\begin{proof}[Proof of Proposition \ref{faithfulness}]
The existence of this morphism follows directly from the definitions. That the representation $\cH$ is faithful follows
from Proposition~\ref{finiteplactic} and Remark~\ref{crystal}, but we give an explicit argument here. Let $X=\sum c_b b$
be a finite linear combination of basis elements $b\in B$ (see Proposition \ref{finiteplactic}) in the plactic algebra $\Pfin$. Assume $X$ acts by zero on $\cH$.
Applying this to the partition containing only one box implies that $c_b=0$ for all basis vectors $b$ which only consist of a single monomial of the form $(a_ra_{r-1}\ldots a_1)$ for arbitrary $r$. More generally, applying $X$ to a column of height $i$, the partition $(1^i)$, then implies that $c_b=0$ for all basis vectors consisting of one monomial of the form $(a_ra_{r-1}\ldots a_{i+1}a_i)$. In order to single out basis words consisting of two monomials consider next the partition $(2,1)$. Because of the special structure of the basis vectors application of $X$ then implies that $c_b=0$ for all summands $b$ which now consist of two monomials, one of the form $(a_ra_{r-1}\ldots a_1)$ for arbitrary $r$ and the other of the form $(a_sa_{s-1}\ldots a_2)$ for arbitrary $s$. Continuing with the partition $(3,1)$ implies that $c_b=0$ for all basis vectors $b$ consisting of two monomials of the form $(a_ra_{r-1}\ldots a_1)$ for arbitrary $r$ and one factor of the form $(a_sa_{s-1}\ldots a_2)$ etc. %
Carrying on like this gives the faithfulness. The faithfulness of the second representation follows then from Proposition~\ref{Standardbasis} and the definition of the affine Lie algebra $\widehat{\mathfrak{sl}}(n)$.

\end{proof}

Henceforth, we shall always identify the local affine plactic algebra $\PL$ with its image in $\limfunc{End}\cH[z]$.
\subsection{Noncommutative polynomials}\label{ncpoly1}
Mimicking the case of the ordinary local plactic algebra we now
introduce noncommutative polynomials in the generators $\{a_{i}\}$ of $\PL$ which are (noncommutative, affine) analogs
of the ordinary elementary and complete symmetric functions. (This approach is similar to \cite{FG}, but unfortunately
does not satisfy their assumption, so that we cannot use their results directly.)
We need the notion of cyclically ordered products $\tprod_{i\in I}^{\circlearrowright
}a_{i}$ and $\tprod_{i\in I}^{\circlearrowleft
}a_{j} $. A monomial $a_{i_1}a_{i_2}\cdots a_{i_r}$ in the variables $\cA$ is {\it clockwise cyclically ordered}
(respectively {\it anticlockwise cyclically ordered}) if for any two indices $i_j$, $i_k$ with $i_k=i_j+1$ modulo $n$,
the variable $a_{i_j}$ occurs to the left (to the right) of $a_{i_k}$. (In case $i_k\neq i_j+1$ the order does not matter because of \eqref{PL1}.)  The origin of the name becomes obvious if we identify $a_i$ with the corresponding point on the Dynkin diagram $\Gamma$: there are two circle segments connecting the two
points. If they are not of the same length we choose the shorter one and the (anti-)clockwise order is the same as the intuitively defined anti-clockwise order with respect
to this segment. For any monomial $a_{i_1}a_{i_2}\cdots
a_{i_r}$ not containing all the generators $a_i$, there is a unique clockwise (resp. anti-clockwise) cyclically monomial which differs only by a permutation of the variables. We denote
it by $\displaystyle{\tprod_{1\leq j\leq r}^{\circlearrowright }a_{i_j}}$, respectively $\displaystyle{\tprod_{1\leq j\leq r}^{\circlearrowleft }a_{i_j}}$.

\begin{definition}\label{edef}{\rm
For $1\leq r\leq n-1$ we define {\em cyclic noncommutative elementary symmetric polynomials} as the following elements of $\PL$
\begin{equation}\label{ncehdef}
e_{r}(\cA)=\sum_{|I|=r}\tprod_{i\in I}^{\circlearrowleft }a_{i}
\end{equation}%
where the sum runs over
all sets $I=\{i_{1},\ldots,i_{r}\}$ with $i_{s}\neq i_{t}$ for $s\neq t$.
}
\end{definition}

\begin{example}
If $n=4$ then $$e_{2}(\cA)=a_{2}a_{1}+a_{3}a_{1}+a_{1}a_{0}+a_{3}a_{2}+a_{0}a_{2}+a_{0}a_{3}\;.$$
\end{example}

\begin{remark}{\rm
To define the noncommutative elementary symmetric polynomials we had to make a choice for the cyclic order. The other choice would give just the adjoint operators with respect to  scalar product $\langle\;,\rangle$ on $\cH$ from Proposition \ref{phasealgebra}.
}
\end{remark}

 The following result realizes the transfer matrix as the generating function of the elementary symmetric polynomials:

\begin{proposition}[Generating function]
\label{genfunction}
Let $T(u)$ denote the transfer matrix of the phase model as before. Then
\begin{equation}
T(u)=A(u)+zD(u)=\sum_{r=0}^{n}e_{r}(\mathcal{A})u^{r},
\label{T is e}
\end{equation}%
where we define $e_{0}(\mathcal{A})=1$ and $e_{n}(\mathcal{A}%
)=z~1.$
\end{proposition}

\begin{proof}
The proof uses the explicit form of the transfer matrix in terms of the
phase algebra stated in Lemma \ref{ABCD}, $T_{r}=A_{r}+zD_{r}$ with%
\begin{equation*}
A_r= \sum_{0\leq s\leq \min\{r,n/2\}}\sum_{\substack{ i_{1}<j_{1}<\cdots
<i_{s}<j_{s}  \\ (j_{1}-i_{1})+\cdots +(j_{s}-i_{s})=r}}a_{i_{1},j_{1}}%
\cdots a_{i_{s},j_{s}}
\end{equation*}
and
\begin{equation*}
D_r=\sum_{0\leq s\leq \min\{n/2,n-r\}}\sum_{\substack{ i_{1}<j_{1}<\cdots
<i_{s}<j_{s}  \\ (j_{1}-i_{1})+\cdots +(j_{s}-i_{s})=n-r}}%
a_{j_{1},i_{1}}\cdots a_{j_{s},i_{s}}\; .
\end{equation*}
For $j>i$, the identities
$a_{i,j}=a_{j-1}a_{j-2}\cdots a_{i}$ and
$a_{j,i}=z^{-1}a_{i-1}\cdots a_{1}a_{0}a_{n-1}\cdots a_{j}$ are easily verified.
The asserted equality between $T_{r}$ and $e_{r}(\mathcal{A})$
follows.
\end{proof}

\begin{corollary}
\label{ecomm}
The cyclic noncommutative elementary symmetric functions pairwise commute,
\begin{equation}
\lbrack e_{r}(\mathcal{A}),e_{r^{\prime }}(\mathcal{A})]=0,\qquad
0\leq r,r^{\prime }\leq n\;.  \label{e's commute}
\end{equation}
and generate a commutative subalgebra $\cA^{\op{sym}}$ of $\cA$.
\end{corollary}

\begin{proof}
The first part is a direct consequence of \eqref{T is e} and the
integrability \eqref{int}. The second part then follows from Proposition \ref{faithfulness}.
\end{proof}

Thanks to Corollary \ref{ecomm}, the following definition makes sense\footnote{Note also the stronger (independent) result of Theorem \ref{eigenbasis}, namely the construction of a simultaneous eigenbasis.}:

\begin{definition}
{\rm Given a partition $\lambda $, the {\it Jacobi-Trudy formula}
\begin{equation}
s_{\lambda }(\mathcal{A})=\det \left(e_{\lambda _{i}^{t}-i+j}(\mathcal{A})\right).
\label{JTformulae}
\end{equation}
gives well-defined polynomials (in the generators of the affine plactic
algebra), which we call the \textit{cyclic noncommutative Schur polynomials}%
.}
\end{definition}
In further analogy with the commutative case we also introduce noncommutative versions of the complete symmetric polynomials.

\begin{definition}{\rm
Define the set of {\em cyclic noncommutative complete symmetric polynomials} via
\begin{equation}\label{nchdef}
h_{r}(\mathcal{A})=\det (e_{1-i+j}(\mathcal{A}))_{1\leq i,j\leq r}\;.
\end{equation}
In particular, we have for $\lambda=(r)$, a horizontal $r$-strip, that $s_{\la}(\cA)=h_r(\cA)$ analogous to the commutative case.
}
\end{definition}

\section{Bethe Ansatz and the isomorphism of rings}

Employing the algebraic Bethe Ansatz (or Quantum Inverse Scattering Method) one can compute the
eigenvectors of the transfer matrix (\ref{tcoeff}), which according to (\ref{T is e}) generates the cyclic noncommutative elementary symmetric
functions. Algebraic Bethe Ansatz calculations for the phase model were first undertaken by Bogoliubov et al. \cite{Bogoliubovetal}. We shall give here a far more detailed analysis of the solutions to the Bethe Ansatz equations and generalize their discussion from periodic (z=1) to quasi-periodic (z generic) boundary conditions. The important novel aspect of our work with regard to the discussion of the phase model in \cite{Bogoliubovetal} is the connection with representation theory and finally the identification of the abelian algebra generated by the commuting transfer matrices of this integrable system with the Verlinde algebra (by showing that the so-called Bethe vectors are given in terms of Weyl characters).

\subsection{Bethe vectors}
\label{sec:Bethevectors}
We introduce the following notation for the special vector already discussed
in \eqref{BetheSchur}
\begin{equation}
b(x):=B(x_{1}^{-1})\cdots B(x_{k}^{-1})\emptyset =\sum_{\hat{\lambda}\in
P_{k}^{+}}s_{\hat{\lambda}^{t}}(x_{1}^{-1},...,x_{k}^{-1})\hat{\lambda}\;.
\label{bethevec}
\end{equation}

The vectors $b(x)$ are by definition in $\cH_k$. We first like to find precise conditions on $x$, ensuring that $b(x)$ is an eigenvector of the transfer matrix $T(u)$, and hence a simultaneous eigenvector of the noncommutative elementary symmetric functions $e_r$ and then also for the whole algebra of noncommutative symmetric functions. We will show that the Bethe vectors form an orthogonal basis of $\cH_k$.
The Ansatz  for the algebraic form of the eigenvectors (the first identity in \eqref{BetheSchur}) and the derivation of the resulting conditions on $x$, usually called {\it Bethe Ansatz equations}, exploiting the quantum Yang-Baxter algebra, is standard in the physics literature, (see for example \cite{KBI}). The following proposition states the result of the algebraic Bethe Ansatz. For convenience we sketch its proof:

\begin{proposition}
\label{eigenvalues}
\begin{enumerate}
\item For the vector $b(x)$ depending on invertible indeterminates $x=(x_{1},\ldots,x_{k})$ to be an eigenvector for the transfer matrix $T(u)$, $x$ must obey the Bethe Ansatz equations
\begin{equation}
x_{1}^{n+k}=\cdots =x_{k}^{n+k}=(-1)^{k-1}ze_{k}(x_{1},\ldots,x_{k}),
\label{bae}
\end{equation}%
where $e_{k}(x)=\prod_{i=1}^{k}x_{i}$ is the (ordinary) $k^{\text{th}}$ elementary
symmetric polynomial.
\item Given a solution $x=(x_{1},\ldots,x_{k})$ of the Bethe Ansatz equations (\ref%
{bae}), $b(x)$ is an eigenvector of the transfer matrix (\ref{tcoeff}) and
of the integrals of motion $H_{r}^{\pm }$ (\ref{imotion}) for the phase model. The eigenvalues are given by
\begin{equation}
T(u)b(x)=\left[ 1+(-1)^{k}ze_{k}(x)u^{n+k}\right] \prod_{i=1}^{k}\frac{1}{%
1-ux_{i}}~b(x)  \label{Tspec}
\end{equation}%
and%
\begin{equation}
H_{r}^{\pm }b(x)=-\frac{h_{r}(x_{1},\ldots,x_{k})\pm
z h_{r}(x_{1}^{-1},\ldots,x_{k}^{-1})}{2}~b(x)\ ,  \label{Hspec}
\end{equation}%
where $e_{r}$, $h_{r}$ is the $r^{th}$ ordinary (commutative) elementary and
complete symmetric polynomials, respectively.
\end{enumerate}
\end{proposition}

\begin{proof}[Sketch of proof:]
It is easy to verify $A(u)\emptyset=\emptyset$ and $D(u)\emptyset
=u^{n}\emptyset$. One can show via induction on $k$ that
\begin{eqnarray*}
T(u)B(x_{1}^{-1})...B(x_{k}^{-1})\emptyset &=&\left( \prod_{i=1}^{k}\frac{1}{%
1-ux_{i}}+zu^{n}\prod_{i=1}^{k}\frac{ux_{i}}{ux_{i}-1}\right)
B(x_{1}^{-1})...B(x_{k}^{-1})\emptyset \\
&&+%
\sum_{i=1}^{k}F_{i}(u)~B(x_{1}^{-1})...B(x_{i-1}^{-1})B(x_{i+1}^{-1})...B(x_{k}^{-1})\emptyset ,
\end{eqnarray*}%
where $F_{i}$ is given by
\begin{eqnarray*}
F_{i}(u) &=&\frac{1}{1-ux_{i}}\left( \prod_{j\not=i}\frac{x_{i}}{x_{i}-x_{j}}%
-zx_{i}^{-n}\prod_{j\not=i}\frac{x_{j}}{x_{j}-x_{i}}\right) \\
&=&\frac{x_{i}^{k-1}}{1-ux_{i}}\left( 1+z(-1)^{k}x_{i}^{-n-k}e_{k}(x)\right)
\prod_{j\not=i}\frac{1}{x_{i}-x_{j}}
\end{eqnarray*}%
(The case $k=1$ is immediate from the commutation relation \eqref{aba1}). The second sum is referred
to as the {\it unwanted terms} which have to vanish in order to turn the Bethe
vector into an eigenstate of the transfer matrix. This leads to the condition%
\begin{equation*}
1+z(-1)^{k}x_{i}^{-n-k}e_{k}(x)=0\;,\qquad i=1,2,\ldots,k,
\end{equation*}%
which are the Bethe Ansatz equations. Hence $b(x)$ is an eigenvector if and only if $x$ satisfies the Bethe Ansatz equations. The eigenvalues can now be simply deduced by a series expansion in the spectral parameter $u$ and using Lemma \ref{hformula} below.
\end{proof}

\begin{lemma}
\label{hformula}
Assume  $x:=(x_1,x_2,\ldots x_n)$ satisfies the Bethe equations \eqref{bae}. Then
\begin{eqnarray}
\label{hformel}
h_{n-r}(x_1,x_2,\ldots x_n)=zh_{r}(x_{1}^{-1},\ldots,x_{k}^{-1})
\end{eqnarray}
for any $0\leq r\leq n$.
\end{lemma}
\begin{proof}
We use
$ze_k(x)=(-1)^{k-1}x_i^{n+k}$ and calculate
\begin{eqnarray*}
\begin{array}[t]{ll}
h_{n-r}(x_{1},\ldots,x_{k})
=\sum_{i=1}^{k}x_{i}^{n-r}\prod_{j,j\neq i}\frac{x_{i}}{x_{i}-x_{j}}
=\sum_{i=1}^{k}x_{i}^{n-r+(k-1)}\prod_{j,j\neq i}\frac{1}{x_{i}-x_{j}}\\
=z\sum_{i=1}^{k}x_{i}^{-r-1}e_{k}(x)\prod_{j,j\neq i}\frac{1}{x_{j}-x_{i}}
=z\sum_{i=1}^{k}x_{i}^{-r-1}e_{k}(x)\prod_{j,j\neq i}\frac{x_i^{-1}x_j^{-1}}{x_{i}^{-1}-x_{j}^{-1}}\\
=z\sum_{i=1}^{k}x_{i}^{-r}\prod_{j,j\neq i}\frac{x_i^{-1}}{x_{i}^{-1}-x_{j}^{-1}}
=zh_{r}(x_{1}^{-1},\ldots,x_{k}^{-1}),
\end{array}
\end{eqnarray*}
where only the first and last equality need some explanation. Following \cite[page 209]{MacDonald} we define $g_{s}(x_{1},\ldots,x_{k};q)$ through the formula
\begin{equation}
\label{HallLittel}
\prod_{i}\frac{1-uq~x_{i}}{1-u~x_{i}}=\sum_{s\geq
0}g_{s}(x_{1},\ldots,x_{k};q)u^{s}
\end{equation}
In particular, $g_0=1$. In case $s>0$, $g_{s}(x_{1},\ldots,x_{k};q)=(1-q)P_{(s)}(x_{1},\ldots,x_{k};q)$ where $P_{(s)}(x_{1},\ldots,x_{k};q)$ is the so-called {\it Hall-Littelwood polynomial}, explicitly
\begin{equation*}
g_{s}(x_{1},\ldots,x_{k};q)=(1-q)\sum_{i=1}^{k}x_{i}^{s}\prod_{j\neq i}\frac{%
x_{i}-qx_{j}}{x_{i}-x_{j}}.
\end{equation*}%
Our two formulae in question follow then from the specialisation $q\mapsto 0$.
\end{proof}

The Bethe Ansatz equations \eqref{bae} can be reformulated as follows

\begin{lemma}
\label{baeequiv}
The equations (\ref{bae}) are equivalent to the conditions%
\begin{equation}
h_{n+1}(x)=\cdots =h_{n+k-1}(x)=h_{n+k}(x)+(-1)^{k}ze_{k}(x)=0,\quad h_n(x)=z \label{bae2}
\end{equation}
\end{lemma}

\begin{proof}
Assume the equations \eqref{bae} hold. Integrability, property \eqref{int}, implies that the transfer matrix $T(u)$ has a simultaneous (i.e. independent of $u$) eigenspace decomposition. The definition of the transfer matrix \eqref{tcoeff} implies that the eigenvalues in \eqref{Tspec} must be polynomial in the spectral parameter $u$ and at most of degree $n$. The eigenvalues in \eqref{Tspec} are of the form $(1+(-1)^{k}ze_{k}(x)u^{n+k})H(u)$, where $H(u)$ is the generating function of the complete symmetric polynomials. The coefficients of $u^j$ for $j>n$ have to vanish, hence $0=h_{n+j}(x)$ for $1\leq j\leq k-1$ and $h_{n+k+j}(x)+(-1)^{k}ze_{k}(x)h_j(x)=0$ for $j\geq 0$ which of course implies \eqref{bae2}. The last identity in \eqref{bae2} follows from Lemma \ref{hformula} for $r=0$.

It is in fact equivalent, because $h_{n+k+j}(x)+(-1)^{k}ze_{k}(x)h_j(x)=0$ implies
\small
$$0=\left(\sum_{i=1}^k x_i\right) \left(h_{n+k+j}(x)+(-1)^{k}ze_{k}(x)h_j(x)\right)=h_{n+k+j+1}(x)+(-1)^{k}ze_{k}(x)h_{j+1}(x).$$
\normalsize
Conversely if we assume
\eqref{bae2} holds, then \eqref{Tspec} is a polynomial in $u$ (of degree $n$). If we choose a generic point $x$, fix $x_i$ and formally replace $u$ by $u-x_i^{-1}$, then \eqref{Tspec} turns into $p(u)q(u)ux_i^{-1}$, where $p(u)$ is a polynomial, and $q(u)$ is a formal power series in $u$. Since this has to be a polynomial in $u$ again, its residue, which is a multiple of $1+(-1)^{k}ze_{k}(x)x_i^{-(n+k)}$, must be zero and the equations \eqref{bae} hold.
\end{proof}

To state the explicit form of the Bethe roots, i.e. solutions to \eqref{bae}, we formally want to allow $n^{\text{th}}$ roots of $z$ and also inverses, hence work in the ring $R=\mC[z^{\mp{1}}, z^{\mp\frac{1}{n}}]=\mC[z^{\mp\frac{1}{n}}]$. We extend the complex conjugation to $R$ by setting $\overline z:=z^{-1}$. The expansion \eqref{bethevec} into Schur polynomials shows that the Bethe vectors $b(x)$  do not really depend on the tuple $x=(x_1,x_2,\ldots , x_k)$, but only on the set $\{x_1,x_2,\ldots , x_k\}$ or its $S_k$-orbit. Furthermore, in case $x_i=x_j$ for some $i\not=j$ we get $b(x)=0$. Therefore we denote $$\op{Sol}(n,k):=\{x\in R^k\mid \text{$x$ solves \eqref{bae} and the $x_i$ are pairwise distinct}\}$$ and want to describe the orbit space $\widehat{\op{Sol}(n,k)}:=\op{Sol}(n,k)/S_k$.
\begin{theorem}[Solutions of the Bethe Ansatz]\label{Solutions}
There is a bijection
\begin{eqnarray}
\mathfrak{P}_{\leq n-1,k}&\cong&\widehat{\func{Sol}(n,k)},\nonumber\\
\sigma&\mapsto&x_\sigma:=z^{\frac{1}{n}}\zeta ^{\frac{|\sigma|}{n}} \left(\zeta ^{I_{1}},\ldots,\zeta^{I_{k}}\right),
\label{Betheroots}
\end{eqnarray}
where $\zeta =\exp \frac{2\pi \iota}{k+n}$ and
\begin{equation}
I=I(\sigma ^{t}):=\left( \tfrac{k+1}{2}+\sigma
_{k}^{t}-k,\ldots ,\tfrac{k+1}{2}+\sigma _{1}^{t}-1\right).  \label{Imap0}
\end{equation}
Moreover, $\widehat{\op{Sol}(n,k)}$ decomposes into the disjoint union of orbits under the $\mathbb{Z}_n$-action given by $\ell.x:=e^\frac{2\pi \iota\ell}{n}x$ with $\ell\in\mathbb{Z}$, $x\in \widehat{\op{Sol}(n,k)}$.
\end{theorem}

\begin{proof}
It is easy to check that $\sigma\mapsto I(\sigma ^{t})$ defines a bijection
between $\mathfrak{P}_{\leq n-1,k}$ and the following set of tuples of
(half)-integers
\begin{equation*}
\mathcal{I}_{n-1,k}:=\left\{(I_{1},...,I_{k})\mid -\tfrac{k-1}{2}\leq
I_{1}<\cdots <I_{k}\leq n-1+\tfrac{k-1}{2},I_{j}\in \tfrac{k+1}{2}+\mathbb{Z}%
,\forall j\right\} .
\end{equation*}

For instance, $\mathfrak{P}_{\leq 2,3}$ from \eqref{PYoung} gets identified
with the tuples
\begin{equation*}
\left\{(-1,0),(-1,1),(-1,{2}),(-1,3),(0,1),(0,2),(0,3),(1,2),(1,3),(2,3)%
\right\}.
\end{equation*}
Hence it is enough to show that the $x_\sigma$ as in \eqref{Betheroots} with
$I\in\mathcal{I}_{n-1,k}$ are a full list of representatives for $\widehat{%
\func{Sol}(n,k)}$.

Assume $x\in R^k$ solves \eqref{bae}. First we describe the dependence on $z$%
. Let $s_i$ be the highest exponent of $z$ in $x_i$. Then the highest
exponent in $x_i^{n+k}$ is $s_i(n+k)$. In particular, $s_i=s_j$ for any $%
1\leq i,j\leq n$ and then also $s_i(n+k)=1+ks_i$, hence $s_in=1$. The same
argument also applies to the lowest exponent, hence any $x_i$ is a monomial
in $z^{\frac{1}{n}}$. In the following it will be enough to consider
therefore the special case $z=1$.

First note that \eqref{bae} implies that all $x_j$ are of norm $1$, hence of
the form $x_{j}=e^{2\pi \iota\alpha _{j}}$ for some $\alpha_j\in\mathbb{R}$.
Then the Bethe Ansatz equations (\ref{bae}) for the phase model can be
rewritten (by taking the logarithm) as
\begin{equation}
\alpha _{j}=\frac{J_{j}}{k+n-1}+\frac{1}{k+n-1}\sum_{l\neq j}\alpha _{l}\;
\label{log bae}
\end{equation}%
with the $J_{j}$ being half-integers in case $k$ is even, and integers in
case $k$ is odd. For each fixed generic configuration $J=(J_{1},\ldots,J_{k})
$ one easily checks (using our general assumption $n+k>2$) that there is
precisely one solution to this system of linear equations, namely $\alpha
_{j}=\frac{1}{(k+n)}(\frac{n+1}{n}~J_{j}+\frac{1}{n}\sum_{l\neq j}J_{l})=%
\frac{1}{(k+n)}(J_{j}+\frac{1}{n}||J||)$ with $||J||=\sum_{l=1}^kJ_{l})$.
Let $x(J)=z^{\frac{1}{n}}\zeta ^{\frac{||J||}{n}}(\zeta
^{J_{1}},\ldots,\zeta^{J_{k}})$ with $\zeta=e^{\frac{2\pi \iota}{k+n}}$ be
the corresponding element in $\func{Sol}(n,k)$. The choice $%
J_j:=I_j(\sigma^t)$ gives $x(J)=x(I(\sigma^t))=x_\sigma$. Hence we get all
the solutions $x_\sigma$ predicted in the theorem, but have to show there
are not more (up to permutations).

The assumption $x_j\not=x_i$ for $j\not=i$ implies $\alpha_i\notin\alpha_j+%
\mathbb{Z}$ and then $J_i\not\equiv J_j\mod (n+k)$. Therefore, modulo $n+k$,
the (half)-integers $J_i$ are contained in the interval $[-\tfrac{k-1}{2},n+%
\tfrac{k-1}{2}]$ and are distinct. Since $k<n+k$ there is at least one
(half)-integer $m$ in the interval which does not occur amongst the $J_i$'s.
Hence we can find some $l\in\mathbb{Z}$ such that $J^{\prime}=J+(l,l,\ldots,
l)$ and all the (half)-integer $J_i^{\prime}$ are contained in the interval $%
[-\tfrac{k-1}{2},n-1+\tfrac{k-1}{2}]$ modulo $n+k$ (in other words we shift
the 'gap' $m$ to the end of the interval). The two operations (taking $J$
modulo $n+k$ and adding some $(l,l,\ldots, l)$) however only rescale the
solutions, because
\begin{eqnarray*}
x(J_{1},\ldots,J_{j}+k+n,\ldots,J_{k})&=&z^{\frac{1}{n}}\zeta ^{\frac{||J||}{%
n}+\frac{k+n}{n}}(\zeta ^{J_{1}},\ldots,\zeta ^{J_{k}}) =e^{\frac{2\pi \iota%
}{n}}x(J) \\
x(J_{1}+1,\ldots,J_{k}+1) &=&z^{\frac{1}{n}}\zeta ^{\frac{||J||}{n}+\frac{k}{%
n}+1}(\zeta ^{J_{1}},\ldots,\zeta ^{J_{k}})=e^{\frac{2\pi \iota}{n}}x(J)
\end{eqnarray*}
Therefore, it is enough to show that multiplication with $\eta:=e^{\frac{%
2\pi \iota}{n}}$ preserves the set $X:=\{x(I(\sigma^t))\}$, in other words
we have to show that
\begin{equation*}
\eta x=z^{\frac{1}{n}}\zeta ^{\frac{||I||+k}{n}}(\zeta ^{I_{1}+1},\ldots
,\zeta ^{I_{k}+1})=z^{\frac{1}{n}}\zeta ^{\frac{||I^{\prime}||}{n}}(\zeta
^{I^{\prime}_{1}},\ldots,\zeta ^{I^{\prime}_{k}})\in X
\end{equation*}
with $I^{\prime}_j=I_j+1$ for $1\leq j\leq k$. In case $I_{k}<n-1+\frac{k-1}{%
2}$ there is nothing to do, hence assume $I_{k}=n-1+\frac{k-1}{2}$. Now note
that $x(I^{\prime})=\eta x(I^{\prime}_{1},\ldots
,I^{\prime}_{k-1},I^{\prime}_k-(n+k))=\eta
x(I^{\prime}_k-(n+k),I^{\prime}_{1},\ldots ,I^{\prime}_{k-1})=\eta^{-1}\eta
x(I^{\prime}_k-(n+k)+1,I^{\prime}_{1}+1,\ldots ,I^{\prime}_{k-1}+1)$. Since $%
n+k-1>k$, applying these two equations repeatedly, we finally obtain an
element $x=x(I^{\prime})\in X$. Hence the elements in $\widehat{\func{Sol}%
(n,k)}$ are all of the form as claimed. This proves that we have a bijection
of sets $\mathfrak{P}_{\leq n-1,k}\cong\widehat{\func{Sol}(n,k)}$ and that $%
\widehat{\func{Sol}(n,k)}$ decomposes into disjoint sets whose elements are
obtained by rescaling with an $n^{th}$ root of unity.
\end{proof}

Before deducing interesting consequences we collect well-known formulae (see e.g. \cite{Rietsch}, \cite[Proposition 38.2]{Bump}):
\begin{lemma}
\label{Schurprop}Let $\mu=\hat{\la} \in \mathfrak{P}_{\leq n,k},\sigma \in \mathfrak{P}_{\leq n-1,k}$. With the
notation from Theorem~\ref{Solutions} and Section \ref{sl}
we have the following identities

\begin{enumerate}
\item \label{easy1} $s_{\mu }(\zeta ^{I(\sigma )})=s_{\mu
^{t}}(\zeta ^{-I(\sigma^{t})})$

\item \label{easy2}$s_{\mu }(\zeta ^{I(\sigma )})=e_{n}(\zeta ^{I(\sigma
)})^{m_{n}(\mu )}s_{\tilde\mu}(\zeta ^{I(\sigma )}),\qquad
{\text where }\quad m_{n}(\mu )=\frac{|\mu |-|\tilde\mu|}{n}$.

\item \label{easy3}$s_{\mu }(\zeta ^{I(\sigma )})=e_{n}(\zeta ^{I(\sigma
)})^{k}s_{\mu ^{\vee }}(\zeta ^{-I(\sigma )})$.

\item \label{easy4} Suppose $m_{n}(\mu )=0$ then $s_{\func{rot}({
\mu})}(\zeta ^{I(\sigma )})=\zeta ^{|\sigma |}s_{\mu }(\zeta
^{I(\sigma )})$.
\end{enumerate}
\end{lemma}

\begin{proof}
Statement \eqref{easy1} follows directly from the formulae
\begin{equation}
\label{rietschformulae}
e_{r}(\zeta ^{I(\sigma )})=h_{r}(\zeta ^{-I(\sigma ^{t})})\qquad \text{and}%
\qquad h_{r}(\zeta ^{I(\sigma )})=e_{r}(\zeta ^{-I(\sigma ^{t})})\;.
\end{equation}%
of \cite{Rietsch}. From \eqref{Betheroots} and \eqref{bae2} it follows
\begin{equation*}
h_{n+1}(\zeta ^{I(\sigma ^{t})})=\cdots =h_{n+k-1}(\zeta ^{I(\sigma
^{t})})=h_{n+k}(\zeta ^{I(\sigma ^{t})})+(-1)^{k}=0
\end{equation*}%

If $\mu =\tilde{\mu}$ then there is nothing to do in \eqref{easy2}.
So assume that $\mu _{1}^{t}=n$. By expanding the determinant of the
Jacobi-Trudy formula $s_{\mu}(y)=\det (h_{\mu _{i}-i+j}(y))$
for $y=\zeta ^{I(\sigma ^{t})}$ with respect to the first row together with  \eqref{bae2} implies
\begin{equation*}
s_{\mu}(\zeta ^{I(\sigma ^{t})})=h_{n}(\zeta ^{I(\sigma
^{t})})s_{(\mu _{2}^{t},\ldots,\mu _{k}^{t})}(\zeta ^{I(\sigma ^{t})})
\end{equation*}%
Repeating the same step until all parts of size $n$ are removed, we obtain
the asserted identity \eqref{easy2}. Now \eqref{easy3} follows
from the general formula (\cite[(4.3)]{Rietsch})
\begin{equation}
\label{dualRietsch}
s_{(\mu ^{\vee })^{t}}(y_{1}^{-1},\ldots,y_{k}^{-1})=s_{(\mu
^{t})^{\vee }}(y_{1}^{-1},\ldots,y_{k}^{-1})=\frac{s_{\mu
^{t}}(y_{1},\ldots,y_{k})}{s_{(n^{k})}(y_{1},\ldots,y_{k})}\;.
\end{equation}%
To obtain \eqref{easy4} we employ the Pieri formula to find%
\begin{equation*}
s_{\func{rot}(\mu )}(\zeta ^{I(\sigma ^{t})})=e_{k}(\zeta ^{I(\sigma
^{t})})s_{\mu}(\zeta ^{I(\sigma ^{t})})=\zeta ^{\frac{k+n}{n}%
|\sigma |}s_{\mu}(\zeta ^{I(\sigma ^{t})})\;.
\end{equation*}%
and then apply \eqref{easy1} once more.
\end{proof}
We need another Schur function identity which we will frequently use in what follows.
\begin{lemma} Adopting the notation from Theorem \ref{Solutions} we set $\tilde{x}=\zeta^{|\sigma|/n}\zeta^{I(\sigma^t)}$ and $\tilde{y}=\zeta^{|\sigma|/n}\zeta^{-I(\sigma)}$, where $\sigma\in\mathfrak{P}_{\leq n-1,k}$. Then we have for any $\lambda\in\mathfrak{P}_{\leq n-1,k}$ the identities
\begin{equation}\label{Rietschbetter}
s_{\lambda^t}(\tilde x)=s_\lambda(\tilde y)\quad\text{and}\quad
\overline{s_\lambda(\tilde y)}=s_\lambda(\tilde y^{-1})=
s_{\hat\lambda ^{\ast }}(\tilde y)
\end{equation}
\end{lemma}
\begin{proof}
The first equality is immediate from Lemma \ref{Schurprop} \eqref{easy1}. %
Employing Lemma \ref{Schurprop} \eqref{easy3} and the following identities $e_n(\zeta^{I(\sigma)})^k=\zeta^{|\sigma|k}$ from \eqref{Imap0} and $|\la^\vee|=|\hat\lambda^\ast|=k n-|\lambda|$ (see Section \ref{autos}) we calculate
\begin{eqnarray}
s_{\lambda }(\zeta ^{-\frac{|\sigma |}{n}%
}\zeta ^{I(\sigma )})&=&\zeta ^{-\frac{|\sigma
||\lambda|}{n}}e_n(\zeta^{I(\sigma)})^k s_{\lambda ^{\vee
}}(\zeta^{-I(\sigma )})\\
&=& s_{\hat\lambda ^{\ast }}(\zeta ^{\frac{|\sigma |}{n}}\zeta
^{-I(\sigma )})=\overline{s_{\lambda}(\zeta ^{\frac{|\sigma |}{n}}\zeta
^{-I(\sigma )})}\notag
\end{eqnarray}
which proves the assertion.
\end{proof}

\subsection{Eigenbasis and Weyl characters} \label{sec:Weyl}
 In this section we show that the Bethe vectors form a complete set of pairwise orthogonal eigenvectors, and determine their eigenvalues in terms of Weyl characters.

 Given $g\in \op{GL}(k,\mC)$ diagonalisable with eigenvalues $t_1,\ldots, t_n$ and $\chi_\la$ the character of the finite dimensional irreducible module corresponding to the partition $\la$, we have the equality $\chi_\la(g)=s_\la(t_1,\ldots,t_n)$. We consider for $\xi\in\mZ^n$ the values $$\chi _{\lambda }(\xi ):=s_{\lambda }(e^{\frac{2\pi \iota\xi _{1}}{k+n}},\ldots,e^{\frac{2\pi \iota\xi _{k}}{k+n}})$$ (where $e$ denotes here the usual exponential function). Recall the scalar product $\langle\;,\rangle$ on $\cH$ from Proposition~\ref{phasealgebra}.

\begin{theorem}
\label{eigenbasis}
\begin{enumerate}
\item
For each $k\in\mathbb{Z}_{\geq 0}$ the Bethe vectors $b(x_\sigma)$, $\sigma \in \mathfrak{P}_{\leq n-1,k}$, form a complete set of pairwise orthogonal eigenvectors for the action of the cyclic noncommutative symmetric functions on $\cH_k[z]$.
The eigenvalues are given by the following formulae
\begin{eqnarray}
e_{r}(\mathcal{A})b(x_\sigma)&=&
\begin{cases}
h_{r}(x_\sigma)b(x_\sigma)&\text{ if $r\not=n$,}\\
h_{n}(x_\sigma)b(x_\sigma)=zb(x_\sigma)&\text{ if $r=n$.}
\end{cases}
\label{Schur spec1}\\
s_{\lambda }(\mathcal{A})b(x_{\sigma })&=&s_{\lambda ^{t}}(x_{\sigma
})b(x_{\sigma })=z^{\frac{|\lambda |}{n}}\chi _{\lambda }(\xi _{\sigma
})b(x_{\sigma }).\label{Schur spec}
\end{eqnarray}
\item The norm of the Bethe vectors is given by the following formula
\begin{eqnarray}
\langle b_\sigma,b_\sigma\rangle&=&\frac{n(n+k)^{n-1}}{|\op{Van}_\sigma|^2},
\end{eqnarray}
where $\op{Van}_\sigma$ denotes the Vandermonde determinant $\displaystyle{\prod_{1\leq i<j\leq n}({\zeta ^{I_{i}(\sigma)}-\zeta ^{I_{j}(\sigma)}}})$.
\end{enumerate}
\end{theorem}

\begin{proof}
Thanks to \eqref{Tspec} and \eqref{T is e} we have for any
solution $x$ of (\ref{bae}) the explicit formula for the eigenvalues $e_{r}(\mathcal{A})b(x)=h_{r}(x)b(x)$ for $r\not=n$ and $e_{n}(\mathcal{A})b(x)=zb(x)$ and $zb(x)=h_n(x)$ by \eqref{hformel}.

The (noncommutative versus the commutative) Jacobi-Trudy formula implies then
$s_{\lambda }(\mathcal{A})b(x)=s_{\lambda ^{t}}(x)b(x)$.
Using \eqref{Betheroots} and Lemma \ref{Schurprop}\eqref{easy1} we get $s_{\lambda ^{t}}(x_{\sigma })=z^{\frac{|\lambda |}{n}}s_{\lambda }(\zeta ^{%
\frac{|\sigma |}{n}}\zeta ^{I_{1}(\sigma )},\ldots,\zeta ^{\frac{|\sigma |}{n}%
}\zeta ^{I_{n}(\sigma )})$ and \eqref{Schur spec} follows then from the definition of $\chi_\la$.

Now consider the action of the commutative algebra $\cA^{\op{sym}}$ (from Corollary \eqref{ecomm}) on $\cH_k$. It decomposes the $\cH_k$ into simultaneous generalised eigenspaces, and the Bethe vectors are simultaneous eigenvectors with the eigenvalues as computed above. To show that they form a complete eigenbasis it is enough to see that they are in different eigenspaces, because we have previously established that there are $\dim {\mathcal{H}}_{k}$ solutions to the Bethe Ansatz equations. Assume $e_{r}(\mathcal{A})b(x_\sigma)=e_{r}(\mathcal{A})b(x_\tau)$ for all $r$. Then $h_{r}(x_\sigma)=h_{r}(x_\tau)$ for all $r$. Now the $h_r$ generate the ring of symmetric functions in $k$ variables which in turn can be identified with the (polynomial) ring of regular functions on the orbit space $\mC^k/S_k$ of the symmetric group $S_k$ acting on $\mC^k$ by permuting the variables. In particular,  $h_{r}(x_\sigma)=h_{r}(x_\tau)$ for all $r$ implies that $x_\sigma$ lies in the same $S_k$ orbit as $x_\tau$. Because of definition \eqref{Imap0} this is only possible if $\sigma=\tau$.

Next we claim that $\langle b_\sigma,b_\tau\rangle=0$ if $\sigma\not=\tau$. The adjointness of $\varphi_i$ with $\varphi_i^*$ from Proposition \ref{phasealgebra} and the definition of the elementary symmetric functions implies that $e_{n-r}(\cA)$ has adjoint  $z^{-1}e_r(\cA)$.
In particular,
\begin{eqnarray*}
0&=&\langle e_{n-r}(\cA)b(x_\sigma),b(x_\tau)\rangle-\langle b(x_\sigma), z^{-1}e_{r}(\cA) b(x_\tau)\rangle\\
&=&(\overline{h_{n-r}(x_\sigma)}-z^{-1}{h_{r}(x_\tau)}\langle b(x_\sigma),b(x_\tau)\rangle\\
&=&(z^{-1}h_{r}(x_\sigma)-z^{-1}{h_{r}(x_\tau)})\langle b(x_\sigma),b(x_\tau)\rangle.
\end{eqnarray*}
where we used Lemma \ref{hformula} for the last equality. The claim follows then as above.

To compute the norm, note that the ordinary Schur polynomials can be expressed (see e.g. \cite[Theorem 38.1]{Bump}) as a ratio $$s_\la(x_1,x_2,\ldots,x_n)=\frac{D_\la(x_1,x_2,\ldots,x_n)}{\op{Van}(x)},$$
where ${\op{Van}}$ denotes the Vandermond determinant and
\begin{eqnarray*}
D_\la(x_1,x_2,\ldots,x_n)=\det
\begin{pmatrix}
x_1^{\la_1+n-1}&x_2^{\la_1+n-1}&\ldots&x_{n-1}^{\la_1+n-1}&x_n^{\la_1+n-1}\\
x_1^{\la_2+n-2}&x_2^{\la_2+n-2}&\ldots&x_{n-1}^{\la_2+n-2}&x_n^{\la_2+n-2}\\
\vdots&\vdots&\ldots&\vdots&\vdots\\
x_1^{\la_{n-1}+1}&x_2^{\la_{n-1}+1}&\ldots&x_{n-1}^{\la_{n-1}+1}&x_n^{\la_{n-1}+1}\\
x_1^{\la_n}&x_2^{\la_n}&\ldots&x_{n-1}^{\la_n}&x_n^{\la_n}
\end{pmatrix}
\end{eqnarray*}
Let $\la\in\mathfrak{P}_{\leq n,k}$. Then $\la\in\mathfrak{P}_{\leq n-1,k}$ if and only if $\la_n=0$ which is if and only if the last row of $D_\la=D_\la(x_1,x_2,\ldots,x_n)$ consists of $1$'s only. Let us assume this is the case and expand $D_\la$ with respect to the last row to obtain $D_\la=\sum_{r=1}^n (-1)^r D_\la^{(r)}$ where the $D_\la^{(r)}$ denotes the determinant of the matrix obtained by removing the $r^{\text{th}}$ column and the $n^{\text{th}}$ row.
On the other hand, for fixed $r$, $1\leq r\leq n$, we get
\begin{eqnarray}
\label{D}
D_\la(x_1,x_2,\ldots,x_n)_{x_r=0}=
\begin{cases}
0&\text{if $\la_n\not=0$},\\
(-1)^r D_\la^{(r)}&\text{if $\la_n=0$}.
\end{cases}
\end{eqnarray}
Thanks to \eqref{bethevec} and Proposition~\ref{Schurprop}\eqref{easy1} we have
\begin{eqnarray*}
\langle b_\sigma,b_\sigma\rangle=\sum_{\lambda \in \mathfrak{P}_{\leq n-1,k}}\zeta ^{\frac{|\sigma |-|\sigma |}{
n}|\lambda |}s_{\lambda ^{t}}(\zeta ^{-I(\sigma ^{t})})s_{\lambda ^{t}}(\zeta
^{I(\sigma ^{t})})
=\sum_{\lambda \in \mathfrak{P}_{\leq n-1,k}}s_{\lambda}(\zeta ^{I(\sigma)})s_{\lambda ^{t}}(\zeta
^{I(\sigma ^{t})})\\
\end{eqnarray*}
Recall the Cauchy identity $\sum_{\lambda \in \mathfrak{P}_{\leq n,k}}s_{\lambda}(z_1,z_2,\ldots z_n)s_{\la^t}(w_1,w_2,\ldots w_k)=\prod_{i=1}^n\prod_{j=1}^k(1+z_iw_j).$
Rietsch \cite[(4.5)]{Rietsch} showed that for $(z_1,z_2,\ldots z_n)=\zeta ^{I(\sigma)}$ and $(w_1,w_2,\ldots w_k)=\zeta
^{I(\sigma ^{t})}$ it is $|\op{Van}(\zeta ^{I(\sigma)})|^2\prod_{j=1}^k(1+z_iw_j)=(n+k)$ for any $1\leq i\leq n$.  Putting everything together we get
\begin{eqnarray*}
&&|\op{Van}_\sigma|^2\sum_{\lambda \in \mathfrak{P}_{\leq n-1,k}}s_{\lambda}(\zeta ^{I(\sigma)})s_{\lambda ^{t}}(\zeta
^{I(\sigma ^{t})})\\
&=&\overline{\op{Van}_\sigma}\sum_{\lambda \in \mathfrak{P}_{\leq n-1,k}}D_{\lambda}(\zeta ^{I(\sigma)})s_{\lambda ^{t}}(\zeta
^{I(\sigma ^{t})})\\
&=&\overline{\op{Van}_\sigma}\sum_{\lambda \in \mathfrak{P}_{\leq n-1,k}}\sum_{r=1}^n (-1)^{r}D_{\lambda}^{(r)}(\zeta ^{I(\sigma)})s_{\lambda ^{t}}(\zeta
^{I(\sigma ^{t})})\\
&=&\sum_{\lambda \in \mathfrak{P}_{\leq n-1,k}}\sum_{r=1}^n (-1)^{r}D_{\lambda}^{(r)}(\zeta ^{I(\sigma)})D_{\lambda ^{t}}(\zeta
^{I(\sigma ^{t})})\\
&\stackrel{\eqref{D}}{=}&\sum_{\lambda \in \mathfrak{P}_{\leq n-1,k}}\sum_{r=1}^n
(-1)^r D_{\lambda}(y_1^{(r)},y_2^{(r)},\ldots y_n^{(r)})D_{\lambda ^{t}}(\zeta^{I(\sigma ^{t})})\\
&\stackrel{\eqref{D}}{=}&\sum_{\lambda \in \mathfrak{P}_{\leq n,k}}\sum_{r=1}^n
(-1)^r D_{\lambda}(y_1^{(r)},y_2^{(r)},\ldots y_n^{(r)})D_{\lambda ^{t}}(\zeta^{I(\sigma ^{t})})\\
&=&\sum_{r=1}^n\sum_{\lambda \in \mathfrak{P}_{\leq n,k}}
(-1)^r D_{\lambda}(y_1^{(r)},y_2^{(r)},\ldots y_n^{(r)})D_{\lambda ^{t}}(\zeta^{I(\sigma ^{t})})
\end{eqnarray*}
where $y_i^{(r)}=\zeta ^{I(\sigma)}_i$ if $i\not=r$ and $y_r^{(r)}=0$, and with Rietsch's formula gives
\begin{eqnarray*}
\langle b_\sigma,b_\sigma\rangle=\sum_{\lambda \in \mathfrak{P}_{\leq n-1,k}}s_{\lambda}(\zeta ^{I(\sigma)})s_{\lambda ^{t}}(\zeta
^{I(\sigma ^{t})})&=&
|\op{Van}_\sigma|^{-2}\sum_{r=1}^n (n+k)^{n-1}\\
&=&|\op{Van}_\sigma|^{-2}n(n+k)^{n-1}.
\end{eqnarray*}
Hence the formula follows.
\end{proof}

\begin{corollary}
\label{allascomm}
The known determinant
formulae from the ring of commutative functions are also true for the
cyclic noncommutative functions, in particular
\begin{equation*}
s_{\lambda }(\mathcal{A})=\det (h_{\lambda _{i}-i+j}(%
\mathcal{A}))_{1\leq i,j\leq n}\;.
\end{equation*}
The cyclic noncommutative Schur polynomials pairwise commute.
\end{corollary}
\begin{proof} Specialising in \eqref{Schur spec} to a horizontal $r$-strip it follows from Definition \ref{edef} that
\begin{eqnarray}
h_{r}(\mathcal{A})b(x)&=&\det (e_{1-i+j}(\mathcal{A}))_{1\leq i,j\leq r}b(x)\notag\\
&=&\det (h_{1-i+j}(x))_{1\leq i,j\leq r}b(x)=e_r(x)b(x)\;.\label{nchspec}
\end{eqnarray}
Employing the familiar relations form the ring of commutative symmetric functions we therefore have
\begin{multline*}
s_\lambda(\cA)b(x)=s_{\lambda^t}(x)b(x)=\\
\det (e_{\lambda _{i}-i+j}(x))_{1\leq i,j\leq n}b(x)=
\det (h_{\lambda _{i}-i+j}(\mathcal{A}))_{1\leq i,j\leq n}b(x)\;.
\end{multline*}
According to Theorem \ref{eigenbasis} the Bethe vectors \eqref{bethevec} form a basis in each subspace $\mathcal{H}_k[z]$ for any $k\in\mathbb{Z}_{\geq 0}$ and, hence, the last identity implies the assertion.
\end{proof}

\begin{corollary}\label{hdef}
For $0\leq r\leq n-1$ the complete symmetric polynomials $h_r(\cA)$ defined in \eqref{nchdef} have the following explicit, simpler form
\begin{equation}
h_{r}(\cA)=\sum_{|J|=r}\tprod_{j%
\in J}^{\circlearrowright }a_{j},  \label{nchid}
\end{equation}%
where the sum runs over all multisets $J=\{j_{1},\ldots,j_{r}\}$ (i.e. $s\neq t$ does not necessarily imply that $j_{s}\neq j_{t}$).
\end{corollary}
\begin{example}
If $n=4$ then
\begin{multline*}
h_{3}(\cA)=
(\sum_{i=0}^{3}a_{i}^{3})+a_{1}^{2}a_{2}+a_{1}a_{2}^{2}+a_{1}^{2}a_{3}+%
a_{1}a_{3}^{2}+a_{0}^{2}a_{1}+a_{0}a_{1}^{2}
+a_{2}^{2}a_{3}+a_{2}a_{3}^{2}\\
+a_{2}^{2}a_{0}+a_{2}a_{0}^{2}+a_{3}^{2}a_{0}+a_{3}a_{0}^{2}
+a_{1}a_{2}a_{3}+a_{0}a_{1}a_{2}+a_{2}a_{3}a_{0}\;.
\end{multline*}
\end{example}

\begin{proof}[Proof of Corollary \ref{hdef}]
Because of Theorem \ref{eigenbasis}, Proposition \ref{faithfulness} and \eqref{nchspec} it suffices to show that for any $k\in\mathbb{Z}_{\geq 0}$ the Bethe vectors are eigenvectors of the polynomials defined on the right hand side of equation \eqref{nchid} and have eigenvalues $e_r(x)$ for $0\leq r\leq n-1$.
The case $r=0$ is obviously true, hence assume $r>0$.
Let us start by describing the action of the polynomial in \eqref{nchid} on a weight $\hat\lambda$. Given a composition $p=(p_{0},p_{1},...,p_{n-1})$ of $r<n$, there must be at least one $p_{i}=0$ and, hence, we may rewrite the polynomial in \eqref{nchid} as
\begin{equation}\label{hcomp}
\sum_{|J|=r}\tprod_{j\in J}^{\circlearrowright}a_{j}%
=\sum_{p\vdash r}\tprod_{0\leq j\leq n-1}^{\circlearrowright}a_{j}^{p_j},\quad \tprod_{0\leq j\leq n-1}^{\circlearrowright}a_{j}^{p_j}{=}
a_{i+1}^{p_{i+1}}\cdots
a_{n-1}^{p_{n-1}}a_{0}^{p_{0}}a_{1}^{p_{1}}\cdots a_{i-1}^{p_{i-1}}\;.
\end{equation}
Set  $\hat{\mu}=a_{i+1}^{p_{i+1}}\cdots a_{0}^{p_{0}}a_{1}^{p_{1}}\cdots
a_{i-1}^{p_{i-1}}\hat{\lambda}$. Let us first assume that $p_{0}=0$.
In this case each $a_{j}$ acts on $\hat{\lambda}$ by adding a box in the $(j+1)^{
\text{th}}$ row (of the associated Young diagram $\hat{\op{P}}(\hat\lambda)$ if allowed).
Then either $\hat{\mu}$ is the null vector or $\nu/\la=(r)$ with $p_{j-1}$ boxes in the $j^{\text{th}}$ row where $\nu\in\mathfrak{P}_{\leq n,k}$ and $\la\in\mathfrak{P}_{\leq n-1,k}$ are obtained from respectively $\hat\mu$ and $\hat\la$ by removing $m_0(\hat\la)$ $n$-columns.

Now assume that $p_{0}>0$ and some other $p_{i}=0$. Then $a_{0}^{p_{0}}$
acts on $\hat\la$ by removing  $p_{0}$ columns of height $n$ (if possible), adding $p_{0}$ boxes in the first row, and then multiplying with $z^{p_{0}}$. It is mapped to zero otherwise.
Together with \eqref{bethevec} we arrive at the formula
\begin{equation}\label{haction}
\sum_{|J|=r}\bigl(\tprod_{j\in J}^{\circlearrowright }a_{j}\bigr)b(x)%
=\sum_{\hat{\lambda}\in P_{k}^{+}}%
\sum_{\substack{ \nu
\in \mathfrak{P}_{\leq n,k} \\ \nu/\lambda =(r)}}
z^{\nu_1-\la_1} s_{\hat{\lambda}^{t}}(x^{-1})\hat{\mu}(\nu)
\end{equation}
where $\hat{\mu}(\nu)$ is the diagram obtained from $\nu$ by adding $(k-\nu_1)$ columns of height $n$.

Now we multiply the Bethe vector $b(x)$ from \eqref{bethevec} by $e_{r}(x)$. Let $\tilde{x}(\sigma )=\zeta ^{\frac{|\sigma |}{n}}\zeta ^{I(\sigma ^{t})}
$ be the Bethe roots \eqref{Betheroots} with the $z$-dependence removed. Using the identity and the Pieri formula for commutative Schur functions we compute with help of \eqref{Rietschbetter}
\begin{eqnarray}
e_{r}(x)b(x) &=&\sum_{\hat{\mu}\in P_{k}^{+}}e_{r}(x)s_{\hat{\mu}%
^{t}}(x^{-1})\hat{\mu}=\sum_{\hat{\mu}\in P_{k}^{+}}z^{-\frac{|\hat{%
\mu}|-r}{n}}e_{r}(\tilde{x})s_{(\mu^{\ast })^{t}}(\tilde{x})\hat{%
\mu} \\
&=&\sum_{\hat{\mu}\in P_{k}^{+}}z^{-\frac{|\hat{\mu}|-r}{n}%
}\sum_{\nu/\mu^{\ast }=(r)}s_{\nu^{t}}(\tilde{x})\hat{\mu},\label{tocompare}
\end{eqnarray}
Since $s_{\nu^{t}}(\tilde x)=0$ if $\nu_{1}>k$ and $\mu^{\ast }\in
\mathfrak{P}_{\leq n-1,k}$, we may assume in \eqref{tocompare} that $\nu\in \mathfrak{P}%
_{\leq n,k}$, in particular $\nu^{\vee }$ is defined. Then $\nu/\mu ^{\ast }=(r)$ is equivalent to $\hat{\mu}/\nu^{\vee }=(r)$ by taking complements and using the equation  $\hat{\mu}=(\mu^{\ast })^{\vee }
$ from \eqref{complement}. So, by \eqref{Rietschbetter}
\begin{eqnarray}
e_{r}(x)b(x)
&=&\sum_{\hat{\mu}\in P_{k}^{+}}\Biggl(\,\sum_{\substack{\nu^\vee\in \mathfrak{P}_{\leq n,k},\\\hat{\mu}/\nu^{\vee
}=(r)}}s_{(\nu ^{\vee })^{t}}(x^{-1})\hat{\mu}\;\Biggr).
\end{eqnarray}
Let $\nu'$ be obtained from $\nu ^{\vee }$ by removing all columns of height $n$. Let $c=c(\nu^\vee)$ be the number of columns removed. Then the condition $\hat{\mu}/\nu ^{\vee }=(r)$ is equivalent to $\mu^{(\nu)}/\nu'=(r)$ if $\mu^{(\nu)}$ denotes the partition obtained from $\hat{\mu}$ by removing $c=k-\mu^{(\nu)}_1$ columns of height $n$. Denote by $\widehat{\nu'}$ the preimage of $\nu'$ under $\hat{\op{P}}$, then the Bethe Ansatz equations imply $s_{(\nu^{\vee
})^{t}}(x^{-1})=z^{\mu_1^{(\nu)}-\nu'_1}s_{(\widehat{\nu'})^{t}}(x^{-1})$ (as $\nu ^{\vee }$ and $\widehat{\nu'}$ differ by $k-\nu^\vee_1=\mu_1^{(\nu)}-\nu'_1$ columns of length $n$, \eqref{hformel} implies
$h_n(x^{-1})=\overline{h_n(x)}=z^{-1}$ and then the claim follows from Lemma \ref{Schurprop}, (\ref{easy1}) and (\ref{easy2}).) Altogether we obtain

\begin{eqnarray}
e_{r}(x)b(x)
&=&
\sum_{\hat{\mu}\in P_{k}^{+}}
\sum_{\substack{ \nu^\vee
\in \mathfrak{P}_{\leq n,k} \\
\mu^{(\nu)}/\nu'=(r)}}z^{\mu_1^{(\nu)}-\nu'_1}s_{(%
\widehat{\nu}')^{t}}(x^{-1})\hat{\mu},\label{etimesb}
\end{eqnarray}
Each summand belonging to $(\hat\mu,\nu^\vee)$ is in one-to-one correspondence with a summand in \eqref{haction} if we identify the pair $(\mu^{(\nu)},\nu')$ in \eqref{etimesb} with the same pair $(\nu,\la)$ in \eqref{haction}.
\end{proof}

Corollary \ref{hdef} prompts us to introduce
the following counterpart of the transfer matrix. Namely, we set%
\begin{equation}
Q(u)=\sum_{r\geq 0}h_{r}(\mathcal{A})u^{r}\;  \label{Q}
\end{equation}%
as the generating operator for the
noncommutative complete symmetric functions where, similar as before, we
choose the convention
$h_{0}(\mathcal{A})=1$.
Note that this operator is well-defined since upon
restricting to a subspace ${\mathcal{H}}_{k}[z]$ with fixed particle number $k$ only
the first $(k+1)$ terms in the sum (\ref{Q}) do not vanish, i.e. $h_{r}(%
\mathcal{A})_{k}=0$ for $r>k$. Here and in the following $F_k$ denotes the restriction of an endomorphism $F$ of $\cH[z]$ to the subspace  $\cH_{k}[z]$. The operator $h_{r=k}(\mathcal{A})_{k}$ becomes the cyclic translation operator corresponding to the Dynkin diagram automorphism $\op{rot}$ if we specialise $z=1$. Thus, with $Q(u)_k$ being a polynomial in $u$ of degree $k$ the result \eqref{Tspec} can be rewritten as the following functional equation,
\begin{equation}\label{TQeqn}
T(u)_k Q(-u)_k =\left[Q(-uq)_k+z u^n q^k Q(-u q^{-1})_k\right]_{q=0}\;.
\end{equation}%
 This equation is the `crystal limit' of what is known as Baxter's $TQ$-equation in the literature on quantum integrable systems. In particular, it is immediate from our results that for any pair $u,v\in \mathbb{C}$ we have
\begin{equation}
\lbrack Q(u),Q(v)]=[Q(u),T(v)]=0\;.  \label{TQ commute}
\end{equation}
 The operator $Q(u)$ can therefore be identified as the analogue of Baxter's {\em Q-operator} for the XYZ and XXZ models \cite{Baxter}; further details will be presented in a forthcoming paper \cite{Korff}.

\subsection{The combinatorial fusion ring and Verlinde algebra}
Employing the cyclic noncommutative symmetric functions generated from the transfer matrix of the phase model, we now show that the $k$-particle space $\mathcal{H}_k[z]$ can be turned into a commutative, associative and unital algebra. To explicitly compute its structure constants we first need the following result:
\begin{proposition}[Transformation matrix]
\label{newlemma}
Let $\hat\sigma\in P_k^+$ and $\sigma=\op{P}(\hat\sigma)$ the corresponding partition.
Let $\xi_\sigma={\sigma+\rho}$, where $\rho=(\frac{n+1}{2}-1,\dots,\frac{n+1}{2}-n)$. Then we have the identity
\begin{eqnarray}
\label{Weyl}
\langle \hat{\lambda},b(x_{\sigma })\rangle=z^{-\frac{|\hat{\lambda}|}{n}}s_{\lambda
}(\zeta ^{-\frac{|\sigma |}{n}}\zeta^{I(\sigma)})=z^{-\frac{|\hat{\lambda}|}{n}}\overline{\chi_{\la}(\xi_\sigma)}
=z^{-\frac{|\hat{\lambda}|}{n}}\chi_{\la^\ast}(\xi_\sigma).
\end{eqnarray}
\end{proposition}

\begin{proof}
From the definition \eqref{bethevec} we have $\langle \hat{\lambda},b(x_{\sigma })\rangle =s_{\hat{\lambda}^{t}}(x_{\sigma
}^{-1})$. We claim that
\begin{equation}\label{schurred}
s_{\hat{\lambda}^{t}}(x_{\sigma }^{-1})=h_n(x_{\sigma}^{-1})^{m_{n}(\hat{\lambda})}s_{\lambda
^{t}}(x_{\sigma }^{-1})=z^{-m_{n}(\hat{\lambda})}s_{\lambda
^{t}}(x_{\sigma }^{-1})=z^{-\frac{|\hat{\lambda}|-|\lambda |}{n}}s_{\lambda
^{t}}(x_{\sigma }^{-1}).
\end{equation}%
The first equality here is just Lemma \ref{Schurprop} \eqref{easy2} and \eqref{rietschformulae}. The identity $h_{n}(x_{\sigma }^{-1})=z^{-1}h_{0}(x_{\sigma })^{-1}=z^{-1}=h_{n}(x_{\sigma })^{-1}$ holds by Lemma \ref{hformula} and \eqref{Schur spec1}. Since the Schur function $s_{\lambda }$ is a homogeneous polynomial of degree $|\lambda |$, Lemma \ref{Schurprop} \eqref{easy1} and the explicit form of the Bethe roots \eqref{Betheroots} imply
$s_{\lambda ^{t}}(x_{\sigma }^{-1})=z^{-\frac{|\lambda |}{n}}s_{\lambda
}(\zeta ^{-\frac{|\sigma |}{n}}\zeta ^{I(\sigma )})$.
Hence,
\begin{eqnarray}
\label{Weyl}
\langle \hat{\lambda},b(x_{\sigma })\rangle=z^{-|\hat{\lambda}|/n}s_{\lambda
}(\zeta ^{-\frac{|\sigma |}{n}}\zeta^{I(\sigma)}).
\end{eqnarray}
Finally we use the identity \eqref{Rietschbetter}.
\end{proof}

The main theorem of this section is the following ring structure on $\mathcal{H}_{k}[z]$:

\begin{theorem}[Combinatorial fusion ring]
\label{combinatorial}
Fix $k\in \mathbb{Z}_{\geq 0}$ and consider the $k$-particle subspace $
\mathcal{H}_{k}[z]\subset\mathcal{H}[z]$. The assignment
\begin{equation}
(\hat{\lambda},\hat{\mu})\mapsto \hat{\lambda}\circledast \hat{\mu}%
:=s_{\hat\lambda }(\mathcal{A})\hat{\mu}  \label{Schurproduct}
\end{equation}%
for basis elements $\hat{\lambda},\hat{\mu}\in P_{k}^{+}$ turns $\mathcal{H}_{k}[z]$ into a commutative, associative and unital  $\mathbb{C}[z]$-algebra $\cF_{\op{comb}}$.
\end{theorem}

\begin{proof}
The unit element is obviously given by the weight corresponding to the empty
partition as $s_{\emptyset }(\mathcal{A})=1_{\mathcal{H}_{k}[z]}$. We have to check the commutativity and associativity of the
product. To achieve this we compute the matrix elements
\begin{equation}
s_{\hat\lambda }(\mathcal{A})\hat{\mu}=\sum_{\hat{\nu}\in P^{+}_{k}}\langle \hat{%
\nu},s_{\hat\lambda }(\mathcal{A})\hat{\mu}\rangle \hat{\nu}
\end{equation}
in the eigenbasis of Bethe vectors $\{b(x_{\sigma
})\}_{\sigma \in \mathfrak{P}_{\leq k,n-1}}$. Note that \eqref{schurred} together with Theorem \ref{eigenbasis} implies the identity
\begin{equation}\label{ncschurred}
s_{\hat\lambda }(\mathcal{A})=
z^{\frac{|\hat\lambda|-|\lambda|}{n}}s_{\lambda }(\mathcal{A}).
\end{equation}
Moreover, thanks to Theorem \ref{eigenbasis} and Proposition \ref{newlemma} we have
\begin{eqnarray}
\langle \hat{\nu},s_{\hat\lambda }(\mathcal{A})\hat{\mu}\rangle &=&\sum_{\sigma
\in \mathfrak{P}_{\leq n-1,k}}\frac{\langle \hat{\nu},s_{\hat\lambda }(\mathcal{A%
})b(x_{\sigma })\rangle \langle b(x_{\sigma }),\hat{\mu}\rangle }{%
\langle b(x_{\sigma }),b(x_{\sigma })\rangle }  \notag \\
&=&\frac{z^{\frac{|\hat\lambda |+|\hat\mu |-|\hat\nu |}{n}}}{n(k+n)^{n-1}}\sum_{\sigma \in
\mathfrak{P}_{\leq k,n-1}}\zeta ^{|\sigma |\frac{|\lambda |+|\mu |-|\nu |}{n}%
}\frac{s_{\lambda }(\zeta ^{-I(\sigma )})s_{\mu }(\zeta ^{-I(\sigma
)})s_{\nu }(\zeta ^{I(\sigma )})}{\prod_{i<j}^{n}|\zeta ^{I_{i}(\sigma
)}-\zeta ^{I_{j}(\sigma )}|^{-2}}
\label{permanentformula}
\end{eqnarray}%
from which it is obvious that the product is commutative as the matrix
element is invariant under exchanging $\hat\lambda $ and $\hat\mu $. Associativity is also clear:
\begin{eqnarray*}
\hat{\lambda}\circledast (\hat{\mu}\circledast \hat{\nu})
&=&\displaystyle\sum_{\hat{\rho}
\in P_{k}^{+}}~\langle \hat{\rho},s_{\hat\mu
}(\mathcal{A})\hat{\nu}\rangle \hat{\lambda}\circledast \hat{\rho}
=
\displaystyle\sum_{\hat{\rho},\hat{\rho}^{\prime }\in
P_{k}^{+}}\langle \hat{\rho}^{\prime },s_{\hat\lambda }(%
\mathcal{A})\hat{\rho}\rangle \langle \hat{\rho},s_{\hat\nu }(\mathcal{A})\hat{%
\mu}\rangle\hat{\rho}^{\prime }\\
&=&\displaystyle\sum_{\hat{\rho}^{\prime }\in P_{k}^{+}}\langle \hat{%
\rho}^{\prime },s_{\hat\lambda }(\mathcal{A})s_{\hat\nu }(\mathcal{A})\hat{\mu}%
\rangle\hat{\rho}^{\prime } =\displaystyle\sum_{\hat{\rho}^{\prime }\in P_{k}^{+}}\langle
\hat{\rho}^{\prime },s_{\hat\nu }(\mathcal{A})s_{\hat\lambda }(\mathcal{A})\hat{\mu}%
\rangle\hat{\rho}^{\prime } \\
&=&\displaystyle\sum_{\hat{\rho},\hat{\rho}^{\prime }\in P_{k}^{+}}\langle \hat{\rho}^{\prime },s_{\hat\nu }(\mathcal{A})\hat{\rho}\rangle \langle
\hat{\rho},s_{\hat\lambda }(\mathcal{A})\hat{\mu}\rangle\hat{\rho}^{\prime
}
=\displaystyle\sum_{\hat{\rho}\in
P_{k}^{+}}\langle \hat{\rho},s_{\hat\lambda }(%
\mathcal{A})\hat{\mu}\rangle\hat{\nu}\circledast \hat{\rho}\\
&=&\displaystyle\sum_{\hat{\rho}\in
P_{k}^{+}}\langle \hat{\rho},s_{\hat\lambda }(%
\mathcal{A})\hat{\mu}\rangle\hat{\rho}\circledast \hat{\nu}
=(\hat{\lambda}\circledast \hat{\mu})\circledast \hat{\nu}\;,
\end{eqnarray*}
where in the first and last line we have exploited commutativity of the product and in the second line we used the non-degeneracy of the bilinear form and the commutativity of the polynomials $%
s_{\hat\lambda }(\mathcal{A})$, $s_{\hat\nu }(\mathcal{A})$ (see Corollary \ref{allascomm}).
\end{proof}

\begin{definition}
\label{modularS}
{\rm
For $\hat{\la}, \hat{\mu},\hat{\sigma}\in P_k^+$ we abbreviate $\mathcal{S}_{0,\hat{\sigma}}=\frac{1}{\sqrt {n(n+k)^{n-1}}}|\op{Van}_{\sigma}|$ and define $\mathcal{S}_{\hat{\lambda}\hat{\sigma}}=\mathcal{S}_{0\hat{\sigma}}\chi_{\lambda }(\xi _{\sigma })$ with $\xi _{\sigma }=-(\op{P}(\sigma)+\rho)$. We call the matrix $\mathcal{S}:=\mathcal{S}_{\hat{\lambda},\hat{\mu}}$ with columns and rows indexed by $P_k^+$(or by $\mathfrak{P}_{\leq n-1,k}$  via \eqref{weight2part}) the {\it modular $S$-matrix}.
}
\end{definition}

\begin{remark}\label{modinv}{\rm
From Theorem \eqref{eigenbasis} it follows directly that the matrix $\mathcal{S}$ is {\it unitary}, i.e. $\sum_{\hat{\lambda}\in P^{+}_{k}}\mathcal{S}%
_{\hat{\lambda}\hat{\sigma}}\mathcal{\bar{S}}_{\hat{\lambda}\hat{\rho}%
}=\delta _{\hat{\sigma},\hat{\rho}}.$ From Lemma \ref{Schurprop} \eqref{easy4} it follows that $\mathcal{S}$ is almost rotation invariant: $\mathcal{S}_{\func{rot}(\hat{\lambda})\hat{\mu}}=e^{-\frac{2\pi \iota
}{n}|\mu |}\mathcal{S}_{\hat{\lambda}\hat{\mu}}$. Moreover, we have the {\it duality} $\mathcal{S}_{\hat{\lambda}\hat{\mu}}=\sqrt{k/n}e^{\frac{2\pi
\iota}{kn}|\lambda ||\mu |}\overline{\cS_{\hat{\lambda}^{t}\hat{\mu}^{t}}}$ and {\it charge conjugation formula} $\mathcal{\bar S}_{\hat\lambda\hat\mu}=\mathcal{ S}_{\hat\lambda^\ast\hat\mu}$.  %
Let $\mathcal{T}$ be the diagonal matrix with entries $e^{2\pi\iota\op{m}(\hat\la)}$ where $\op{m}(\hat\la):=\frac{|\la+\rho|^2}{2(k+n)}-\frac{|\rho|^2}{2n}$ is the so-called {\em modular anomaly}. Then $\op{S}=\left(
\begin{smallmatrix} 0 & -1 \\ 1 & 0 \end{smallmatrix}\right)\mapsto\mathcal{S}$, $\op{T}=\left(\begin{smallmatrix}1 & 1 \\0 & 1\end{smallmatrix}\right)\mapsto\mathcal{T}$ provides a representation of $\op{PSL}(2,\mathbb{Z})$, while $\op{S}\mapsto\mathcal{\bar S}$, $\op{T}\mapsto\mathcal{T}$ yields a representation of $\op{SL}(2,\mathbb{Z})$. Thus, together with complex conjugation this turns $\cF_{\op{comb}}$ %
into an abstract Verlinde algebra in the sense of \cite[0.4.1]{Cherednik}. While defined here in a combinatorial setting these actions of the modular group correspond to the familiar modular invariance in the context of conformal field theory; see e.g. \cite[Chapters 10, 14 and 16]{CFTbook}.
}
\end{remark}

The following key result relates our combinatorial fusion ring to the $\mathfrak{\widehat{sl}}(n)_k$ Verlinde algebra and what is usually called the modular $S$-matrix in the literature. It states that our modular $S$- matrix satisfies the famous Kac-Peterson formula \cite{KacPeterson}.
\begin{proposition}[Kac-Peterson formula]
\label{KacPeterson}
Let ${\hat{\lambda},\hat{\sigma}}\in P_k^+$. Then (in the notation of Section \ref{sl})
\begin{equation*}
\mathcal{S}_{\hat{\lambda}\hat{\sigma}}=\frac{e^{\iota\pi n(n-1)/4}}{\sqrt{%
n(k+n)^{n-1}}}\sum_{w\in S_{n}}(-1)^{\ell (w)}e^{-\frac{2\pi \iota}{k+n}%
(\sigma +\rho,w(\lambda +\rho ) )}
\end{equation*}
Here $\hat{\lambda}=k\hat{\omega}_{0}+\lambda$ with finite part $\la=\sum_{i=1}^{n}(\lambda _{i}-|\lambda |/n)\varepsilon _{i}$
and $\rho =\sum_{i=1}^{n}\left( \frac{n+1}{2}-i\right)\varepsilon _{i}\;$.
\end{proposition}

\begin{remark}
\label{sin}
{\rm
Using Weyl's denominator formula one can also deduce an explicit expression
\begin{eqnarray}
\mathcal{S}_{0\hat{\sigma}} &=&\frac{1}{\sqrt{n(k+n)^{n-1}}}~\prod_{\alpha
>0}2\sin \frac{\pi (\rho +\sigma,\alpha )}{k+n} \\
&=&\frac{e^{\iota\pi n(n-1)/4}}{\sqrt{n(k+n)^{n-1}}}\sum_{w\in S_{n}}(-1)^{\ell
(w)}e^{-\frac{2\pi \iota}{k+n}(\sigma +\rho,w\rho)}\;.\label{S0}
\end{eqnarray}
}
\end{remark}

\begin{remark}[Uniqueness theorem]\label{GoodmanWenzl}{\rm
Goodman and Wenzl (\cite[Theorem 3.2]{GoodmanWenzl}, see also \cite{GoodmanNakanishi}) established a characterisation of the fusion ring $\cF$ by the data of its basis $P_k^+$, the associativity and the explicit formulae for $\hat\la\circledast\hat\mu=\hat\mu\circledast\hat\la$ where $\hat\la,\hat{\mu}\in P_k^+$ with $\hat\mu$ of the form $(1^r)$. Theorem \ref{combinatorial} and the Pieri-type formulae from Proposition~\ref{Pieri} can then be used to verify the assumption of \cite{GoodmanWenzl} and so to deduce an isomorphism of rings $\cF_{\op{comb}}\cong\cF$. We will establish such an isomorphism (in Theorem \ref{combinatorialfusion}) in a different way, involving the Verlinde formula and the modular $S$-matrix.
}
\end{remark}

\begin{proof}[Proof of Proposition \ref{KacPeterson}]
Let  $\delta =(n-1,n-2,\ldots,1,0)$ and for $j=1,\ldots,n$ set
\begin{equation*}
y_{j}=e^{-\frac{2\pi \iota}{k+n}~(\sigma +\rho ,\varepsilon _{j})}
=\zeta ^{\frac{|\sigma |}{n}-I_{j}(\sigma )},\quad y^\mu=y_1^{\mu_1}y_2^{\mu_2}\cdots y_n^{\mu_n}.
\end{equation*}
Using \eqref{Weyl}, with $z=1$, and relating Schur polynomials with Weyl's character formula \cite[p.40 in particular (3.1)]{MacDonald} we get that the term $\langle b(x_{\sigma }),\hat{\lambda}\rangle$ equals
\begin{eqnarray*}
s_{\lambda }(\zeta ^{\frac{|\sigma |}{n}}\zeta ^{-I(\sigma )})=\frac{%
{\displaystyle\sum_{w\in S_{n}}(-1)^{\ell (w)}w(y^{\la+\delta})}}{\displaystyle{\sum_{w\in
S_{n}}(-1)^{\ell (w)}w(y^\delta)}}=\frac{\displaystyle\sum_{w\in S_{n}}(-1)^{\ell
(w)}e^{-\frac{2\pi \iota}{k+n}\sum_{j=1}^n\left(x ,\varepsilon _{j}\right)(n+\lambda
_{j}-j)}}{\displaystyle\sum_{w\in S_{n}}(-1)^{\ell (w)}e^{-\frac{2\pi \iota}{k+n}\sum_{j=1}^n\left(x ,\varepsilon _{j}\right)(n-j)}}
\end{eqnarray*}
with $x=w(\sigma +\rho)$. On the other hand $\left(w(x),\epsilon_j\right)(n+\la_j-j)$ equals
\begin{eqnarray*}
\left(w(x),\left(\la_j-\frac{|\la|}{n}\right)\epsilon_j\right)
-\left(w(x),\left(j-\frac{n+1}{2}\right)\epsilon_j\right)+
\left(w(x),\left(\frac{n-1}{2}+\frac{|\la|}{n}\right)\epsilon_j\right).
\end{eqnarray*}
From the definitions, and using $\sum_{j}I_{j}(\sigma )=|\sigma |$, we get
\begin{eqnarray*}
\langle b(x_{\sigma}),\hat{\lambda}\rangle
&=&\frac{\sum_{w\in S_{n}}(-1)^{\ell (w)}e^{-\frac{2\pi \iota}{k+n}(\sigma +\rho
,w(\rho +\lambda ))}}{\sum_{w\in S_{n}}(-1)^{\ell (w)}e^{-\frac{2\pi \iota}{k+n}%
(\sigma +\rho ,w\rho )}}
\end{eqnarray*}
The assertion then follows from \eqref{S0}.
\end{proof}

\subsection{The Verlinde algebra and the modular S-matrix}
The Verlinde algebra is the fusion algebra of the integrable highest weight modules of level $k$. It plays a prominent role in conformal field theory (for details see e.g. \cite{CFTbook}, \cite{Kac}, \cite{KacPeterson}).

 More precisely, the {\it Verlinde algebra} $V_{k}=V_{k}(\widehat{\mathfrak{sl}}(n),\mathbb{C})=\mC\otimes_\mathbb{Z}\cF$ is the $\mC$-algebra with basis indexed by the elements from $P_{k}^{+}$ together
with the multiplication (called {\it fusion product})
\begin{equation}
\hat{\lambda}\ast \hat{\mu}:=\sum_{\hat{\nu}\in P^{+}_{k}}\mathcal{N}_{\hat{%
\lambda}\hat{\mu}}^{(k),\hat{\nu}}\hat{\nu},  \label{fusionp}
\end{equation}%
where the structure constants, known as {\it fusion coefficients} are given in terms of the {\it Verlinde formula} \cite{Verlinde}
\begin{equation}
\label{Verlinde}
\mathcal{N}_{\hat{\lambda}\hat{\mu}}^{(k),\hat{\nu}}=\sum_{\hat{\sigma}\in
P^{+}_{k}}\frac{\mathcal{S}_{\hat{\lambda}\hat{\sigma}}\mathcal{S}_{\hat{\mu}%
\hat{\sigma}}\mathcal{S}_{\hat{\nu}^\ast\hat{\sigma}}}{\mathcal{S}_{0\hat{%
\sigma}}}
\end{equation}
and the $\mathcal{S}_{\hat{\lambda}\hat{\sigma}}$ are given by the Kac-Peterson formula (precisely the formula from Proposition \eqref{KacPeterson}). We already saw that one could alternatively define $\mathcal{S}_{\hat{\lambda}\hat{\sigma}}=\mathcal{S}_{0\hat{\sigma}}\chi
_{\lambda }(\xi _{\sigma })$, $\xi _{\sigma }=-(\op{P}(\hat\sigma)+\rho),$
with $\mathcal{S}_{0\hat{\sigma}}$ as above.

\begin{theorem}[Combinatorial description of the Verlinde algebra]
\label{combinatorialfusion}
The combinatorial algebra $\cF$ specialises to the integral Verlinde (or fusion) algebra, i.e.
\begin{equation}
\cF_{\op{comb}}/\langle z-1\rangle \cong \cF.
\end{equation}%
In particular, we have for $z=1$ the following equalities%
\begin{equation}
\hat{\lambda}\ast \hat{\mu}=\hat{\lambda}\circledast \hat{\mu}=s_{\hat\lambda }(%
\mathcal{A})\hat{\mu}=s_{\lambda }(%
\mathcal{A})\hat{\mu},\qquad \forall \hat{\lambda},\hat{\mu}\in P^{+}_{k}\;.
\label{result}
\end{equation}
identifying the fusion product (\ref{fusionp}) with the combinatorially defined product (\ref{Schurproduct}).
\end{theorem}

\begin{proof}
We relate the matrix elements of the cyclic noncommutative Schur
functions to the fusion constants:
\begin{eqnarray}
\langle \hat{\nu},s_{\hat\lambda }(\mathcal{A})\hat{\mu}\rangle
&{=}&\frac{z^{\frac{|\hat\lambda |+|\hat\mu |-|\hat\nu |}{n}}}{n(k+n)^{n-1}}\sum_{\sigma\in \mathfrak{P}_{\leq k,n-1}}
\frac{\chi _{\lambda }(\xi _{\sigma })\chi _{\mu }(\xi _{\sigma })%
\chi _{\nu^\ast }(\xi _{\sigma })}{\prod_{i<j}^{k}|\zeta ^{I_{i}(\sigma
)}-\zeta ^{I_{j}(\sigma)}|^{-2}}  \label{Verlindeexplicit}  \\
&=&z^{\frac{|\hat\lambda |+|\hat\mu |-|\hat\nu |}{n}}\sum_{\hat{\sigma}\in P^{+}_{k}}%
\frac{\mathcal{S}_{\hat{\lambda}\hat{\sigma}}\mathcal{S}_{\hat{\mu}
\hat{
\sigma}}\mathcal{S}_{\hat{\nu^\ast}\hat{\sigma}}}{\mathcal{S}_{0\hat{\sigma}%
}}{=}z^{\frac{|\hat\lambda |+|\hat\mu |-|\hat\nu |}{n}}\mathcal{N}_{\hat{\lambda}\hat{\mu}%
}^{(k),\hat{\nu}}\;.\notag
\end{eqnarray}%
The first equality here holds thanks to \eqref{permanentformula}, \eqref{Weyl} and \eqref{Rietschbetter}, the second by Definition \ref{modularS}, and the last one is the Verlinde formula \eqref{Verlinde}.  Specialising $z\mapsto 1$ and employing \eqref{ncschurred} gives the result.
\end{proof}

\begin{remark}{\rm
Note that it is known that the fusion coefficients $\mathcal{N}_{\hat{\lambda%
}\hat{\mu}}^{(k),\hat{\nu}}$ are nonnegative integers which coincide with
the dimension of certain moduli spaces of generalized $\theta $-functions
\cite{Beauville}. From the combinatorial nature of the affine plactic algebra we can conclude completely combinatorially that $\langle \hat{\nu},s_{\hat\lambda }(\mathcal{A})\hat{\mu}\rangle \in
\mathbb{Z}$. A purely combinatorial proof of non-negativity seems to be missing. First steps in this direction will appear in a forthcoming paper.}
\end{remark}

\begin{theorem}
\label{mainresult}
Let $\Lambda^{(k)}=\mathbb{Z}[e_{1},\ldots,e_{k}]$ be the ring of symmetric polynomials in $k$ variables. The assignment $\hat{\la}\mapsto s_{\la^t}$ defines an isomorphism of rings
\begin{equation}
\cF=\cF(\widehat{\mathfrak{sl}}(n))_k\cong \mathbb{Z}[e_{1},\ldots,e_{k}]/\langle h_n-1,h_{n+1},\ldots,h_{n+k-1},h_{n+k}+(-1)^{k}e_{k}\rangle
,  \label{viso}
\end{equation}%
where $h_{r}$ denotes the $r^{\text{th}}$ complete symmetric function in $k$ variables.
\end{theorem}

\begin{remark}
{\rm We will see later (Corollary \ref{qcohom}) that $V_{k}(\widehat{\mathfrak{sl}}(n),\mathbb{Z})$ is a quotient of the quantum cohomology ring  of the Grassmannian by factoring out the extra relations $q=e_k$ and $h_n=1$.
}
\end{remark}

\begin{proof}
The notation $\Lambda^{(k)}=\mathbb{Z}[e_{1},\ldots,e_{k}]$ makes sense, since $\Lambda^{(k)}$ is generated by the $e_i$'s and they are algebraically independent (\cite[Prop. 35.1]{Bump}). Let $\Lambda_\mC^{(k)}=\mC\otimes_\mZ\Lambda^{(k)}$ be the complexification. Denote by $I$ be the ideal generated by $f_i$, $0\leq i\leq k$, where $f_0=h_n-1$, $f_k=h_{n+k}+(-1)^{k}e_{k}$ and $f_i=h_{n+i}$ otherwise; and let $R:=\Lambda^{(k)}_\mC/I$ be the quotient, which is of course exactly the complexification of the right hand side of \eqref{viso}.\\
{\it Claim 1: $I$ is a radical ideal, i.e $I=\sqrt{I}$.}\hfill\\
First note that we could also write $\Lambda^{(k)}=\mathbb{Z}[h_{1},\ldots,h_{k}]$ \cite[Prop. 35.1]{Bump}, in particular $e_k\in\mathbb{Z}[h_{1},\ldots,h_{k}]$. Now consider $\Lambda=\varprojlim\Lambda^{(k)}=\mZ[e_1,e_2,\ldots]=\mZ[h_1,h_2,\ldots]$, the ring of symmetric functions with its natural projection $p:\Lambda\surj\Lambda^{(k)}$. Set $\tilde{h}_i=h_i\in\Lambda$ if $i\not\in\{n, n+k\}$, $\tilde{h}_n=h_n-1$ and $\tilde{h}_{n+k}=h_{n+k}+(-1)^ke_k$. For $\mu$ a partition define $\tilde{h}_\mu=\tilde{h}_{\mu_1}\tilde{h}_{\mu_2}\ldots$. The $h_i$, $i\in\{1,2,\ldots, k, n,n+1, \ldots, n+k\}$ are algebraically independent, and therefore also the corresponding $\tilde{h}_i$. Consequently, the $\tilde{h}_\mu$'s form a linearly independent subset of $\Lambda$. Assume now $f\in\Lambda^{(k)}\subset\Lambda$ is not in the ideal $I'$ generated by the $\tilde{h}_i$, $n\leq i\leq n+k$. If we expand $f$ in the $\tilde{h}_\mu$'s there must appear some $\tilde{h}_\mu$ where $\mu_j\not\in\{n,n+1, \ldots, n+k\}$ for all $j$. Let $s$ be a positive integer. Then $f^s$ must contain $\tilde{h}_{\mu^s}$  in its basis expansion, where $\mu^s$ is the partition which contains each part $\mu_j$ exactly $s$ times. In particular $f^s\not\in I'$ for any $s$. Now one can use the projection $p$ to deduce the result.\\
{\it Claim 2: Let $f,g\in\Lambda$. Then $f=g\in R$ if and only if $f(x)-g(x)=0$ for any solution $x$ of the Bethe Ansatz equations \eqref{bae}.}\hfill\\
Using the alternative version \eqref{bae2} of \eqref{bae} (with $z=1$), the claim is a direct consequence of Claim 1 and Hilbert's Nullstellensatz, since then $I$ equals the vanishing ideal of the zero set of $I$.\\

{\it Claim 3: Sending $\hat{\la}$ to the class of $s_{{\la}^t}$ in $R$ defines a ring homomorphism $$F:\quad V_{k}(\widehat{\mathfrak{sl}}(n),\mathbb{C})\rightarrow R.$$}
It  is enough to show that for any $\hat{\la}, \hat{\mu}\in P_{k}^{+}$ we have  $F(\hat{\la}\ast\hat{\mu})=F(\hat{\la})F(\hat{\mu})$ (thanks to \eqref{result}).
Let $b(x)=b(x_\sigma)$ be a (non-zero) Bethe vector. Using \eqref{bethevec} and \eqref{Rietschbetter} for $z=1$ we expand $b(x)=\sum_{\hat{\nu}\in P_{k}^{+}}s_{{\nu^t}^*}(x)\hat\nu$ and obtain from Theorem \ref{eigenbasis}

\begin{eqnarray}
s_{\lambda }(\mathcal{A})b(x) =s_{\lambda ^{t}}(x)b(x)=\sum_{\hat{\mu}\in P^{+}_{k}}s_{\lambda ^{t}}(x)s_{\mu ^{t}}(x)\hat{\mu}^{\ast }
\label{vecexpansion}
=\sum_{\hat{\mu},\hat{\nu}\in P^{+}_{k}}\langle \hat{\mu}^\ast,s_{\lambda }(\mathcal{A})\hat{\nu}^{\ast }\rangle s_{\nu ^{t}}(x)\hat{\mu}%
^{\ast }\;.  \notag
\end{eqnarray}
Comparing coefficients of each basis element in the last two terms we obtain
\begin{equation}
\label{comparecoeffs}
s_{\lambda ^{t}}(x)s_{\mu ^{t}}(x)=\sum_{\hat{\nu}\in P^{+}_{k}}\langle \hat{%
\mu}^{\ast },s_{\lambda }(\mathcal{A})\hat{\nu}^{\ast }\rangle s_{\nu
^{t}}(x)=\sum_{\hat{\nu}\in P^{+}_{k}}\langle \hat{%
\nu},s_{\lambda }(\mathcal{A})\hat{\mu}\rangle s_{\nu
^{t}}(x)
\end{equation}
Note that our convention $s_{(1^n)}(\cA)=e_n(\cA)=1$ (with $z=1$) precisely corresponds under $F$ to the extra imposed relation
 $s_{n}=h_n=1$. To see the equality of the matrix elements recall from Section \ref{autos} that
 $$\langle \hat{%
\mu}^{\ast },s_{\lambda }(\mathcal{A})\hat{\nu}^{\ast }\rangle=\overline{\langle \hat{%
\nu},(\op{flip}\circ s_{\lambda }(\mathcal{A})^\ast\circ\op{flip})\hat{\mu}\rangle},$$
where $s_{\lambda }(\mathcal{A})^\ast$ denotes the right adjoint of $s_{\lambda }(\mathcal{A})$. One easily verifies that $\func{flip}\circ a_{i}\circ \func{flip}=a_{n-i}^{\ast }$ which implies $\op{flip}\circ e_{r }(\mathcal{A})^\ast\circ\op{flip}=e_{r }(\mathcal{A})$ and, thus, $\op{flip}\circ s_{\lambda }(\mathcal{A})^\ast\circ\op{flip}=s_{\lambda }(\mathcal{A})$. Together with the reality of the matrix elements, which is immediate from the action of the affine plactic algebra this gives the desired identity \eqref{comparecoeffs}.

Now, fom the definition of $F$ one obtains
$$A:=F(\hat{\la}\ast\hat{\mu})=F\left(s_{\lambda }(\mathcal{A})\hat{\mu}\right)=
F\left(\sum_{\hat{\nu}\in P_{k}^{+}}\langle \hat{\nu},s_{\lambda }
(\mathcal{A})\hat{\mu}\rangle {\hat\nu}\right)=\sum_{\hat{\nu}\in P_{k}^{+}}\langle \hat{\nu},s_{\lambda }(\mathcal{A})\hat{\mu}\rangle s_{\nu^t},
$$ whereas
$$B:=F(\hat{\la})F(\hat{\mu})=s_{\la^t}s_{\mu^t}.$$

By \eqref{comparecoeffs} we have $A(x)=B(x)$ for any Bethe root $x=x_\sigma$ or any other in the $S_k$-orbit of $x_\sigma$.
Moreover, $A(x)=B(x)$ for any Bethe root where not all the entries are different, since both sides of the equality just vanish. Hence $A(x)=B(x)$ for any solution of the Bethe Ansatz equations. The Claim follows now from Claim 2.\\
{\it Claim 4: $F$ is an isomorphism.}
Since via \eqref{weight2part} the $\la\in\mathfrak{P}_{\leq n-1,k}$ form a basis of the complexification of $V_{k}(\widehat{\mathfrak{sl}}(n),\mathbb{Z})$ it is enough to show that the $s_{{\la}^t}$, ${\la}\in\mathfrak{P}_{\leq n-1,k}$ form a basis of $R/I$. Recall that the the $\tilde{h}_\mu$'s form a basis of $\Lambda^{(k)}$, if $\mu$ runs through the set  of partitions with at most $k$ parts \cite[p 73]{FultonYT}. In particular, the $\tilde{h}_\mu$'s with $\mu\in\mathfrak{P}_{\leq n-1,k}$ span $R/I$. To see the linearly independence assume that $\sum_{\mu\in\mathfrak{P}_{\leq n-1,k}}a_\mu \tilde{h}_\mu=0$ in $R/I$ for some $a_\mu\in\mC$, equivalently $\sum_{\mu\in\mathfrak{P}_{\leq n-1,k}}a_\mu \tilde{h}_\mu=\sum_{i=n}^{n+k}r_i\tilde{h_i}$ for some $r_i\in R$. Expanding this equality in the basis of the $\tilde{h}_\mu$'s we see that the basis vectors on the left hand side are all of the form $\tilde{h}_\mu$, where $\mu\in\mathfrak{P}_{\leq n-1,k}$, whereas the basis vectors on the
right hand side are all of the form $\tilde{h}_\mu$ with $\mu\not\in\mathfrak{P}_{\leq n-1,k}$. In particular, $a_\mu=0$ for all $\mu\in\mathfrak{P}_{\leq n-1,k}$ and so the $\tilde{h}_\mu$ form a basis. Moreover, $\tilde{h}_\mu=h_\mu$ in this case.
Now the transformation matrix between the  $h_\mu$'s and the $s_\mu$'s is triangular with $1$'s on the diagonal \cite[p 75 (4)]{FultonYT}, hence the linearly independence follows.

Since all the constructions are defined over the integers, $F$ induces an isomorphism for the integral version as well.
\end{proof}

\begin{example}
\label{example3and1}
{\rm
Consider the case $n=3$ and $k=1$. The ring $\mC[e_1]/\langle h_4+e_1\rangle$ is isomorphic to $\mC[e_1]/\langle e_1^4-e_1\rangle$
(since $h_4=e_1^4-3e_2e_1^2+2e_3e_1+e_2^2$ evaluated at $x=(x_1,0,0,0,\ldots)$ gives $h_4(x)=e_1^4(x)$), and has basis $1, e_1, e_1^2, e_1^3=h_3$. If we impose the additional relation $h_3=1$, we get a ring on basis $1,e_1, e_1^2$ and multiplication $e_1e_1=e_1^2$, $e_1e_1^2=e_1^3=1$. Under the map $\Phi$, the three partitions $\emptyset$, $(1)$, $(1^2)$ are mapped to $s_\emptyset=1$, $s_{(1)}=h_1=e_1$, $s_{(2)}=h_2=e_1^2$, and the multiplication in the combinatorial fusion ring is precisely the given one.}
\end{example}

\section{Applications}
\subsection{Symmetries of fusion coefficients}

The following (known) symmetry formulae for the fusion coefficients are now easily verified using \eqref{result}:

\begin{proposition}
For $\hat{\la}$, $\hat{\mu}$, $\hat{\nu}\in P_k^+$ set $\mathcal{N}_{\hat\lambda \hat\mu \hat\nu}:=\mathcal{N}_{\hat\lambda \hat\mu
}^{(k),\hat\nu ^{\ast }}$ then

\begin{enumerate}[(i)]
\item \label{one}
$\mathcal{N}_{\pi (\hat\lambda )\pi (\hat\mu )\pi (\hat\nu )}=\mathcal{N%
}_{\hat\lambda \hat\mu \hat\nu }$ for any permutation $\pi \in S_{3}$ of $(\hat\lambda
,\hat\mu ,\hat\nu )$

\item \label{two}
$\mathcal{N}_{\func{rot}(\hat\lambda )\hat\mu \hat\nu }=\mathcal{N}%
_{\hat\lambda \func{rot}(\hat\mu )\hat\nu }=\mathcal{N}_{\hat\lambda \hat\mu \func{rot}%
(\hat\nu )}$

\item \label{three}
$\overline{\mathcal{{N}}_{\hat\lambda \hat\mu \hat\nu }}=\mathcal{N}_{\hat\lambda
^{\ast }\hat\mu ^{\ast }\hat\nu ^{\ast }}=\mathcal{N}_{\hat\lambda \hat\mu \hat\nu }$
\end{enumerate}
\end{proposition}

\begin{proof}
According to \eqref{Verlindeexplicit} and \eqref{permanentformula} (with $\hat\nu$ replaced by $\hat\nu^*$) we have
\begin{eqnarray*}
\mathcal{N}_{\hat\lambda \hat\mu \hat\nu }
&=&\frac{1}{n(k+n)^{n-1}}\sum_{\sigma \in \mathfrak{P}_{\leq n-1,k}}\zeta
^{|\sigma |\frac{|\lambda |+|\mu |+|\nu |}{n}}\frac{s_{\lambda }(\zeta
^{-I(\sigma )})s_{\mu }(\zeta ^{-I(\sigma )})s_{\nu }(\zeta ^{-I(\sigma )})}{%
\prod_{i<j}^{n}|\zeta ^{I_{i}(\sigma )}-\zeta ^{I_{j}(\sigma )}|^{-2}},
\end{eqnarray*}
and the statement \eqref{one} follows. The fusion coefficients do not depend on $z$, therefore we may set $z=1$ and calculate
\begin{eqnarray*}
\mathcal{N}_{\func{rot}(\hat\lambda )\hat\mu \hat\nu }
&\stackrel{\eqref{one}}{=}&\mathcal{N}_{\hat\mu\hat\nu \func{%
rot}(\lambda ) }
\stackrel{\eqref{Verlindeexplicit}}{=}\langle \func{rot}(\hat\lambda )^{\ast },s_{\mu }(\mathcal{A})\hat\nu \rangle
=\langle \func{rot}^{-1}(\hat\lambda ^{\ast }),s_{\mu }(\mathcal{A})\hat\nu
\rangle \\
&\stackrel{(\dagger)}{=}&\langle \hat\lambda ^{\ast },s_{\mu }(\mathcal{A})\func{rot}(\hat\nu )\rangle
\stackrel{\eqref{Verlindeexplicit}}{=}\mathcal{N}_{\hat\mu \func{rot}(\hat\nu )\la}
\stackrel{\eqref{one}}{=}\mathcal{N}_{\hat\lambda \hat\mu
\func{rot}(\hat\nu )},
\end{eqnarray*}
(where ($\dagger$) follows from Corollary \ref{allascomm} using the fact that $\func{rot}=h_{k}(\mathcal{A})_{k}$ when $z=1$) and then apply again part \eqref{one} repeatedly to deduce  $\mathcal{N}_{\func{rot}(\hat\lambda )\hat\mu \hat\nu }=\mathcal{N}_{\la\func{rot}(\hat\mu)\hat\nu}$ as well. To prove \eqref{three} we make once more use of the relation $\op{flip}\circ s_{\la}(\cA)^\ast= s_{\la}(\cA)\circ\op{flip}$ already discussed when deriving \eqref{comparecoeffs}. Together with $s_{\la}(\cA)^\ast=s_{\la^\ast}(\cA)$, which follows from \eqref{Schur spec} and \eqref{Rietschbetter} and the definitions in Section \ref{autos}, we calculate
\begin{eqnarray*}
\mathcal{N}_{\hat\lambda ^{\ast }\hat\mu ^{\ast }\hat\nu ^{\ast }}
=\langle \hat\nu
,s_{\lambda^\ast}(\mathcal{A})\hat\mu ^{\ast }\rangle
=\langle \hat\nu^*,\op{flip}\circ s_{\lambda}(\mathcal{A})^\ast\circ\op{flip}\hat\mu\rangle
=\langle \hat\nu^*,s_{\lambda}(\mathcal{A})\hat\mu\rangle=\mathcal{N}_{\hat\lambda\hat\mu\hat\nu}.
\end{eqnarray*}
\end{proof}

\subsection{Noncommutative Cauchy identities}

Exploiting the combinatorial definition (\ref{result}) of the fusion
product, we now state new identities for the fusion coefficients.

\begin{corollary}[Noncommutative Cauchy identities]
\label{CorCauchy}
Let $\ell\in\{1,\ldots,n-1\}$, $u_i\in\mC$ for $1\leq i\leq \ell$ and $k\in\mZ_{>0}$. Recalling the analogue of Baxter's $Q$-operator from \eqref{Q} and the transfer matrix \eqref{tcoeff} we have the following equalities of endomorphisms of $\mathcal{H}_{k}$:
\begin{eqnarray}
Q(u_{1})_{k}\cdots Q(u_{\ell })_{k}&=&\sum_{\lambda \in \mathfrak{P}_{\leq
\ell ,k}}s_{\lambda }(u_{1},\ldots,u_{\ell})s_{\lambda }(\mathcal{A})_{k},\;
\label{CauchyQ}\\
T(u_{1})_k\cdots T(u_{\ell })_k&=&\sum_{\lambda \in \mathfrak{P}_{\leq n,\ell
}}s_{\lambda ^{t}}(u_{1},\ldots,u_{\ell })s_{\lambda }(\mathcal{A})_k\;.
\label{CauchyT}
\end{eqnarray}%
In particular,
\begin{eqnarray}
h_{\alpha _{1}}(\mathcal{A})_{k}\cdots h_{\alpha _{\ell }}(\mathcal{A}%
)_{k}&=&\sum_{\lambda \in \mathfrak{P}_{\leq \ell ,k}}K_{\lambda \alpha
}s_{\lambda }(\mathcal{A})_{k}\label{Kostka1}\\
e_{\alpha _{1}}(\mathcal{A})_{k}\cdots e_{\alpha _{\ell }}(\mathcal{A}%
)_{k}&=&\sum_{\lambda \in \mathfrak{P}_{\leq n,\ell }}K_{\lambda ^{t}\alpha
}s_{\lambda }(\mathcal{A})_{k},\label{Kostka2}
\end{eqnarray}
with $K_{\lambda \alpha }$ being the Kostka numbers, i.e. the number of
semi-standard tableaux of shape $\lambda $ and weight $\alpha =(\alpha
_{1},\ldots,\alpha _{\ell })\in\mathbb{Z}_{>0}^\ell$.
\end{corollary}

\begin{proof} We employ once more the eigenbasis of Bethe states; see Theorem \ref{eigenbasis}. The formula \eqref{CauchyQ} is a consequence of
the Cauchy-identity $\sum_{\lambda }s_{\lambda }(x)s_{\lambda ^{t}}(y)=\prod_{i,j}(1+x_{i}y_{j})$
for commutative Schur functions. Namely, the definition of
the $Q$-operator in \eqref{Q} together with \eqref{nchspec} gives $Q(u)b(x)=\prod_{j=1}^{k}(1+ux_{j})b(x)$ for any Bethe state $b(x)\in\mathcal{H}_k$. Therefore, it follows that
\begin{eqnarray*}
Q(u_{1})\cdots Q(u_{\ell})b(x)
&=&\prod_{i=1}^{\ell}\prod_{j=1}^{k}(1+u_{i}x_{j})b(x)
=\sum_{\lambda \in \mathfrak{P}_{\leq \ell,k}}s_{\lambda }(u)s_{\lambda
^{t}}(x)b(x)\\
&=&\sum_{\lambda \in \mathfrak{P}_{\leq \ell,k}}s_{\lambda
}(u)s_{\lambda }(\mathcal{A})b(x)\;.
\end{eqnarray*}%
Since the Bethe states form a basis of $\mathcal{H}_{k}$ (Theorem \ref{eigenbasis}), we have proved \eqref{CauchyQ}. To derive \eqref{Kostka1}, recall the well-known expansion
$s_{\lambda }(u)=\tsum_{\alpha }K_{\lambda \alpha }m_{\alpha }(u)$
with $m_{\alpha }$ being the basis of symmetric monomial functions (see e.g. \cite[Chapter I, Section 6, Table 1 and (6.4)]{MacDonald}). Expanding both sides of \eqref{CauchyQ} with respect to the basis of symmetric monomial functions, we obtain the asserted formula \eqref{Kostka1}.\\
The proof of \eqref{CauchyT} and \eqref{Kostka2} is completely analogous to the known one for the ring of commutative symmetric functions, since we have shown in Corollary \ref{allascomm} that the cyclic noncommutative Schur polynomials satisfy all the familiar relations of their commutative counterparts.
\end{proof}

As an application of the generalised Cauchy identities we now derive novel identities for the fusion coefficients. Let $\ell$ and $\alpha$ be as in Corollary \ref{CorCauchy} and $\mu, \nu\in\mathfrak{P}_{\leq n-1,\ell}$. Then denote by $\hat\mu,\hat\nu$ and $\widehat{(\alpha_{i})},\,\widehat{(1^{\alpha _{i}})}$ the affine weights in $P^+_k$ corresponding to the pre-images of $\mu,\nu$ and a horizontal and vertical strip of length $\alpha_i$ under the bijection \eqref{weight2part}.

\begin{corollary}[Kostka numbers and fusion coefficients]\label{Kostkaid}
Given a  sequence $S=(\nu=\mu^{(0)}, \mu^{(1)}, \ldots, \mu^{(\ell)}=\mu)$ of partitions in $\mathfrak{P}_{\leq n-1,\ell}$ set $\cN_i(S)=\mathcal{N}_{\widehat{(\alpha_{i})}\,\widehat{\mu^{(i)}}}^{(k)\widehat{\mu^{(i-1)}}}$ and $\cN_i'(S)=\mathcal{N}_{\widehat{(1^{\alpha _{i}})}\,\widehat{\mu^{(i)}}}^{(k)\widehat{\mu^{(i-1)}}}$. Then we have the identities
\begin{eqnarray}
\sum_{S}\cN_1(S)\cN_2(S)\ldots \cN_\ell(S)&=&\sum_{\lambda \in \mathfrak{P}_{\leq n-1,\ell
}}K_{\lambda \alpha }\mathcal{N}_{\hat\lambda \hat\mu }^{(k)\hat\nu }  \label{Kostkaold1}\\
\sum_{S}\cN'_1(S)\cN'_2(S)\ldots \cN'_\ell(S)&=&\sum_{\lambda \in \mathfrak{P}_{\leq
n,\ell }}K_{\lambda ^{t}\alpha }\mathcal{N}_{\hat\lambda \hat\mu }^{(k)\hat\nu }\;.\label{Kostkaold2}
\end{eqnarray}
where the left sums run through the set of all sequences $S$ as above and $\hat\lambda\in P^+_k$ is the pre-image of \eqref{weight2part}. If $\lambda\in\mathfrak{P}_{\leq n,\ell }$ we first delete all columns of length $n$.
\end{corollary}

\begin{proof} Consider the matrix element $\langle \hat\nu, h_{\alpha _{1}}(\mathcal{A}
)\cdots h_{\alpha _{\ell }}(\mathcal{A})\hat\mu \rangle$ with $\hat\mu,\hat\nu\in P^+_k$. Employing \eqref{Kostka1} we obviously have
$$\langle \hat\nu, h_{\alpha _{1}}(\mathcal{A}
)\cdots h_{\alpha _{\ell }}(\mathcal{A})\hat\mu \rangle=\sum_{\lambda \in \mathfrak{P}_{\leq n-1,\ell
}}K_{\lambda \alpha }\langle\hat\nu,s_\la(\cA)\hat\mu\rangle\;$$

and \eqref{Kostka1} now follows from \eqref{Verlindeexplicit}. The second equality is proved along the same lines by considering the matrix element $\langle\hat
\nu, e_{\alpha _{1}}(\mathcal{A})\cdots e_{\alpha _{\ell }}(\mathcal{A}
)\hat\mu \rangle$ and applying the identity \eqref{Kostka2}.
\end{proof}

Using the phase algebra generators we now derive a relation between fusion coefficients of level $k$ and level $k+1$. Such a recursion relation seems to be new and is a direct consequence of our combinatorial particle description of the fusion ring.

\begin{corollary}[Recursion relation for fusion coefficients]\label{recfusion}
Adopt the same conventions as in the previous corollary. Choose any integer vector ${\bf j}=(j_{1},\ldots,j_{\ell })\in\mathbb{Z}^\ell_n$ and define for level $k+1$ the fusion coefficients $\cN_i'(S,{\bf j})=\mathcal{N}_{\widehat{(1^{\alpha _{i}})}\,\varphi_{j_{i}^*}\widehat{\mu^{(i)}}}^{(k+1)\,\varphi_{j_{i}^*}\widehat{\mu^{(i-1)}}}$, then we have the recursion relation
\begin{equation}\label{fusionrec}
\sum_{S}\cN'_1(S,{\bf j})\cN'_2(S,{\bf j})\ldots \cN'_\ell(S,{\bf j})=\sum_{\lambda
\in \mathfrak{P}_{\leq n,\ell }}K_{\lambda ^{t}\alpha }\mathcal{N}_{\hat\lambda
\hat\mu }^{(k)\hat\nu }\;.
\end{equation}
In particular, for $\alpha=(r,0,\ldots,0)$ we simply have
\begin{equation}
\mathcal{N}_{\widehat{(1^{r})}\,\varphi_{j}^\ast\hat\mu}^{(k+1)\,\varphi_{j}^\ast\hat\nu}= \mathcal{N}_{\widehat{(1^{r})}\,\hat\mu}^{(k)\,\hat\nu}
\end{equation}
for any $j\in\mathbb{Z}_n$.
\end{corollary}

\begin{proof}
Employing Lemma \ref{ABCD} one easily verifies the following
relations between the generators of the Yang-Baxter and the phase algebra,%
\begin{equation*}
A(u)\varphi _{1} =\varphi _{1}A(u),\quad \varphi _{1}D(u) =u C(u),\quad D(u)=uC(u)\varphi _{1}^{\ast }\;.
\end{equation*}%
From the latter identities one now obtains for the transfer matrix \eqref{tcoeff},
\begin{equation*}
\varphi _{1}T(u)\varphi _{1}^{\ast }=\varphi_1 A(u)\varphi_1^\ast+z \varphi_1 D(u)\varphi_1^\ast
=A(u)\varphi_1\varphi_1^\ast+z u C(u)\varphi_1^\ast=T(u).
\end{equation*}
But since for $z=1$ one has $\func{rot}\circ T(u)=T(u)\circ \func{rot}$, this is obvious as each $e_{r}(%
\mathcal{A})$ is invariant under a cyclic rotation when $z=1$, we must have
\begin{equation*}
\varphi _{i}T(u)\varphi _{i}^{\ast }=T(u)
\end{equation*}%
for all $i=0,1,\ldots,n-1$. Therefore, we have $\varphi _{i}e_{r}(\mathcal{A})\varphi _{i}^{\ast }=e_{r}(\mathcal{A})$ and the assertion follows from (\ref{Kostka2}),
\begin{eqnarray*}
\langle\hat\nu, e_{\alpha _{1}}(\mathcal{A})\cdots e_{\alpha _{\ell }}(\mathcal{A}
)\hat\mu \rangle&=&\langle\hat
\nu, \varphi_{j_1} e_{\alpha _{1}}(\mathcal{A})\varphi_{j_1}^\ast\cdots \varphi_{j_\ell}e_{\alpha _{\ell }}(\mathcal{A}
)\varphi_{j_\ell}^\ast\hat\mu \rangle\\
&=&\sum_{\lambda \in \mathfrak{P}_{\leq
n,\ell }}K_{\lambda ^{t}\alpha }\mathcal{N}_{\hat\lambda \hat\mu }^{(k)\hat\nu }\;.
\end{eqnarray*}
\end{proof}

\begin{example}
\label{ExFusion}
{\rm
Let $n=3$ and $k=2$. Choose $\ell=k=2$ and set $\mu =(2,1)$. Then it is not difficult to compute, using the combinatorial formula \eqref{result}, the fusion product
$$ \Yvcentermath1\yng(2,1)\ast\yng(2,1)=
\Yvcentermath1\yng(2,1)+ \emptyset \;.$$

Fixing $\nu=\emptyset$ and making the choices $\alpha=(1,1,1),\;{\bf j}=(2,3,1)$ one finds that on the left hand side of \eqref{fusionrec} only two sequences $S$ contribute, namely
$$\Yvcentermath1 \left(\nu=\emptyset,\yng(1,1)\;,\,\yng(1)\;,\,\yng(2,1)\;=\mu\right)\;\text{and}
\; \left(\nu=\emptyset,\yng(1,1)\;,\,\yng(2,2)\;,\,\yng(2,1)\;=\mu\right)\;.$$
Again one computes combinatorially
$$\mathcal{N}_{\widehat{(1)}\,\varphi_{2}^{\ast}\widehat{(1,1)}}^{(k+1)\,
\varphi_{2}^\ast\hat\emptyset}=
\mathcal{N}_{\widehat{(1)}\,\varphi_{3}^{\ast}\widehat{(1,0)}}^{(k+1)\,
\varphi_{3}^\ast\widehat{(1,1)}}=
\mathcal{N}_{\widehat{(1)}\,\varphi_{3}^{\ast}\widehat{(2,2)}}^{(k+1)\,
\varphi_{3}^\ast\widehat{(1,1)}}=
\mathcal{N}_{\widehat{(1)}\,\varphi_{1}^{\ast}\hat{\mu}}^{(k+1)\,
\varphi_{1}^\ast\widehat{(1,0)}}=
\mathcal{N}_{\widehat{(1)}\,\varphi_{1}^{\ast}\hat{\mu}}^{(k+1)\,
\varphi_{1}^\ast\widehat{(2,2)}}=1$$
leading to the identity
$$\sum_{S}\cN'_1(S,{\bf j})\cN'_2(S,{\bf j})\cN'_3(S,{\bf j})=2\;.$$
The only possible semi-standard tableaux of weight $\alpha$ lying within the bounding box of size $3\times 2$ are
$\tiny\Yvcentermath1\young(1,2,3)$ and $\tiny\Yvcentermath1\young(12,3)\,,\;\young(13,2)$. %
Obviously, $\mathcal{N}_{\hat{\emptyset}\hat\mu }^{(k)\hat\emptyset }=0$ and, thus, we compute for the right hand side of \eqref{fusionrec}
$$\sum_{\lambda
\in \mathfrak{P}_{\leq n,k }}K_{\lambda ^{t}\alpha }\mathcal{N}_{\hat\lambda
\hat\mu }^{(k)\hat\nu }= 2 \mathcal{N}_{\widehat{(2,1)}
\hat\mu }^{(k)\hat\emptyset }=2$$
which yields the desired equality.
}
\end{example}

\section*{Part II: The quantum cohomology ring: Fermions on a circle}
\section{Free fermion formulation of the small quantum cohomology ring}

\label{sec:qcoho} We want to develop structures parallel to the ones discussed before for the fusion ring, but now for the quantum cohomology
ring and finally connect the two sides. Building on earlier results of for instance Bertram \cite{Bertram}, Rietsch \cite{Rietsch}, and Postnikov \cite{Postnikov} we present a free fermion formulation of the small quantum cohomology ring $qH^{\bullet }(\op{Gr}_{k,k+n})$\
of the complex Grassmannian. Here $k$ and $n$ are analogous to the level and rank in the context of the fusion ring.

 More precisely, let $n, k\in\mZ_+$ and denote $N=k+n$. Let $\op{Gr}_{k,n+k}$ be the Grassmannian of $k$-planes inside $\mC^N=\mC^{n+k}$. The integral cohomology ring $H^\bullet(\op{Gr}_{k,k+n})$ can explicitly be described using Schubert calculus (see for example \cite[9.4]{FultonYT}): the Schubert classes $[\Omega_\la]$ form a basis, the structure constants are the intersection numbers, which coincide with the Littlewood-Richardson coefficients. (Here $\la$ runs through $\mathfrak{P}_{\leq k,n}$.)

 The small quantum cohomology $qH^\bullet(\op{Gr}_{k,n+k})$ is a deformation of the usual (singular) cohomology. More precisely, it is a $\mZ[q]$-algebra which is isomorphic to $\mZ[q]\otimes_\mZ H^\bullet(\op{Gr}_{k,n+k})$ as a $\mZ[q]$-module. The Schubert classes give a $\mZ[q]$-basis $\{1\otimes [\Omega_\la]\}$. This module can be equipped (in a non-obvious way) with a ring structure where the structure constants $C_{\la,\mu}^{\nu}(q)=\sum C_{\la,\mu}^{\nu,d}q^d$ are given by the so-called 3-point, genus 0  Gromov-Witten invariants $C_{\la,\mu}^{d,\nu^\vee}$ which count the number of rational curves of degree $d$ passing through generic translates of $\Omega_\la$, $\Omega_\mu$, $\Omega_{\nu}$. (In the cases $|\la|+|\mu|+|\nu|\not=kn+d(k+n)$, where the number of curves could be infinite, one just puts $C_{\la,\mu}^{d,\nu}=0$.) For a general overview we refer to \cite{FP} and \cite{Kock}. Note that for $q=0$ we obtain the ordinary cohomology $H^\bullet(\op{Gr}_{k,n+k})$, with the ring structure given by  the cup or intersection product. For a description of the quantum cohomology ring from the Schubert calculus point of view we refer to \cite{Bertram}, \cite{Buch}, \cite{BKT}, \cite{Tamvakis}.

In the rest of the paper we will develop a combinatorial (fermionic model) for the quantum cohomology similar to our combinatorial description of the fusion algebra.

\subsection{01-words, partitions and Schubert cells}\label{prelim}
Throughout our discussion we will use the following well-known bijections between
01-words of length $N$ and partitions (see e.g. \cite{Postnikov}) as well as the basis vectors of the wedge space $\bigwedge \mC^N$.

\begin{definition}{\rm
Let $w=w_{1}w_2\cdots w_{N}$ be a word with $N$ letters $w_i$ which
are either $0$ or $1$. The {\it weight} of $w$ is the sum $\sum_i w_i=k$. Let
\begin{equation*}
W_{k,N}=\left\{ w=w_{1}w_2\cdots w_{N}~\left\vert ~|w|:=\tsum_{i}w_{i}\right.
=k\right\}
\end{equation*}
be the set of words of weight $k$. Let $\op{n}_{i}(w):=w_{1}+\cdots
+w_{i}$ the number of $1$-letters $w_j$ for $j$ in the closed interval $[0,i]$.}
\end{definition}

For $w\in W_{k,N}$ denote by $1\leq\ell _{1}<\ldots<\ell _{k}\leq N$ the positions of the $1$-letters in $w$ counting from the left. Then there is obviously a bijection $
\mathfrak{P}_{\leq k,n}\rightarrow W_{k,N}$ given by
\begin{equation*}
\lambda \mapsto w(\lambda )=0\cdots 0\underset{\ell _{1}}{1}0\cdots 0%
\underset{\ell _{k}}{1}0\cdots 0,
\end{equation*}%
where the positions $\ell _{i}$ of the letters 1 are determined according to
the formula
\begin{equation}
\ell _{i}=\lambda _{k+1-i}+i\;.  \label{ppositions}
\end{equation}
We shall denote the inverse image under this bijection $w\mapsto \lambda
(w)$. Graphically this correspondence is simply assigning to each $\lambda $ the
path which is traced out by its Young diagram in the $k\times n$ rectangle.
Namely, starting from the left bottom corner in the $k\times n$ rectangle
each move by a box to the right corresponds to a letter $0$ and each move by
a box up to a letter 1. For instance, the partitions displayed by the Young diagrams from \eqref{PYoung} correspond under this bijection to the set of words
$$\{11000,10100,10010,10001,01100,01010,00110,010001,00101,00011\}.$$

Recall the well-known correspondence (see e.g. \cite{FultonYT}, \cite{Postnikov}) between
01-words of weight $k$ and Schubert classes in $qH^{\bullet }(\limfunc{Gr}%
\nolimits_{k,N})$: fix a full standard flag $\mathbb{C}\subset \mathbb{C}^{2}\subset \cdots \subset \mathbb{C}^{N}$
and assign to each word $w$ the Schubert cell
\begin{equation}
{\Omega}_{w}^\circ=\{V\in \limfunc{Gr}\nolimits_{k,N}~|~\dim (V\cap \mathbb{C}%
^{l})=w_{N}+w_{N-1}\cdots +w_{N+1-l}\}\;.
\end{equation}%
The associated Schubert class $[\Omega _{w}]$ is the fundamental
cohomology class of the Schubert variety $\Omega _{w}$ (i.e. the closure
of the Schubert cell.)

Siebert and Tian \cite{ST} gave an explicit presentation of $qH^\bullet(\op{Gr}_{k,N})$ in terms of the ring of symmetric polynomials:
\begin{equation}\label{SiebertTian}
qH^{\bullet }(\op{Gr}_{k,N})\cong \mathbb{Z}[q]\otimes_\mZ\mZ[e_1,e_2,\ldots e_k]/
\left\langle h_{n+1},\ldots,h_{N-1},h_{N}+(-1)^{k}q\right\rangle \;
\end{equation}
where the $e_i$'s are the (commutative) elementary symmetric functions in $k$ variables and the $h_i$'s are the corresponding complete symmetric functions. A $\mZ[q]$-basis is given by the images of the Schur polynomials $s_\la$ for $\la\in\mathfrak{P}_{\leq k,n}$. The following is then a direct consequence of Theorem \ref{mainresult}:
\begin{corollary}[Verlinde algebra as a quotient of $qH^\bullet$]\hfill\\
\label{qcohom}
The Verlinde algebra $V_{k}(\widehat{\mathfrak{sl}}(n),\mathbb{Z})$ is isomorphic to the quotient of the quantum cohomology ring $qH^{\bullet }(\op{Gr}_{k,k+n})$ by imposing the extra relations $q=e_k$ and $h_n=1$.
\end{corollary}

Finally the words of weight $k$ (or the partitions from $\mathfrak{P}_{\leq k,n}$) are also in bijection to basis vectors of $\bigwedge\mC^N$ by sending a partition $\la$ to the vector $v_{\ell_k}\wedge\cdots\wedge v_{\ell_1}$ with $v_1,\ldots,v_N$ the standard basis vectors of $\mC^N$. (This will turn out to be a small `shadow' of the so-called Fock space which appears in the Fermion-Boson correspondence \cite{Kac}, \cite{KacRaina}; see the table in Part III.)

\subsection{State space and fermions}

Paralleling our previous discussion for the fusion ring we now introduce a vector space of
states for $qH^{\bullet }(\limfunc{Gr}\nolimits_{k,N})$ which will again lead to a
description of the ring structure in terms of quantum particles hopping on circle but now with $N$ sites.
Choose $N=n+k\in \mathbb{Z}_{>0}$ to be fixed and define the vector space
\begin{equation}
\mathcal{F}=\tbigoplus_{k=0}^{N}\mathcal{F}_{k},\qquad \mathcal{F}_{k}=%
\mathbb{C}\mathfrak{P}_{\leq k,N-k}=\mathbb{C}W_{k,N}\label{fermionspace}
\end{equation}%
with $\mathcal{F}_{0}=\mathbb{C}\{\emptyset \}$ as before. If we denote the basis vector in  $\mathcal{F}_k$ given by the empty partition by $\emptyset=\emptyset_k$, then $\emptyset_0$ has the physical interpretation of the vacuum, i.e. no particles in the system.  Note that in contrast to the fusion ring we are now considering a finite-dimensional vector space.

We want to construct an analogue of the phase algebra for our vector space $\mathcal{F}
$.
For $1\leq i\leq N$ define the following linear endomorphisms $\psi _{i}^{\ast }$, $\psi _{i}$ of $\mathcal{F}$:
\begin{eqnarray*}
\psi _{i}^{\ast }(w)&:=&
\begin{cases}
(-1)^{\op{n}_{i-1}(w)}w', & \text{if $w_i=0$. Here $w'$ differs from $w$ only by $w_i'=1$}\\
0 & \text{if $w_i\not=0$}
\end{cases}\\
\psi _{i}(w)&:=&
\begin{cases}
(-1)^{\op{n} _{i-1}(w)}w', & \text{if $w_i=1$. Here $w'$ differs from $w$ only by $w_i'=0$}\\
0 & \text{if $w_i\not=1$}.
\end{cases}
\end{eqnarray*}

The associated Young diagram $\psi _{i}^{\ast }\lambda $ is obtained by
adding the top row of the Young diagram $\lambda $ to itself (thereby
increasing its height) and then subtracting a boundary ribbon starting in
the $(i-k-1)$-diagonal and ending in the top row. In contrast $\psi
_{i}\lambda $ is the Young diagram which is obtained from $\lambda $ by adding a
boundary ribbon with same start and end points as in the previous case
 and then subtracting the top row (thereby decreasing the height) - it corresponds to reading the picture below in reverse order.

\begin{example}{\rm
To visualize the action of $\psi_i^\ast$ consider the special case $n=k=4$ and $\mu=(4,3,3,1)$: $\psi_3^\ast\mu$ is depicted in the figure below, where the entries in the diagram label the diagonals. The $(3-k-1)=-2$-diagonal determines the start of the boundary ribbon (the shaded boxes) which has to be subtracted:
\begin{equation*}
\includegraphics[scale=1.28]{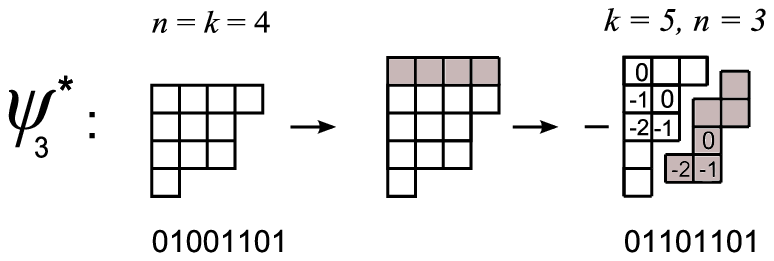}
\end{equation*}
}
\end{example}

\begin{proposition}[Fermion/Clifford algebra]\label{Clifford}
\begin{enumerate}
\item The endomorphisms $\psi _{i},\psi _{i}^{\ast }$ create a subalgebra $\mathfrak{C}\subset \limfunc{End}\mathcal{F}$ isomorphic to a $2^{2N}$-dimensional Clifford algebra with the following defining relations:
\begin{equation*}
\psi _{i}\psi _{j}+\psi _{j}\psi _{i}=\psi _{i}^{\ast }\psi _{j}^{\ast
}+\psi _{j}^{\ast }\psi _{i}^{\ast }=0,\quad \psi _{i}\psi _{j}^{\ast }+\psi
_{j}^{\ast }\psi _{i}=\delta _{ij}\;.
\end{equation*}%
\item If we introduce the following scalar product on $\mathcal{F},$
$\langle \alpha w,\beta w'\rangle =\overline{\alpha}\beta\delta _{w,w'}=\tprod_{i=1}^{N}\delta _{w_{i},w_{i}'}~,$
then $\langle \psi_{i}\lambda ,\lambda ^{\prime }\rangle =\langle \lambda ,\psi _{i}^{\ast
}\lambda ^{\prime }\rangle$.
\item If we introduce a grading on $C$ by putting $\psi _{i}$ in degree $-1$ and $\psi _{i}^{\ast }$ in degree $1$, we get a triangular decomposition $C=C_-\oplus C_0\oplus C_+$. Then $C_-$ annihilates the vacuum vector $\emptyset_0$, and $\cF$ is a free $C_+$-module of rank $1$, with graded decomposition $\cF=\oplus_{k}\mathcal{F}_k$. In particular,  $\mathcal{F}$ becomes an irreducible lowest weight Verma module module for the Clifford algebra $C$.
\end{enumerate}
\end{proposition}

\begin{proof}
All assertions can be easily verified by straightforward calculations.
\end{proof}

If we consider a word $w$ as an $N$-tuple of integers, similar to the phase algebra the maps $\psi _{i}^{\ast },\psi _{i}$ have the
interpretation of particle creation and annihilation operators at site $i$
of a circular lattice but this time of size $N=n+k$. Moreover, they now
describe free fermions which obey the Pauli exclusion principle, i.e. only
one particle per site is allowed. The dimension (level) $k$ is again the
total number of particles which is given by the operator $K=\sum_{i=1}^{N}K_{i},K_{i}=\psi _{i}^{\ast }\psi _{i}$.
The codimension $n=N-k$ is the number of unoccupied sites, i.e. holes in the
system. Again we scalar extend our phase space \eqref{fermionspace}) to
\begin{equation}
\mathcal{F}[q]=\mathbb{C}[q]\otimes _{\mathbb{C}}\mathcal{F}\, ,
\end{equation}
because later on we introduce quasi-periodic boundary conditions on the circle
$\psi_{N+1}^{\ast }=q\psi _{1}^{\ast }$ and $\psi _{N+1}=q^{-1}\psi _{1}$.

Note that our interpretation of $\mathcal{F}[q]$ as state space of a {\it quantum mechanical} particle system is also valid when setting $q=0$. In particular setting $q=0$ does not correspond to the classical limit of our particle system (in the physical sense) but simply alters the boundary conditions. The adjective {\it quantum} in {\it quantum cohomology} is in our context simply synonymous with {\it $q$-deformed} and has nothing to do with its physical meaning.
\subsection{Discrete symmetries: parity and time reversal, particle-hole
duality, shift}\label{trans}

We now introduce a number of mappings which arise naturally in the context
of the physical particle picture. Later on we relate them to symmetries of Gromov-Witten invariants.

\begin{description}
\item[Parity Reversal/Poincare duality] We introduce the parity reversal operator $\mathcal{P}
$ (an involution) which reverses the order of letters in a word by setting,%
\begin{equation*}
\mathcal{P}w=w^{\vee },\quad \quad w^{\vee }:=w_{N}w_{N-1}\cdots w_{1}\;.
\end{equation*}%
This operator corresponds to Poincar\'{e} duality, since $w(\lambda ^{\vee
})=w(\lambda )^{\vee }$. The transformation property of the fermion
operators under parity reversal is (with the notation $e^{\iota\pi K}(?)=(-1)^{{K}(?)}$)
\begin{equation}
\mathcal{P}\psi _{i}\mathcal{P}=e^{\iota\pi K}\psi _{N+1-i},\qquad \mathcal{P}%
\psi _{i}^{\ast }\mathcal{P}=e^{\iota\pi K}\psi _{N+1-i}^{\ast }\;.
\end{equation}

\item[Time reversal/complex conjugation] In continuum quantum physics time reversal is related
to complex conjugation of the wave function through the Schr\"{o}dinger
equation. In the present context of a discrete system we generalize this
notion in terms of the antilinear involution
$$\mathcal{T}\tsum_{w}\alpha _{w}w=\tsum_{w}\bar{\alpha}_{w}w$$
with $\overline{q}=q^{-1}$. In particular,
\begin{equation}
\mathcal{T}\psi _{i}\mathcal{T}=\psi _{i},\qquad \mathcal{T}\psi _{i}^{\ast }%
\mathcal{T}=\psi _{i}^{\ast }.
\end{equation}

\item[Particle-Hole Duality/transpose] This involution, which we shall denote by the
symbol $\mathcal{C}$, is associated with the bijection $\mathfrak{P}_{\leq
k,n}\rightarrow \mathfrak{P}_{\leq n,k}$ or to the  swapping of the letters in each
word $w$: $0$-letters become $1$-letters and vice versa. In terms of partitions we define
\begin{equation}
\mathcal{C}\lambda =(\lambda ^{\vee })^{t}=(\lambda ^{t})^{\vee },
\end{equation}%
i.e. we take the transpose and the Poincar\'{e} dual.
The transformation properties for the Clifford algebra generators are%
\begin{equation}
\mathcal{C}\psi _{i}\mathcal{C}=(-1)^{i-1}\psi _{i}^{\ast },\qquad \mathcal{C%
}\psi _{i}^{\ast }\mathcal{C}=(-1)^{i-1}\psi _{i}\;.
\end{equation}
Together with parity reversal this gives the known isomorphism of rings
\begin{equation}
\omega :qH^{\bullet }(\limfunc{Gr}\nolimits_{k,N})\rightarrow qH^{\bullet }(%
\limfunc{Gr}\nolimits_{n,N}),\qquad 1\otimes[\Omega _{\lambda }]\mapsto 1\otimes[\Omega _{\lambda^t }]
\;,  \label{duality2}
\end{equation}%
which in the language of symmetric functions just amounts to swapping the bases of
elementary and complete symmetric functions, $\omega :e_{r}\mapsto h_{r}$.

\item[Shift operator] There is an analogue $\op{Rot}$ of the cyclic shift operator $\op{rot}$ for the fusion ring and defined on words as follows
\begin{equation}
\limfunc{Rot}(w)=w_{2}\cdots w_{N}w_{1}\;.
\end{equation}%
This operator appears already in the context of Gromov-Witten invariants in \cite{Postnikov}%
. In terms of Young diagrams its action corresponds to shifting the entire
diagram within the $k\times n$ bounding box by one column to the left if $%
w_{1}=0$ and by one row down if $w_{1}=1$. Expressed in formulae this means,%
\begin{equation}
\limfunc{Rot}(\lambda )=\left\{
\begin{array}{cc}
(\lambda _{1}-1,\ldots,\lambda _{k}-1), & \text{if }w_{1}(\lambda )=0 \\
(n,\lambda _{1},\lambda _{2},\ldots,\lambda _{k-1}), & \text{if }w_{1}(\lambda
)=1%
\end{array}%
\right. \;.  \label{RotN}
\end{equation}
Under the identification of partitions with particle configurations on a circle via \eqref{weight2part} it simply rotates the particle positions on the $N$-circle, alias the Dynkin diagram, but this time of $\widehat{\mathfrak{sl}}(N)$.
\end{description}

\section{The affine nil-Temperley-Lieb algebra and noncommutative
polynomials}

Mimicking the passing from the phase algebra to the local plactic algebra we define the following weight preserving endomorphism of $\mathcal{F}$:
\begin{eqnarray}\label{nTLrep}
u_{i}=\psi _{i+1}^{\ast }\psi _{i}~,\quad\text{if } 1\leq i\leq N-1,&& u_{N}\mapsto -qe^{\iota\pi \op{n}_N}\psi _{1}^{\ast }\psi _{N}\;.
\end{eqnarray}

\begin{proposition}
Let $N\geq 2$.
\begin{enumerate}
\item The subalgebra $E_{\op{fin}}$ of $\END(\mathcal{F})$ generated by the $u_i$'s $1\leq i\leq N-1$ is isomorphic to the nil-Temperley-Lieb algebra $0$-${TL}_{N}$.
\item The subalgebra $E_{\op{aff}}$ of $\END(\mathcal{F}[q])$ generated by the $u_i$'s $1\leq i\leq N$ is isomorphic to the affine nil-Temperley-Lieb algebra $0$-$\widehat{TL}_{N}$.
\end{enumerate}
\end{proposition}
\begin{proof} Recall (see e.g. \cite{FG}, \cite{Postnikov}) that $0$-$\widehat{TL}_{N}$ is the unital associative algebra with generators  $\{\hat{u}_{1},\ldots,\hat{u}_{N}\}$ and relations
\begin{equation}\label{nTLrel}
\hat{u}_{i}^{2}=\hat{u}_{i}\hat{u}_{i+1}\hat{u}_{i}=\hat{u}_{i+1}\hat{u}_{i}\hat{u}_{i+1}=0,\qquad
\hat{u}_{i}\hat{u}_{j}=\hat{u}_{j}\hat{u}_{i}\text{\quad if }i-j\neq \pm 1\func{mod}N
\end{equation}%
where all indices are understood modulo $N$. The subalgebra generated by the $\hat{u}_i$, $i>0$ is the ordinary nil-Temperley-Lieb algebra $0$-${TL}_{N}$. Using the Clifford relations it is straight-forward to check that the $u_i$'s satisfy the nil-Temperley-Lieb relations, hence $\hat{u}_i\mapsto u_i$ defines surjective algebra morphisms $0$-${TL}_{N}\surj E_{\op{fin}}$  and $0$-$\widehat{TL}_{N}\surj E_{\op{aff}}$. By \cite[Proposition 2.4.1]{BFZ}, $\mathcal{F}$ is a faithful representation of $0$-${TL}_{N}$ from which one easily deduces that $\mathcal{F}[q]$ is a faithful representation of $0$-$\widehat{TL}_{N}$.
\end{proof}

In the following we will consider the nil-Temperley algebras as algebras of endomorphisms and will not distinguish in the notation between $u_i$ and $\hat{u}_i$.

\begin{remark}[Fermion hopping]{\rm
In terms of our particle description the operators $u_{i}$ simply correspond
to one fermion hopping in clockwise direction on the circle. The action of the affine nil-Temperley-Lieb algebra leaves the $k$%
-particle space invariant, i.e. $u_{i}\mathcal{F}_{k}[q]\subset \mathcal{F}%
_{k}[q] $. In terms of Young diagrams $u_{i}$ adds a box in the $(i-k)^{%
\text{th}}$-diagonal for $i=1,\ldots,N-1$. (The $p^{\text{th}}$%
-diagonal of a Young diagram is formed by all squares $(i,j)$ with content $%
p=j-i$). The action of the affine generator $u_{N}$ corresponds to removing
(adding) a rim hook of size $N-1$. If either of these actions is not allowed
the image is $0$; see \cite{Postnikov} for details.}
\end{remark}

\begin{lemma}
We have the following transformation properties of the affine nil-Temperley-Lieb generators, $1\leq i\leq N-1$, where $u_{i}^*=\psi _{i}^{\ast }\psi _{i+1}$ denotes the right adjoint of $u_i$:
\begin{eqnarray}
\begin{array}{lllllll}
\mathcal{P}u_{N}\mathcal{P}&=&-qe^{\iota\pi\op{n}_N}\psi _{N}^{\ast }\psi _{1},
&&\mathcal{P}u_{i}\mathcal{P}&=&u_{N-i}^{\ast }\label{PTL} \\
\mathcal{T}u_{N}\mathcal{T}&=&-q^{-1}e^{\iota\pi\op{n}_N}\psi _{1}^{\ast }\psi _{N},
&&\mathcal{T}u_{i}\mathcal{T}&=&u_{i},
\label{TTL}\\
\mathcal{C}u_{N}\mathcal{C}&=&-qe^{\iota\pi\op{n}_N}\psi _{N}^{\ast }\psi _{1},&&
\mathcal{C}u_{i}\mathcal{C}&=&u_{i}^{\ast }.\label{CTL}
\end{array}
\end{eqnarray}
\end{lemma}

Analogously to Section \ref{ncpoly1}, we again consider noncommutative symmetric polynomials:

\begin{definition}
\label{defehs}
For $r=1,2,\ldots,N-1$ let
\begin{equation}
\label{ncehdefTL}
e_{r}(\mathcal{U})=\sum_{|I|=r}\tprod_{i\in I}^{\circlearrowright}u_{i}, \quad h_{r}(\mathcal{U})=\sum_{|I|=r}\tprod_{i\in
I}^{\circlearrowleft}u_{i}.
\end{equation}
We also set $e_{N}(\mathcal{U})=-qe^{\iota\pi\op{n}_N}$ and $h_{N}(\mathcal{U})_{\vert\mathcal{F}_{N}[q]}=q$. Note that we have $h_{N}(\mathcal{U})_{|\mathcal{F}_{k}[q]}=0$ for $k<N$ and $h_{r}(\mathcal{U})=0$ for $r>N$.
For $\lambda\in\mathfrak{P}_{\leq N,N}$ define the noncommutative Schur polynomial
\begin{equation}
\label{defsNil}
s_{\lambda }(\mathcal{U})\stackrel{(1)}{=}\det (e_{\lambda _{i}^{t}-i+j}(\mathcal{U}
))_{1\leq i,j\leq N}\stackrel{(2)}{=}\det (h_{\lambda _{i}-i+j}(\mathcal{U}))_{1\leq i,j\leq
N}\;.
\end{equation}
\end{definition}
\begin{remark}{\rm
Note that in comparison to \eqref{ncehdef} the cyclic ordering is reversed and that the complete symmetric polynomials are now square-free because of \eqref{nTLrel}.}
\end{remark}
\begin{remark}{\rm
The definition of $s_{\lambda }(\mathcal{U})$ only makes sense after we have proved that the $e_r(\cU)$ resp. $h_r(\cU)$ pairwise commute (see Corollary \ref{Nilecomm} below). We also have to show that the two definitions from \eqref{defsNil} agree (Theorem \ref{freeeigenbasis}). However, we would like to accept the definitions and draw a few consequences first. We also want to point out that both, the definitions \eqref{ncehdefTL} as well as Corollary \ref{Nilecomm} are due to Postnikov \cite{Postnikov}. Our approach provides a new point of view of these results.
}
\end{remark}

\begin{lemma}
We have the following transformation properties
\begin{eqnarray}
\mathcal{PC}e_{r}(\mathcal{U})=h_{r}(\mathcal{U})\mathcal{PC},&&\mathcal{PC}s_{\la}(\mathcal{U})=s_{\la^t}(\mathcal{U})\mathcal{PC}
\label{PCduality}
\end{eqnarray}
\end{lemma}

\begin{proof}
Using \eqref{PTL} and $\mathcal{PC}u_{N}\mathcal{PC}=u_{N}$
one verifies that each clockwise cyclically ordered monomial is transformed
into an anti-clockwise cyclically ordered one and vice versa,
$\mathcal{PC}\tprod_{i\in I}^{\circlearrowright }u_{i}\mathcal{PC}%
=\tprod_{i\in I}^{\circlearrowleft }u_{i}.$
The first assertion follows then from \eqref{ncehdefTL}, the second from the first with \eqref{defsNil}
\end{proof}

\section{Construction of a common eigenbasis}

We now construct a common eigenbasis for the polynomials $e_{r}(\mathcal{U}%
),h_{r}(\mathcal{U})$, then consider the algebra of noncommutative symmetric functions and establish in this way Postnikov's result that the noncommutative Schur polynomial $s_\la(\cU)$ are well-defined. Our construction and discussion of the eigenbasis will naturally lead to the Bertram-Vafa-Intrilligator formula for Gromov-Witten invariants allowing us
to make contact with the results in \cite{Rietsch} and also set it in parallel to the Verlinde formula \eqref{Verlindeexplicit}.

Analogous to the definition of the Bethe vectors in the context of the fusion ring we now introduce (free) $k$-fermion states $\tilde{b}(x)$ through the following identity which resembles equation \eqref{bethevec}.

\begin{lemma}\label{kfermions}
Let $y=(y_{1},\ldots,y_{k})\in\bigl(\mC[q^{\pm\frac{1}{N}}]\bigr)^k$ with pairwise distinct $y_i$'s. Set $\hat{\psi}^{\ast}(y_{i})=\sum_{p=1}^{N}y_{i}^{-p}\psi _{p}^{\ast }$. Then we have the identity
\begin{equation}
\tilde{b}(y):=\eta(y)
\;\hat{\psi}^{\ast }(y_{1})\cdots \hat{\psi}^{\ast }(y_{k})\emptyset_0 
=\sum_{\lambda \in \mathfrak{P}_{\leq k,n}}\! s_{\lambda
}(y_{1}^{-1},\ldots,y_{k}^{-1})~\lambda,  \label{kfermions2}
\end{equation}
where the normalisation factor $\eta(y)$ is given by
\[
\eta(y)=\frac{(-1)^{\frac{k(k-1)}{2}}e_{k}(y)}{\tprod_{i<j}(y_{i}^{-1}-y_{j}^{-1})}\,.
\]
\end{lemma}

\begin{proof}
Inserting the definition of the creation operators $\hat\psi^\ast(y_i)$ one derives the expansion%
\begin{eqnarray}
\tilde{b}(y)&=&\eta(y)\,\hat{\psi}^{\ast }(y_{1})\cdots \hat{\psi}^{\ast }(y_{k})\emptyset_0\notag\\
&=&\eta(y)\sum_{\pi \in S_{k}}(-1)^{\ell(\pi) }\sum_{\ell _{1}<\ldots<\ell _{k}}y_{\pi
(1)}^{-\ell _{1}}\cdots y_{\pi(k)}^{-\ell _{k}}\psi _{\ell _{1}}^{\ast }\cdots
\psi _{\ell _{k}}^{\ast }\emptyset_0 \label{kfermions3}\\
&=& \sum_{\pi \in S_{k}}(-1)^{\ell(\pi)+\frac{k(k-1)}{2}}\sum_{\ell _{1}<\ldots<\ell _{k}}\frac{y_{\pi(1)}^{-\ell _{1}+1}\cdots y_{\pi(k)}^{-\ell _{k}+1}}{\tprod_{i<j}(y_{i}^{-1}-y_{j}^{-1})}\;
\psi _{\ell _{1}}^{\ast }\cdots
\psi _{\ell _{k}}^{\ast }\emptyset_0\;.\notag
\end{eqnarray}%
Recalling the relationship (\ref{ppositions}) between the particle positions
and the associated partition, $\ell _{k+1-j}-1=k+\lambda _{j}-j$
the statement follows after applying the permutation $w _{0}=(k\ldots21)$ from the
following determinant formula for Schur functions (see e.g. \cite{MacDonald}),%
\begin{equation*}
s_{\lambda }(y^{-1})=\frac{\det (y_{i}^{-k-\lambda _{j}+j})_{1\leq i,j\leq k}%
}{\tprod_{i<j}(y_{i}^{-1}-y_{j}^{-1})}\;.
\end{equation*}
\end{proof}

In the context of the fusion ring we used the relations of the quantum
Yang-Baxter algebra to derive the Bethe Ansatz equations; see Proposition \ref{eigenvalues} and the sketched proof. For the quantum
cohomology ring we will employ the Clifford algebra relations from Proposition \ref{Clifford} (1) and \eqref{kfermions3} to derive a set of
conditions on $y=(y_{1},\ldots ,y_{k})$ such that the vector $\tilde{b}(y)$ is
an eigenvector of the noncommutative symmetric functions.
\begin{lemma}
\label{BetheAnsatzcondition}
Provided $y=(y_{1},\ldots,y_{k})$, with pairwise distinct $y_i$'s, is a solution to the equations%
\begin{equation}
y_{1}^{N}=\cdots =y_{k}^{N}=(-1)^{k-1}q  \label{freebae}
\end{equation}%
the associated vector $\tilde{b}(y)$ from \eqref{kfermions} is a common eigenvector of the
noncommutative elementary and complete symmetric polynomials \eqref{ncehdefTL}. In particular, the eigenvalues are of the following form,
\begin{equation}\label{ncehTLspec}
e_r(\mathcal{U})\tilde{b}(y)=e_r(y)\tilde{b}(y)\qquad\text{and}\qquad
h_r(\mathcal{U})\tilde{b}(y)=h_r(y)\tilde{b}(y)\,.
\end{equation}
\end{lemma}

\begin{proof}
To prove our claim, we first compute the action of $e_r(\mathcal{U})$ on $\tilde{b}(y)$ and then show that the result matches with multiplying $\tilde{b}(y)$ by $e_{r}(y)$ if we impose the conditions \eqref{freebae}. We employ once more the expansion \eqref{kfermions3} of the vector $\tilde b(y)$. Define $v_{\ell _{1},\ldots ,\ell _{k}}=\psi _{\ell _{1}}^{\ast
}\cdots \psi _{\ell _{k}}^{\ast }\emptyset _{0}$ with $v_{\ell _{1},\ldots
,\ell _{k}}=0$ if and only if $\ell _{i}=\ell _{j}$ for $i\neq j$ because of the Clifford algebra relations. From equation \eqref{nTLrep} it follows that%
\begin{multline*}
u_{N}v_{\ell _{1},\ldots ,\ell _{k-1},N}=(-1)^{k+1}\psi _{1}^{\ast }\psi
_{N}~v_{\ell _{1},\ldots ,\ell _{k-1},N}= \\
-q\psi _{1}^{\ast }\psi _{\ell _{1}}^{\ast }\cdots \psi _{\ell _{k-1}}^{\ast
}\psi _{N}\psi _{N}^{\ast }\emptyset _{0}=-qv_{1,\ell _{1},\ldots ,\ell
_{k-1}},
\end{multline*}%
whence we have the quasi-periodic boundary conditions
$$v_{\ell _{1},\ldots,\ell _{k-1},N+1}=-q v_{1,\ell _{1},\ldots ,\ell _{k-1}}.$$ %
For convenience, we define in addition $v_{0,\ell _{2},\ldots ,\ell _{k}}:=-q^{-1}v_{\ell_{2},\ldots ,\ell _{k},N}$. The computation is further simplified by introducing the idempotent $\op{A}_{k}=\sum_{\pi \in S_{k}}(-1)^{\ell(\pi)
}~\pi $ which projects onto the totally antisymmetric subspace; it allows us to rewrite \eqref{kfermions3} as follows,
\begin{equation*}
\tilde{b}(y)=\eta(y)\op{A}_{k}\sum_{0\leq \ell _{1}<\cdots <\ell
_{k}\leq N}y_{1}^{-\ell _{1}}\cdots y_{k}^{-\ell _{k}}v_{\ell _{1},\ldots
,\ell _{k}}\,.
\end{equation*}%
Without loss of generality we may restrict to $r\leq k$, since $e_r(\mathcal{U})_k=0$ for $r>k$. Via induction one then shows that for $i_{1}<\cdots <i_{r}$ (in the cyclic sense) one has%
\begin{equation*}
u_{i_{1}}\cdots u_{i_{r}}v_{\ell _{1},\ldots ,\ell _{k}}=\sum_{1\leq
j_{1}<\cdots <j_{r}\leq k}\delta _{i_{1}\ell _{j_{1}}}\cdots \delta
_{i_{r}\ell _{j_{r}}}v_{\ell _{1},\ldots ,\ell _{j_{1}}+1,\ldots ,\ell
_{j_{r}}+1,\ldots ,\ell _{k}},
\end{equation*}%
Hence, we arrive at the following expression%
\begin{multline}
e_{r}(\mathcal{U})\hat{\psi}^{\ast }(y_{1})\cdots \hat{\psi}^{\ast
}(y_{k})\emptyset _{0} =  \\
\op{A}_{k}\sum_{0\leq \ell _{1}<\cdots <\ell _{k}\leq N}y_{1}^{-\ell
_{1}}\cdots y_{k}^{-\ell _{k}}\sum_{i_{1}<\cdots <i_{r}}^{\circlearrowleft
}u_{i_{1}}\cdots u_{i_{r}}v_{\ell _{1},\ldots ,\ell _{k}}= \\
\op{A}_{k}\sum_{0\leq \ell _{1}<\cdots <\ell _{k}\leq N}y_{1}^{-\ell
_{1}}\cdots y_{k}^{-\ell _{k}}\sum_{1\leq j_{1}<\cdots <j_{r}\leq k}v_{\ell
_{1},\ldots ,\ell _{j_{1}}+1,\ldots ,\ell _{j_{r}}+1,\ldots ,\ell _{k}}\;.\label{fermispec1}
\end{multline}%
In comparison, multiplication of $\tilde{b}(y)$ by $e_r(y)$ yields the formula
\begin{multline}\label{fermispec2}
e_{r}(y)\hat{\psi}^{\ast }(y_{1})\cdots \hat{\psi}^{\ast }(y_{k})\emptyset
_{0}=  \\
\left( \sum_{1\leq j_{1}<\cdots <j_{r}\leq k}y_{j_{1}}\cdots
y_{j_{r}}\right) ~\left( \op{A}_{k}\sum_{0\leq \ell _{1}<\cdots <\ell _{k}\leq N}%
y_{1}^{-\ell _{1}}\cdots y_{k}^{-\ell _{k}}v_{\ell
_{1},\ldots ,\ell _{k}}\right) = \\
\op{A}_{k}\sum_{\substack{ 0\leq \ell _{1}<\cdots <\ell _{k}\leq N \\ 1\leq
j_{1}<\cdots <j_{r}\leq k}}y_{1}^{-\ell _{1}}\cdots y^{-\ell
_{j_{1}}+1}\cdots y^{-\ell _{j_{r}}+1}\cdots y_{k}^{-\ell _{k}}v_{\ell
_{1},\ldots ,\ell _{k}}\;.
\end{multline}%
Notice that the summands for which $\ell _{j_{a}}+1=\ell _{j_{a}+1}$ vanish
under multiplication with $y_{j_{a}}$ due to the antisymmetrizer. This ensures compatibility with the Clifford algebra relations which imply $%
v_{\ell _{1},\ldots ,\ell _{k}}=0$ if two indices are equal. Comparing coefficients in both expressions, \eqref{fermispec1} and \eqref{fermispec2}, we are led to the equations, $%
y_{i}^{N}=(-1)^{k-1}q$.


The proof for $h_{r}(\mathcal{U})$ is slightly more complicated, since now shifts in the particle positions larger than one can occur, but follows along similar lines. We therefore omit it.
\end{proof}

\begin{remark}{\rm
The vectors \eqref{kfermions2} are well-known in the physics literature: via a  {\it Jordan-Wigner transformation} they diagonalize the Hamiltonian $H_{XX}=e_{1}(\mathcal{U})+e_{1}(\mathcal{U})^\ast$ of the so-called {\em XX spin-chain}, mapping it to a free fermion system.
Our computation shows that the eigenvectors of $H_{XX}$ are also eigenvectors of the algebra of noncommutative functions.
}
\end{remark}

We use analogous conventions and notation as in Theorem \ref{Solutions} and denote by
$\widehat{\func{Sol}(n,k)}_{\op{fermi}}$ the solutions $y\in \bigl(\mC[q^{\pm\frac{1}{N}}]\bigr)^k$ of \eqref{freebae}
up to permutation and giving rise to non-trivial $\tilde b(y)$.
\begin{proposition}
\label{freeBetherootsresult}
Let $\zeta =\exp \frac{2\pi\iota}{k+n}$. Then there is a bijection
\begin{eqnarray}
\mathfrak{P}_{\leq n,k}&\cong&\widehat{\func{Sol}(n,k)}_{\op{fermi}}\nonumber\\
\sigma&\mapsto&y_\sigma:=q^{\frac{1}{N}}\left(\zeta^{I_{1}},\ldots,\zeta^{I_{k}}\right),
\label{freeBetheroots}\\
\text{where }I=I(\sigma^t)&:=&\left( \tfrac{k+1}{2}+\sigma^t_{k}-k,\ldots ,\tfrac{k+1}{2}+\sigma^t _{1}-1\right).  \label{freeImap0}
\end{eqnarray}
\end{proposition}
\begin{proof}
In contrast to the Bethe Ansatz equations \eqref{bae} we discussed in the context of the fusion ring, the equations \eqref{freebae} are decoupled and we simply have to take the $N^{\text{th}}$ root in each of the $k$ equations. The $q$ dependence is obvious. Thus, for each single equation there are $N$ solutions.
Given $\sigma$ as in \eqref{freeBetheroots} the corresponding $y_\sigma$ has pairwise distinct components, because for $i>j$ the equality $(y_\sigma)_i=(y_\sigma)_j$ would imply $\sigma_i^t>\sigma^t_j$ which contradicts $\sigma$ being a partition. So $\tilde b(y_\sigma)$  is well-defined via \eqref{kfermions2}. Since we are only interested in solutions $y=(y_1,\ldots,y_k)$ up to permutation of the components, we are left with $\binom{N}{k}=|\mathfrak{P}_{\leq n,k}|=|\mathfrak{P}_{\leq k,n}|=\op{dim}\mathcal{F}_k$ possible solutions to \eqref{freebae}. Proceeding similar to the case of the fusion ring we define a bijection, \eqref{freeImap0}, but this time from $\mathfrak{P}_{\leq n,k}$ into the set of integers
\begin{equation*}
\mathcal{I}_{n,k}:=\left\{(I_{1},...,I_{k})\mid -\tfrac{k-1}{2}\leq
I_{1}<\cdots <I_{k}\leq n+\tfrac{k-1}{2},I_{j}\in \tfrac{k+1}{2}+\mathbb{Z}%
,\forall j\right\} .
\end{equation*}
Taking the $N^{\text{th}}$ power of each component in \eqref{freeBetheroots} one verifies that these indeed solve \eqref{freebae} giving rise to $\binom{N}{k}$ distinct solutions.
\end{proof}

\begin{remark}
\label{baeequiv2}
{\rm
The solutions to (\ref{freebae}) are discussed in \cite{Rietsch} in the
context of the Bertram-Vafa-Intrilligator formula. Employing the well-known recursion relation $h_{r}(y_{1},\ldots ,y_{k})=h_{r}(y_{1},\ldots
,y_{k-1})+y_{k}h_{r-1}(y_{1},\ldots ,y_{k}) $ one can show (compare with Lemma \ref{baeequiv}) that the system of equations \eqref{freebae} is equivalent to
\begin{equation}
h_{N-k+1}(y)=\cdots =h_{N-1}(y)=h_{N}(y)+(-1)^{k}q=0.  \label{ideal}
\end{equation}%
This provides a new perspective: the set of equations \eqref{freebae} coincide with the Bethe Ansatz equations of free fermions on a circular
lattice with quasi-periodic boundary conditions. (Note that these relations are precisely the condition $dW=0$ in the Landau-Ginzburg formulation of the quantum cohomology, see e.g. \cite[p.125]{McDS}.)
 The integers $I$ correspond
to the momenta $p_{r}$ of the individual particles via $p_{r}=2\pi I_{r}/N$
and the total momentum $P=\sum_{r=1}^{k}p_{r}$ is simply given by
\begin{equation}
P=\frac{2\pi }{N}\sum_{r=1}^{k}I_{r}(\sigma )=\frac{2\pi |\sigma |}{N}.
\label{momentum}
\end{equation}
}
\end{remark}

\begin{theorem}
\label{freeeigenbasis}
\begin{enumerate}
\item
The Bethe vectors $\tilde b(y_\sigma)$, $\sigma \in \mathfrak{P}_{\leq k,n}$, form a complete set of pairwise orthogonal eigenvectors for the action of the cyclic noncommutative symmetric functions on $\mathcal{F}_k[q]$.
The eigenvalues are given by the following formulae
\begin{eqnarray}
e_{r}(\mathcal{U})\tilde b(y_\sigma)&=&
\begin{cases}
e_{r}(y_\sigma)\tilde b(y_\sigma)&\text{ if $r\not=N$,}\\
0&\text{ if $r>k$.}
\end{cases}\\
h_{r}(\mathcal{U})\tilde b(y_\sigma)&=&
\begin{cases}
h_{r}(y_\sigma)\tilde b(y_\sigma)&\text{ if $r\not=N$,}\\
-qe^{\iota\pi k}\tilde b(y_\sigma)&\text{ if $r=N$.}\\
\end{cases}
\label{freeSchur spec}
\end{eqnarray}
\item The norm of the Bethe vectors is given by the following formula
\begin{eqnarray}
\langle \tilde b_\sigma,\tilde b_\sigma\rangle&=&\frac{n(n+k)^{n}}{|\op{Van}_\sigma|^2},
\end{eqnarray}
where $\op{Van}_\sigma$ denotes the Vandermonde determinant $\displaystyle{\prod_{1\leq i<j\leq n}({\zeta ^{I_{i}(\sigma)}-\zeta ^{I_{j}(\sigma)}}})$.
\end{enumerate}
\end{theorem}

\begin{proof}
The eigenvalue expressions for $e_r(\mathcal{U}),\;h_r(\mathcal{U})$ follow from  Lemma \ref{BetheAnsatzcondition}.
The special cases  $e_r(\mathcal{U}),\;r>k$ and $h_N(\mathcal{U})$ are immediate from the fact that there are only $k$-particles in the system and \eqref{ideal}, respectively.

Orthogonality and the norm of the Bethe vectors follows from \eqref{kfermions2} and the summation identity for Schur functions in \cite{Rietsch}; see Proposition 4.3 (3) therein.
\end{proof}

We arrive at the following statement, compare with \cite{Postnikov} where the commutativity of $\{e_{r}(\mathcal{U}),h_{s}(\mathcal{U})\}$ was shown by different methods.
\begin{corollary}
\label{Nilecomm}
The elements in the set $\{e_{r}(\mathcal{U}),h_{s}(\mathcal{U})\}$ pairwise
commute. Thus, the noncommutative Schur polynomials $s_\la(\cU)$ are well defined and one has the eigenvalue equation
\begin{equation}
s_{\lambda }(\mathcal{U})\tilde b(y_{\sigma })=s_{\lambda}(y_{\sigma
})\tilde b(y_{\sigma })=q^{\frac{|\lambda |}{N}}s _{\lambda}(\zeta^{I(\sigma)
})\tilde b(y_{\sigma }).
\end{equation}
Hence, both determinant formulae in \eqref{defsNil} for $s_\la(\cU)$ are valid.
\end{corollary}

\begin{proof}
Theorem \ref{freeeigenbasis} and Proposition \ref{freeBetherootsresult} imply that if we formally set $q\cdot \bar{q}=q\cdot
q^{-1}=1$ then the non-zero vectors (\ref{kfermions}) form a complete set of pairwise orthogonal eigenvectors, hence an eigenbasis of
the noncommutative polynomials $e_{r}(\mathcal{U}),h_{s}(\mathcal{U})$. The eigenvalue equation for the noncommutative Schur polynomial is then immediate from the definition (1) in \eqref{defsNil}. Moreover, the equality (2) in \eqref{defsNil} is then deduced from the corresponding equality in the commutative case.
\end{proof}

\begin{theorem}[Combinatorial quantum cohomology ring]
\label{freecombinatorial}
Fix $k\in \mathbb{Z}_{\geq 0}$ and consider the $k$-particle subspace $
\mathcal{F}_{k}[q]\subset\mathcal{F}[q]$. The assignment
\begin{equation}
({\lambda},{\mu})\mapsto {\lambda}\star {\mu}%
:=s_{\lambda }(\mathcal{U}){\mu}  \label{freeSchurproduct}
\end{equation}%
for basis elements ${\lambda},{\mu}\in \mathfrak{P}_{\leq k,n}$ turns $\mathcal{F}_{k}[q]$ into a commutative, associative and unital  $\mathbb{C}[q]$-algebra $H^\bullet_{\op{comb}}$.
\end{theorem}

\begin{proof}
The proof is completely analogous, with the obvious replacements, to the one of Theorem \ref{combinatorial}.
\end{proof}

Define the structure constants  $C_{\la,\mu}^{\nu,d}\in\mZ$ of $H^\bullet_{\op{comb}}$ through $$\la\star\mu=\sum_{\nu\in\mathfrak{P}_{\leq k,n},d\in\mZ} C_{\la,\mu}^{\nu,d}q^d\nu.$$

\begin{corollary}
\label{BVIcomb}
The following formula describes the structure constants $C_{\la,\mu}^{\nu,d}$ for $\la,\mu,\nu\in \mathfrak{P}_{\leq k,n}$:
\begin{equation}
C_{\lambda ,\mu }^{\nu ,d}q^d=\frac{q^{\frac{|\lambda |+|\mu |-|\nu |}{N}}%
}{N^{k}}\sum_{\sigma\in \mathfrak{P}_{\leq k,n}}\frac{s_{\lambda }(\zeta ^{I(\sigma)})s_{\mu }(\zeta ^{I(\sigma)})s_{\nu^\vee%
}(\zeta ^{-I(\sigma)})}{\zeta^{k|\sigma|}\prod_{i<j}^{k}|\zeta ^{I_{i}(\sigma)}-\zeta
^{I_{j}(\sigma)}|^{-2}}  \label{BVI}
\end{equation}%
In particular, $C_{\lambda ,\mu }^{\nu ,d}=0$ unless $|\lambda |+|\mu |-|\nu |=dN$ and $d\geq 0$.
\end{corollary}

\begin{proof}
From \eqref{nTLrep} and \eqref{freeSchurproduct} we can conclude that $d$ is non-negative as the noncommutative Schur polynomial $s_\la(\mathcal{U})$ is polynomial in the generators of the affine nil-Temperley-Lieb algebra; see \eqref{defsNil}. According to \eqref{freeSchurproduct} the coefficients in the product expansion $\la\star\mu$ are given by the matrix elements of $s_\la(\mathcal{U})$ which we compute in the eigenbasis \eqref{kfermions}, completely analogous to the derivation of \eqref{permanentformula}.
\end{proof}

We illustrate the last formula by an example:

\begin{example}
\label{ExGW}
{\rm
Consider the case $k=4$ and $%
n=3$, hence $N=7$. Set $\lambda =(3,3,2,1)$ and $\mu =(2,2,1,0)$. Then we have $|\la|+|\mu|=9+5=14$ and therefore $C_{\lambda ,\mu }^{\nu ,d}=0$ unless $7d+|\nu|=14$. Therefore $d\in\{0,1,2\}$. The case $d=0$ gives a partition $\nu$ with $14$ boxes which is impossible, $d=2$ implies $\nu=\emptyset$ and $d=1$ implies that $\nu$ has $7$ boxes, hence we get a formula of the form
\begin{equation*}
\Yvcentermath1\yng(3,3,2,1)\star\yng(2,2,1)= \\
\Yvcentermath1 a_1q~\yng(2,2,2,1)+a_2q~\yng%
(3,2,1,1)+a_3q~\yng(3,2,2)+a_4q~\yng(3,3,1)%
+a_5q^{2}\emptyset\;.
\end{equation*}
for some $a_i\in\mathbb{Z}$. We used a computer to evaluate the formula \eqref{BVI} and obtained $a_i=1$ except for $a_2=2$.
}
\end{example}

The following reproves a result of Postnikov \cite{Postnikov}, and at the same time connects his approach with the results of Rietsch \cite{Rietsch}:

\begin{theorem}
The map $s_\la(\cU)\mapsto s_\la$ defines an isomorphism of rings   $H^\bullet_{\op{comb}}\cong qH^{\bullet }(\limfunc{Gr}\nolimits_{k,N})$.
In particular,
\begin{equation}\label{isoH}
H^\bullet_{\op{comb}}\cong \mathbb{Z}[q]\otimes_\mZ\mZ[e_1,e_2,\ldots e_k]/
\left\langle h_{n+1},\ldots,h_{n+k-1},h_{n+k}+(-1)^{k}q\right\rangle
\end{equation}
and the Gromov-Witten invariants are given by
the matrix elements of the noncommutative Schur function,
\begin{equation}
q^dC_{\lambda ,\mu }^{\nu ,d}=\langle \nu ,s_{\lambda}(\mathcal{U}%
)\mu \rangle \;.  \label{Presult}
\end{equation}
\end{theorem}

\begin{proof}
The expansion \eqref{BVI} is the celebrated Bertam-Vafa-Intrilligator formula for Gromow-Witten invariants. In the stated form it has appeared in \cite{Rietsch}. The isomorphism \eqref{isoH} follows then from Remark \ref{baeequiv2} by the same arguments as explained in detail in the proof of Theorem \ref{mainresult}. The remaining statements are then direct consequences from the presentation \eqref{SiebertTian} of the quantum cohomology ring.
\end{proof}

\begin{remark}{\rm
We want to point out here that the definitions and basic properties of the combinatorial ring  $H^\bullet_{\op{comb}}$ only rely on \eqref{freeSchurproduct}. From there one derives in an elementary way the associativity, the Bertram-Vafa-Intrilligator formula \eqref{BVI} and the presentation \eqref{isoH}.
All of these are non-trivial results in the quantum cohomology of the Grassmannian.}
\end{remark}

\section{Symmetries, Recursion Formulae and an Algorithm for Gromov-Witten Invariants}

The combinatorial description of the quantum cohomology ring can be exploited to derive in a very simple manner non-trivial identities, referred
to as  `symmetries', between Gromov-Witten
invariants. The main idea is to employ the combinatorial description \eqref{freeSchurproduct} of the product and then exploit certain invariance or transformation properties of the noncommutative Schur polynomials. %

To explain our approach in more detail, assume we have a linear
endomorphism $F$ of $\mathcal{F}[q]$ which sends a basis vector to a basis vector. For notational convenience we shall abbreviate throughout this section $S_\la=s_\la(\mathcal{U}),\, E_r =e_r(\mathcal{U})$ etc. In general we cannot expect the endomorphism $F$ to be compatible with the product structure in the sense that $S_\la\circ F=F\circ S_\la$. However, we might be able to choose $F$ such that
$S_\la\circ F=\sum_\mu F_\mu\circ S_\mu^\prime$   for  certain (well-understood) endomorphisms $F_\mu$ and
$S_\mu^\prime=S_\mu$ at least after
replacing $q$ by $q^\prime=-q$ or $q^{-1}$. Applying this equality to a basis vector in $\cF[q]$ we obtain identities of Gromov-Witten invariants by comparing coefficients in the respective product expansions,
$$
\sum_\nu q^d C_{{\la},F(\sigma)}^{\nu,d}\nu=S_\la (F(\sigma))=\sum_\mu F_\mu (S_\mu^\prime(\sigma))=\sum_\nu\sum_\mu
(q^\prime)^d C_{\mu,\sigma}^{\nu,d}F_\mu(\nu)\,.
$$

In case that the endomorphisms $F,F_\mu$ preserve the dimension $k$
and $q'=q$, we obtain relations between structure constants of the same quantum cohomology ring. Our standard examples are the discrete symmetry operators from Section \ref{trans}. In this way we deduce in particular (mostly well-known)
symmetry properties of Gromov-Witten invariants, albeit providing simpler and shorter proofs (Proposition \ref{GWsym} below).  By choosing for $F$ particle generation and
annihilation operators $\psi_i$ and $\psi_i^\ast$ we obtain, however, more interesting (new) results, namely formulae relating Gromov-Witten invariants from {\em different} quantum cohomology rings. In particular, one obtains an explicit recursion formula and an inductive algorithm to compute all the structure constants.

We start with the following (apart from (\ref{GWsym3}) probably) well-known symmetry formulae:

\begin{proposition}\label{GWsym}
Let $C_{\lambda ,\mu ,\nu }(q)=q^{d}C_{\lambda ,\mu }^{\nu ^{\vee },d}$ with $d=%
\frac{|\lambda |+|\mu |+|\nu |-kn}{N}~,
$ %
then we have the following identities:

\begin{enumerate}
\item\label{GWsym1} {\em $S_3$ symmetry:} for any $p\in S_3$ we have
$C_{\lambda,\mu \,\nu }(q)=C_{p(\lambda ),p(\mu ),p(\nu )}(q)\,.$

\item\label{GWsym2} {\em Level-rank duality:} $C_{\lambda ,\mu ,\nu }(q)=C_{\lambda
^{t},\mu ^{t},\nu ^{t}}(q)$.

\item\label{GWsym3} {\em Rotation invariance:} let $\op{R}=\limfunc{Rot}$ be as in \eqref{RotN} and extend the definition of $\op{n}_i$ from Section \ref{prelim} to any $i\in\mathbb{Z}$ by setting $\op{n} _{i+N}=\op{n} _{i}+k$. Then
\begin{equation*}
C_{\op{R}^{a}(\lambda ),\mu ,\nu }(q)=q^{\op{n} _{a}(\mu
)-\op{n} _{a}(\lambda )}C_{\lambda ,\op{R}^{a}(\mu ),\nu
}(q)=q^{\op{n} _{a}(\nu )-\op{n} _{a}(\lambda )}C_{\lambda ,\mu ,\op{R}^{a}(\nu )}(q),  
\end{equation*}%
and in particular for $a+b+c=0$ we have%
\begin{equation*}
C_{\op{R}^{a}(\lambda ),\op{R}^{b}(%
\mu ),\op{R}^{c}(\nu )}(q)=q^{\op{n} _{a}(\lambda )+\op{n}
_{b}(\mu )+\op{n} _{c}(\nu )}C_{\lambda ,\mu ,\nu }(1)\;.
\end{equation*}

\item\label{GWsym4} {\em Postnikov's `curious duality':}%
\begin{equation*}
C_{\op{R}^{a}(\lambda ),\op{R}^{b}(%
\mu ),\op{R}^{c}(\nu )}(q)=q^{\op{n} _{a}(\lambda )+\op{n}
_{b}(\mu )+\op{n} _{c}(\nu )}C_{\lambda^\vee,\mu^\vee,\nu^\vee}(q^{-1})
\end{equation*}%
with $a+b+c=N-k$.
\end{enumerate}
\end{proposition}

\begin{proof}
The first identity is obvious
from the Bertram-Vafa-Intrilligator formula \eqref{BVI}. The
second equality is easily obtained by employing parity-reversal and
particle-hole duality,
\begin{equation*}
C_{\lambda ,\mu ,\nu }(q)=\langle \nu ^{\vee },S_\lambda\mu \rangle =\langle \left( \nu ^{t}\right) ^{\vee },\mathcal{PC}%
S_\lambda\mathcal{PC}\mu \rangle =\langle \left( \nu
^{t}\right) ^{\vee },S_{\lambda^t }\mu ^{t}\rangle =C_{\lambda
^{t},\mu ^{t},\nu ^{t}}(q)\;.
\end{equation*}%
To prove statement \eqref{GWsym3} we first specialise to $q=1$, denoting the corresponding noncommutative Schur polynomial by $\tilde S_\la$. Under this condition the noncommutative elementary and symmetric polynomials in \eqref{ncehdefTL} become invariant under a cyclic permutation of the letters $u_i$. Together with \eqref{defsNil} this implies %
\begin{equation}\label{RotS}
\op{R}\circ \tilde{S}_{\lambda}=\tilde{S}_{\lambda}\circ \op{R}\,.
\end{equation}
Hence, the following equalities hold true%
\begin{multline*}
\langle \nu ^{\vee },\op{R}\circ\tilde S_\la(\mu)\rangle
=\langle\nu ^{\vee },\tilde S_\la(\op{R}\mu)\rangle
=\langle \nu ^{\vee },\op{R}\circ\tilde S_\la(\mu)\rangle=\\
\langle \nu ^{\vee },\op{R}\circ\tilde S_\mu(\la)\rangle
=\langle \nu ^{\vee },\tilde S_\mu(\op{R}\la)\rangle
=\langle \nu ^{\vee },\op{R}\circ\tilde S_\mu(\la)\rangle
\end{multline*}%
Using the $S_3$-symmetry from \eqref{GWsym1} statement \eqref{GWsym3} follows for $q=1$. Because of Corollary \ref{BVIcomb} it suffices to compute the difference in the degrees of the
respective Gromov-Witten invariants in order to treat the case of generic $q$. The respective degrees in \eqref{GWsym3} can be deduced with help of the general formula
\begin{equation}\label{Rotweight}
|\op{R}^{a}(\mu )|=|\mu |+N\op{n} _{a}(\mu )-ak\;,
\end{equation}%
which is easily verified using the action of $\op{R}$ on Young diagrams as described in \eqref{RotN}. For instance, the degrees for the last equality in \eqref{GWsym3} follow from the computation %
\begin{multline*}
|\op{R}^{a}(\lambda )|+|\op{R}^{b}(%
\mu )|+|\op{R}^{c}(\nu )|= \\
|\lambda |+|\mu |+|\nu |+N(\op{n} _{a}(\lambda )+\op{n} _{b}(\mu )+\op{n} _{c}(\nu
))+(a+b+c)k
\end{multline*}%
and, hence, we arrive at the asserted identity under the stated assumption
that $a+b+c=0$.

To show the last identity \eqref{GWsym4} we again specialise to $q=1$. Employing the parity and time reversal operator of Section \ref{trans} one shows that the eigenvector (\ref{kfermions2}) satisfies
\[
\mathcal{PT}\tilde{b}(y)=\sum_{\mu\in\mathfrak{P}_{\leq k,n}}s_{\mu}(y)\mu^\vee
=e_k(y)^n\sum_{\mu\in\mathfrak{P}_{\leq k,n}}s_{\mu^\vee}(y^{-1})\mu^\vee
=e_k(y)^n\tilde b(y),
\]
see Lemma \ref{Schurprop} \eqref{easy3}. Applying Theorem \ref{freeeigenbasis}, Corollary \ref{Nilecomm} and once more Lemma \ref{Schurprop} \eqref{easy3} one obtains the identity
$
\mathcal{PT}\tilde{S}_\lambda\mathcal{PT}=\op{R}^{n}\circ \tilde{S}_{\lambda ^{\vee }}\,.
$  Thus, we calculate
\begin{multline*}
C_{\lambda ,\mu ,\nu } =\langle \nu ^{\vee },\tilde{S}_\lambda\mu \rangle =\langle \nu ,\mathcal{PT}\tilde{S}_\lambda\mathcal{PT}\mu ^{\vee}\rangle=
\langle \nu ,\op{R}^{n}\tilde{S}_{\lambda ^{\vee }}\mu ^{\vee }\rangle =\\
\langle
\op{R}^{-c}(\nu ),\tilde{S}_{\op{R}^{a}(\lambda ^{\vee })}\op{R}^{-b}(\mu
)^{\vee }\rangle
=C_{\op{R}^{-a}(\lambda )^{\vee },\op{R}^{-b}(\mu )^{\vee },\op{R}^{-c}(\nu )^{\vee
}}
\end{multline*}%
for any $a,b,c$ with $a+b+c=n=N-k$. Here we have used again \eqref{RotS} in the penultimate equality. This proves \eqref{GWsym4} when $q=1$. The general
case now follows by comparing degrees on both sides of \eqref{GWsym4}. We find%
\begin{equation*}
C_{\lambda ^{\vee },\mu ^{\vee },\nu ^{\vee }}(q)=q^{d}C_{\lambda ^{\vee
},\mu ^{\vee },\nu ^{\vee }}(1),\qquad d=-\frac{|\lambda |+|\mu |+|\nu |-2k(N-k)%
}{N}
\end{equation*}%
and
\begin{equation*}
C_{\op{R}^{a}(\lambda ),\op{R}^{b}(%
\mu ),\op{R}^{c}(\nu )}(q)=q^{d^{\prime }}C_{\op{R}^{a}(\lambda ),\op{R}^{b}(\mu ),\op{R}^{c}(\nu )}(1)
\end{equation*}%
with%
\begin{equation*}
d^{\prime }=\frac{|\lambda |+|\mu |+|\nu |-k(N-k)}{N}+\op{n} _{a}(\lambda
)+\op{n} _{b}(\mu )+\op{n} _{c}(\nu )-\frac{a+b+c}{N}k\;
\end{equation*}
and the proposition follows.
\end{proof}

\begin{example}{\rm
Let us illustrate the second equality in the first formula of Proposition \ref{GWsym} \eqref{GWsym3}. Recall from Example \ref{ExGW} the product expansion for $\la\star\mu$,
\begin{equation*}
\Yvcentermath1\underset{0101011}{\yng(3,3,2,1)}\star\underset{1010110}{\yng(2,2,1)}= \\
\Yvcentermath1q~\underset{0101110}{\yng(2,2,2,1)}+2q~\underset{0110101}{\yng%
(3,2,1,1)}+q~\underset{1001101}{\yng(3,2,2)}+q~\underset{1010011}{\yng(3,3,1)%
}+q^{2}~\underset{1111000}{\emptyset }\;.
\end{equation*}%
Exploiting that $\op{Rot}(\nu^\vee)=(\op{Rot}^{-1}(\nu))^\vee$, Proposition \ref{GWsym} \eqref{GWsym3} can be rewritten as %
\begin{eqnarray*}
\limfunc{Rot}(\lambda )\star \mu  &=&\sum_{\nu }C_{\lambda,
\mu,\op{Rot(\nu^\vee)} }(q)\nu
=\sum_{\nu }q^{d^{\prime }}C_{\lambda
\mu }^{\limfunc{Rot}\nolimits^{-1}(\nu ),d^{\prime }}\nu  \\
&=&\sum_{\nu }q^{d-k+\op{n} _{1}(\lambda )+\op{n} _{6}(\limfunc{Rot}(\nu
))}C_{\lambda \mu }^{\nu ,d}\limfunc{Rot}(\nu )\;,
\end{eqnarray*}%
where $d^{\prime
}(\nu)=d+k-\op{n} _{N-1}(\nu )+\op{n} _{1}(\lambda )$ according to the definition of $C_{\la,\mu,\nu}(q)$ in Proposition \ref{GWsym} and \eqref{Rotweight}. Since, $\op{n} _{1}(\lambda )=0$ and $\op{Rot}(\la)=\mu$ we find by applying $\op{Rot}$ to each diagram in the product expansion%
\begin{equation*}
\Yvcentermath1\yng(2,2,1)\star \yng(2,2,1)\stackrel{\eqref{RotN}}{=}
\Yvcentermath1q~\yng(1,1,1)+2q~\yng(2,1)+\yng(3,3,2,2)+\yng(3,3,3,1)+q~\yng%
(3)\;.
\end{equation*}
}
\end{example}
We now derive recursion formulae for Gromov-Witten invariants. As a preliminary step we need to derive the commutation relations between $\psi _{i}^{\ast },\psi _{i}$ and the noncommutative Schur polynomials. We start with the elementary symmetric polynomials and then, once more, make use of the Jacobi-Trudy formula. Henceforth let $S_\la^\prime$ and $E_r^\prime$ denote $S_\la$ and $E_r$ with $q$ replaced by $-q$, respectively.

\begin{lemma}
We have the following formula%
\begin{equation*}
E_{r}\psi _{j}^{\ast }=\psi _{j}^{\ast }E_{r}^\prime+\psi _{j+1}^{\ast }E_{r-1}^\prime,
\end{equation*}%
where for $j=N$ we have $\psi _{N+1}^{\ast }=-qe^{i\pi K}\psi _{1}^{\ast }$.
\end{lemma}

\begin{proof}
We start from the relations%
\begin{eqnarray*}
u_{i}\psi _{l}^{\ast } &=&\psi _{l}^{\ast }u_{i}+\delta _{i,l}\psi
_{i+1}^{\ast },\qquad i=1,2,\ldots,N-1 \\
u_{N}\psi _{l}^{\ast } &=&\psi _{l}^{\ast }u_{N}^\prime+\delta _{N,l}~q\psi
_{1}^{\ast }e^{i\pi K}
\end{eqnarray*}%
which are easily verified. More generally we claim that for any $i_{1}>\cdots >i_{r}$ (with
respect to the cyclic order) the following identity is true,%
\begin{equation*}
u_{i_{r}}\cdots u_{i_{1}}\psi _{j}^{\ast }=\psi _{j}^{\ast
}u_{i_{r}}^\prime\cdots u_{i_{1}}^\prime
+\psi _{j+1}^{\ast }\sum_{s=1}^{r}\delta _{j,i_{s}}u_{i_{r}}^\prime\cdots \NEG%
{\,u_{i_{s}}^\prime}\cdots u_{i_{1}}^\prime\;,
\end{equation*}%
where $\NEG{\,u}$ means we leave out the factor $u$. Assume this is true for $r$, then multiplication on both sides with $u_{i_{r+1}}$ yields%
\begin{multline*}
u_{i_{r+1}}\cdots u_{i_{1}}\psi _{j}^{\ast }=\psi _{j}^{\ast
}u_{i_{r+1}}^\prime\cdots u_{i_{1}}^\prime\\+\delta _{j,i_{r+1}}\psi _{j+1}^{\ast
}u_{i_{r}}^\prime\cdots u_{i_{1}}^\prime
+u_{i_{r+1}}\psi _{j+1}^{\ast }\sum_{s=1}^{r}\delta
_{j,i_{s}}u_{i_{r}}^\prime\cdots \NEG{\,u}_{i_{s}}^\prime\cdots u_{i_{1}}^\prime\;.
\end{multline*}%
Provided that $i_{r+1}<i_{r}$ one has $u_{i_{r+1}}\psi _{j+1}^{\ast
}=\psi _{j+1}^{\ast }u_{i_{r+1}}^\prime$ which concludes the induction step and, thus, proves the claim.
Noting that $\psi _{j+1}^{\ast }u_{j}=0$ the assertion follows from the
definition \eqref{ncehdefTL} of the cyclic noncommutative symmetric functions.
\end{proof}

\begin{proposition}
The following commutation relations are satisfied,
\begin{eqnarray}
S_{\lambda }\psi _{i}^{\ast } &=&\psi _{i}^{\ast }S_{\lambda
}^\prime+\sum_{r=1}^{\lambda _{1}}\psi _{i+r}^{\ast
}\sum_{\lambda /\mu =(r)}S_{\mu }^\prime  \label{schurcom1} \\
S_{\lambda }\psi _{i} &=&\psi _{i}S_{\lambda }^\prime%
+\sum_{r=1}^{\ell (\lambda )}(-1)^{r}\psi _{i-r}\sum_{\lambda /\mu
=(1^{r})}S_{\mu }^\prime  \label{schurcom2}
\end{eqnarray}%
where we again impose the quasi-periodic boundary conditions $\psi _{j-N}=-qe^{i\pi K}\psi _{j}$ and $\psi _{j+N}^{\ast
}=-qe^{i\pi K}\psi _{j}^{\ast }$.
\end{proposition}

\begin{proof}
We start with the first identity and perform an induction on the length of $\la$ with  the previous lemma as starting point. Assume as induction
hypothesis that (\ref{schurcom1}) holds true for $\la$ and let $%
\tilde{\lambda}=(\lambda _{0},\lambda _{1},\lambda
_{2},\ldots)$ be $\la$ with an additional row added to the top. For the induction step we will exploit $S_{\tilde{\lambda}^t}=\det (E_{\tilde{\lambda}_{i}-i+j})$. %
Expanding with respect to the first column one finds %
$
S_{\tilde{\lambda}^t}=\sum_{i=0}^{\ell }(-1)^{i}E%
_{\lambda _{i}-i}S_{(\tilde{\lambda}^{(i)})^t},
$ %
where $\tilde{\lambda}^{(0)}=\lambda $ and for $i>0$ we set %
$
\tilde{\lambda}^{(i)}=(\lambda _{0}+1,\lambda _{1}+1,\ldots,\lambda
_{i-1}+1,\lambda _{i+1},\ldots,\lambda _{\ell })\;.
$ %
The induction hypothesis then yields%
\begin{eqnarray*}
S_{\tilde{\lambda}^t}\psi _{j}^{\ast } &=&\sum_{i=0}^{\ell
}(-1)^{i}E_{\lambda _{i}-i}\sum_{r=0}^{\ell }\psi _{j+r}^{\ast
}\sum_{\tilde{\lambda}^{(i)}/\tilde{\mu}=(1^{r})}S^\prime_{%
\tilde{\mu}^t} \\
&=&\sum_{i,r=0}^{\ell }(-1)^{i}\left\{ \left( \psi _{j+r}^{\ast }%
E^\prime_{\lambda _{i}-i}+\psi _{j+r+1}^{\ast }E^\prime_{\lambda
_{i}-i-1}\right) \tsum_{\tilde{\lambda}^{(i)}/\tilde{\mu}=(1^{r})}%
S^\prime_{\tilde{\mu}^t}\right\}
\end{eqnarray*}%
We now rearrange the sum over $r$ and find%
\begin{multline*}
S_{\tilde{\lambda}^t}\psi _{j}^{\ast } =\psi _{j}^{\ast }%
S^\prime_{\tilde{\lambda}^t}+ \\
+\sum_{r=1}^{\ell }\psi _{j+r}^{\ast }\sum_{i=0}^{\ell }(-1)^{i}\left\{
E^\prime_{\lambda _{i}-i}\tsum_{\tilde{\lambda}^{(i)}/\tilde{\mu}%
=(1^{r})}S^\prime_{\tilde{\mu}}+E^\prime_{\lambda
_{i}-i-1}\tsum_{\tilde{\lambda}^{(i)}/\tilde{\mu}=(1^{r-1})}%
S^\prime_{\tilde{\mu}^t}\right\} \\
+\sum_{i=0}^{\ell }(-1)^{i}\psi _{j+\ell +1}^{\ast }E^\prime%
_{\lambda _{i}-i-1}S^\prime_{(\lambda _{0},\lambda
_{1},\ldots,\lambda _{i-1},\lambda _{i+1}-1,\ldots,\lambda _{\ell }-1)^t}
\end{multline*}

The two sums in the second line give the two possible contributions when one
subtracts a vertical $r$-strip from $\tilde{\lambda}$ with $1\leq r\leq \ell
$: either a single box or no box is deleted in the $i^{\text{th}}$ row. The
last term is the removal of a vertical ($\ell +1$)-strip from $\tilde{\lambda%
}$. Thus, we arrive at
\begin{equation*}
S_{\tilde{\lambda}^t}\psi _{j}^{\ast }=\psi _{j}^{\ast }%
S^\prime_{\tilde{\lambda}^t}+\sum_{r=1}^{\ell }\psi _{j+r}^{\ast
}\sum_{\tilde{\lambda}/\tilde{\mu}=(1^{r})}S^\prime_{\tilde{\mu}^t%
}+\psi _{j+\ell +1}^{\ast }S^\prime_{(\lambda _{0}-1,\lambda
_{1}-1,\ldots,\lambda _{\ell }-1)^t}
\end{equation*}%
which upon taking the transpose in all partitions is the desired result.

The second formula (\ref{schurcom2}) now follows from applying the
transformations
\begin{equation*}
\mathcal{PC}\psi _{j}^{\ast }\mathcal{PC}=(-1)^{j-1}e^{i\pi K}\psi _{N+1-j}\quad
\text{and}\quad\mathcal{PC}S_\lambda\mathcal{PC}=S_{\lambda
^{t}}
\end{equation*}%
which we proved earlier; see Section \ref{trans} and Lemma \ref{PCduality}.
\end{proof}

\begin{remark}[Algorithm]\label{GWalgo}{\rm
The above commutation relations imply the following product formulae%
\begin{eqnarray}
\lambda \star \psi _{i}^{\ast }(\mu ) &=&S_{\lambda }\psi
_{i}^{\ast }(\mu )=\sum_{r=0}^{\lambda _{1}}\sum_{\lambda /\nu =(r)}\psi
_{i+r}^{\ast }(\nu ~\star^\prime~\mu ), \\
\lambda \star \psi _{i}(\mu ) &=&S_{\lambda }\psi
_{i}(\mu )=\sum_{r=0}^{\lambda _{1}}(-1)^{r}\sum_{\lambda /\nu =(1^{r})}\psi
_{i-r}(\nu ~\star^\prime~\mu ),
\end{eqnarray}%
where $\star^\prime$ denotes the product where we replace the deformation
parameter $q$ with $-q$. Thus, one can successively create for $k=0,1,\ldots ,N$ the entire hierarchy $qH^{\bullet }(\func{Gr}_{k,N})$ of quantum cohomology rings
starting with either $k=0$ or $k=N$.
}
\end{remark}

\begin{example}{\rm
We illustrate the algorithm by computing the (quite trivial) example of projective space and some product in $qH^{\bullet }(\func{Gr}_{2,5})$ starting from the quantum cohomology of a point. First of all we create the basis vectors of $qH^{\bullet }(\func{Gr}_{1,N})$ according to the formula%
\begin{equation*}
c_{j}=\underset{10000}{\emptyset }\star \psi _{j}^{\ast }(\underset{00000}{%
\emptyset })=\psi _{j}^{\ast }(\emptyset )=\Yvcentermath1\underset{j-1}{%
\underbrace{\yng(4)}}~,\quad j=1,\ldots ,N\,.
\end{equation*}%
Then the product of two arbitrary basis vectors is given by%
\begin{eqnarray*}
\Yvcentermath1 c_i\star c_j=\underset{i}{\underbrace{\yng(4)}}\star \psi _{j}^{\ast }(%
\underset{00000}{\emptyset }) &=&\sum_{i^{\prime }=0}^{i}\psi _{i^{\prime
}+j}^{\ast }(h_{i-i^{\prime }}(\mathcal{U},-q)\underset{00000}{\emptyset }%
)=\psi _{i+j}^{\ast }(\underset{00000}{\emptyset }) \\
&=&\Yvcentermath1q^{p}~\underset{i+j-pN}{\underbrace{\yng(4)}}\;,
\end{eqnarray*}%
where $p=0$ if $i+j<N=n+1$ and $p=1$ otherwise (by invoking the
quasi-periodic boundary conditions $\psi _{N+r}^{\ast }=(-1)^{K+1}q\psi
_{r}^{\ast }$). Now we can pass to the case $k=2$ and compute in  $qH^{\bullet }(\func{Gr}_{2,5})$ for instance %
\begin{multline*}
\Yvcentermath1\underset{00101}{\yng(3,2)}\star \psi _{2}^{\ast }\left(
\underset{00010}{\yng(3)}\right)  =\Yvcentermath1\underset{00101}{\yng(3,2)%
}\star \underset{01010}{\yng(2,1)} \\
=\Yvcentermath1\psi _{2+2}^{\ast }\left( \yng(3)~\star^\prime~\yng%
(3)~\right) +\psi _{2+3}^{\ast }\left( \yng(2)~\star^\prime~\yng(3)~\right)  \\
=\Yvcentermath1-q\psi _{2+2}^{\ast }\left( \underset{01000}{\yng(1)}%
\right) -q\psi _{2+3}^{\ast }\left( \underset{10000}{\emptyset }\right)
=\Yvcentermath1q~\yng(2,1)+q~\yng(3)~.
\end{multline*}%
Note that the cases $k=3,4,5$ and $N=5$ are obtained via level-rank duality, see Proposition \ref{GWsym} \eqref{GWsym2}.
}
\end{example}

Using the identity $\psi _{i}^{\ast }\psi _{i}+\psi _{i}\psi _{i}^{\ast }=1$
one now derives the corresponding recursion formulae for the Gromov-Witten
invariants.

\begin{corollary}[Recursion formulae]\label{GWrec}
Let $j=1,2,\ldots ,N$. One has the following recursion formulae for
Gromov-Witten invariants: if $\psi _{j}\mu =0$ then%
\begin{equation}\label{GWrec1}
C_{\lambda \mu }^{\nu ,d}(k,N)=\sum_{r=0}^{\ell (\lambda )}(-1)^{d+r+\op{n}
_{j-1}(\mu )+\op{n} _{j-r-1}(\nu )}\sum_{\lambda /\rho =(1^{r})}C_{\rho ,\psi
_{j}^{\ast }\mu }^{\psi _{j-r}^{\ast }\nu ,d_{r}}(k+1,N)
\end{equation}%
else%
\begin{equation}\label{GWrec2}
C_{\lambda \mu }^{\nu ,d}(k,N)=\sum_{r=0}^{\lambda _{1}}(-1)^{d+\op{n}
_{j-1}(\mu )+\op{n} _{j+r-1}(\nu )}\sum_{\lambda /\rho =(r)}C_{\rho ,\psi
_{j}\mu }^{\psi _{j+r}\nu ,d_{r}^{\prime }}(k-1,N),
\end{equation}%
where $\op{n} _{j+N}=\op{n} _{j}+k$ and
\begin{equation}
d_{r}=\left\{
\begin{array}{cc}
d-1, & j<r \\
d, & \text{else}%
\end{array}%
\right. ,\qquad d_{r}^{\prime }=\left\{
\begin{array}{cc}
d-1, & j+r>N \\
d, & \text{else}%
\end{array}%
\right. \;.
\end{equation}
\end{corollary}

\begin{proof}
We start from the identity%
\begin{equation*}
q^{d}C_{\lambda ,\mu }^{\nu ,d}(k,N)=\langle \nu ,S_\lambda\mu \rangle =\langle \nu ,S_\lambda\psi _{j}\psi
_{j}^{\ast }\mu \rangle +\langle \nu ,S_\lambda\psi
_{j}^{\ast }\psi _{j}\mu \rangle \;.
\end{equation*}%
Only one of the matrix elements on the right hand side of the last equation
can be nonzero. Assume $\psi _{j}\mu =0$, then it follows from the previous
proposition that%
\begin{eqnarray*}
q^{d}C_{\lambda ,\mu }^{\nu ,d}(k,N) &=&\langle \nu ,S_\lambda\psi _{j}\psi _{j}^{\ast }\mu \rangle  \\
&=&\sum_{r=0}^{\lambda _{1}^{t}}(-1)^{r}\sum_{\lambda ^{t}/\rho
^{t}=(r)}\langle \nu ,\psi _{j-r}S^\prime_{\rho }\psi
_{j}^{\ast }\mu \rangle  \\
&=&q^{d}\sum_{r=0}^{\ell (\lambda )}(-1)^{r+d+\op{n} _{j-1}(\mu )+\op{n}
_{j-r-1}(\nu )}\sum_{\lambda /\rho =(1^{r})}C_{\rho ,\psi _{j}^{\ast }\mu
}^{\psi _{j-r}^{\ast }\nu ,d_{r}}(k+1,N)\;.
\end{eqnarray*}%
Here we have used that%
\begin{equation*}
|\lambda |=\sum_{i=1}^{k}\ell _{i}-\frac{k(k+1)}{2}\qquad \text{and\qquad }%
|\psi _{j}^{\ast }\lambda |=|\lambda |+j-k-1,
\end{equation*}%
from which one deduces that for $j-r>0$%
\begin{equation*}
|\rho |+|\psi _{j}^{\ast }\mu |-|\psi _{j-r}^{\ast }\nu |=|\lambda |+|\mu
|-|\nu |=dN\;.
\end{equation*}%
If $j-r<0$ it follows from the quasi-periodic boundary conditions that
\begin{equation*}
\langle \nu ,\psi _{j-r}S^\prime_{\rho }\psi _{j}^{\ast }\mu
\rangle =(-1)^{k-1}q\langle \nu ,\psi _{j-r+N}S^\prime_{\rho
}\psi _{j}^{\ast }\mu \rangle
\end{equation*}%
and the degree is decreased by one. Since we set $\op{n} _{a+N}=\op{n} _{a}+k$
the overall sign works out correctly.

The second case, $\psi _{j}^{\ast }\mu =0$, leads to the second recursion
formula via a similar computation. To determine the equality of the degrees
one now uses the relation $|\psi _{j}\lambda |=|\lambda |-j+k$.
\end{proof}

\begin{example}{\rm
Consider once more the product expansion from Example \ref{ExGW}. Choose $\nu
=(3,2,1,1)$, then $C_{\lambda \mu }^{\nu ,1}=2$. Since $w(\mu )=1010110$ and $w(\nu )=0110101$, $\psi
_{j}^{\ast }\mu $ is only non-vanishing for $j=2,4,7$,%
\begin{equation*}
\psi _{2}^{\ast }\mu =-(1,1,0,0,0),\quad \psi _{4}^{\ast }\mu
=(1,1,1,1,0),\quad \psi _{7}^{\ast }\mu =(2,2,2,1,0)\;.
\end{equation*}%
In contrast, $\psi _{j-r}^{\ast }\nu \neq 0$ only if $j-r=1,4,6,$%
\begin{equation*}
\psi _{6}^{\ast }\nu =-(2,2,2,1,1),\quad \psi _{1}^{\ast }\nu
=(2,1,0,0,0),\quad \psi _{4}^{\ast }\nu =(2,1,1,1,1)\;,
\end{equation*}%
whence the only allowed value for $j=2,7$ is $r=3$, while for $j=4$ we have
the possible values $r=0,3$. The corresponding partitions $\rho $ in the
non-trivial skew diagrams $\lambda /\rho $ are $\rho _{1}=(2,2,1,1,0)$ and $%
\rho _{2}=(2,2,2,0,0)$. We obtain the following
identities between Gromov-Witten invariants, where the first equality in each line follows from \eqref{GWrec1}, while the actual numerical values are calculated with the help of a computer using \eqref{BVI}:%
\begin{eqnarray*}
j =2:\; C_{\lambda \mu }^{\nu ,1}(4,7) &=& C_{\rho _{1}\psi _{2}^{\ast
}\mu }^{\psi _{6}^{\ast }\nu ,0}(5,7)+C_{\rho _{2}\psi _{2}^{\ast }\mu
}^{\psi _{6}^{\ast }\nu ,0}(5,7)=1+1=2, \\
j =4:\; C_{\lambda \mu }^{\nu ,1}(4,7) &=& C_{\lambda \psi _{4}^{\ast
}\mu }^{\psi _{4}^{\ast }\nu ,1}(5,7)+C_{\rho _{1}\psi _{4}^{\ast }\mu
}^{\psi _{1}^{\ast }\nu ,1}(5,7)+C_{\rho _{2}\psi _{4}^{\ast }\mu }^{\psi
_{1}^{\ast }\nu ,1}(5,7)=0+1+1=2, \\
j =7:\; C_{\lambda \mu }^{\nu ,1}(4,7)&=& C_{\rho _{1}\psi _{7}^{\ast
}\mu }^{\psi _{4}^{\ast }\nu ,1}(5,7)+C_{\rho _{2}\psi _{7}^{\ast }\mu
}^{\psi _{4}^{\ast }\nu ,1}(5,7)=1+1=2\;.
\end{eqnarray*}%
Note that for $j=2$ we have used the
quasi-periodic boundary conditions, $\psi _{-1}=q\psi _{6}e^{\iota\pi K}$ to
determine the overall sign factor and degree.
}
\end{example}

\section*{Part III: The `Boson-Fermion-correspondence'}
The following table summarizes the analogy of the constructions on both sides. We call it the `Boson-Fermion-correspondence' (although a precise representation theoretic statement will appear in \cite{KS}).\\

\begin{tabular}{l|l}
fusion ring $\cF(\mathfrak{\hat{sl}}(n))_k$&quantum cohomology $qH^\bullet(\op{Gr}_{k,k+n})$\\
\hline\\
{\bf integrable system:}\\
(bosonic) configuration of $k$ particles&(fermionic) configurations of $k$ particles\\
on a circle with $n$ sites&on a circle with $N=n+k$ sites\\
\\
{\bf creation and annihilation operators:}\\
phase algebra&Clifford algebra\\
\\
{\bf noncommuting variables:}\\
Local affine plactic algebra& affine nil-Temperley-Lieb algebra\\
\\
{\bf crystal structure:}\\
Crystal of the $k^{\text{th}}$ symmetric power&Crystal of the $k^{\text{th}}$ exterior power
\\
\\
{\bf Bethe Ansatz equations:}\\
$x_{1}^{n+k}=\cdots =x_{k}^{n+k}=(-1)^{k-1}ze_{k}(x_{1},\ldots,x_{k})$
&$
y_{1}^{N}=\cdots =y_{k}^{N}=(-1)^{k-1}q.
$
\\
\\
{\bf Bethe roots} ($\zeta =\exp \frac{2\pi i}{k+n}$)\\
$x_\sigma=z^{\frac{1}{n}}\zeta ^{\frac{|\sigma|}{n}} \left(\zeta ^{I_{1}},\ldots,\zeta^{I_{k}}\right), \sigma\in\mathfrak{P}_{\leq n-1,k}$
&
$y_\sigma=q^{\frac{1}{N}}\left(\zeta ^{I_{1}},\ldots,\zeta^{I_{k}}\right),
\sigma\in\mathfrak{P}_{\leq k,n}$,\\
\\
{\bf ring structure}\\
$\la\circledast\mu=s_\la(a_0,a_1,\ldots a_{n-1})\mu$
&
$\la\star\mu=s_\la(u_0,u_1,\ldots u_{n-1})$\\
Verlinde formula&Bertram-Vafa-Intrilligator formula
\end{tabular}

\section*{Appendix: local plactic algebra and tableaux}
\label{Appendix}
\label{algorithm}
We found the following explicit algorithm which produces a sequence of (generalised) tableaux, $\mathcal{T}\ni T\rightsquigarrow \hat{T}\rightsquigarrow  T^\prime\rightsquigarrow  D(T)\in\mathcal{T}_{\text str}$ (using the notation from Proposition \ref{finiteplactic} \eqref{2}) such that the following statement is true.

\begin{proposition}
The algorithm associates to each tableau $T$ a tableau $D(T)$ in $\mathcal{T}_{\text str}$ with equivalent column word under the relations \eqref{PLfin1} and \eqref{PLfin2}.
\end{proposition}

In particular, from $D(T)$ one can read off the corresponding basis element in $B$.
\begin{definition}
Assume we are given a tableau $T$ (or any more general diagram, where the following definition makes sense). We say a column $C$ of T {\em has property (S)}  if there are no two numbers $i<j$, appearing in $C$ in two consecutive rows and satisfying $j-i>1$. (That means the numbers in C form an increasing sequence of consecutive numbers).
\end{definition}
{\bf The algorithm}:
\subsubsection*{Step 0}
If the tableau is already in $\mathcal{T}_{\text str}$ then there is nothing to do.
\subsubsection*{Step 1: Creating $\hat{T}$}
Read the tableau $T$ column by column from right to left (and from top to bottom)  until the first column $C$ which has not property (S). Let $i,j$ be the first pair of numbers destroying property (S).

Now do the following procedure {\em basic procedure} $P(C,j)$ (associated with the pair $(C,j)$): cross out the number $j$ and consider the column, $D$, next to $C$ on the right, deal with the first applicable situation from the following list, where the term `contains' stands for `contains an uncrossed'.
\begin{enumerate}[(1)]
\item \label{A1} $D$ is empty. Place $j$ in the first row of $D$.
\item \label{A2a} $D$ contains $j-1$ and $j$. Do $P(D,j)$.
\item \label{A2b} $D$ contains $j-1$, not $j$. Place $j$ below $j-1$, (replacing any crossed out number there) and pushing down all greater numbers by one row.
\item \label{A3a} $D$ contains $j+1$, not $j$. Then $j+1$ is necessarily on the top of $D$ (due to the choice of $i$ and $j$). Replace $j+1$ by $j$, and do $P(D,j+1)$.
\item \label{A3b} $D$ contains $j+1$ and $j$. Do $P(D,j)$.
\item  $D$ contains neither $j+1$ nor $j-1$ and is not empty.
\begin{enumerate}
\item \label{A4a} If it contains just $j$ then insert an extra copy of the column $D$ to the right of $D$.
\item \label{A4b} If all entries in $D$ are greater than $j$ then replace the entry, $x$, in the first row by $j$, and do $P(D,x)$.
\item \label{A4c} If all entries in $D$ are less than $j$ do $P(D,j)$.
\end{enumerate}
\end{enumerate}
Since $T$ has a finite number of rows and columns, we are done after applying finitely many basic procedures. The multi-set of numbers which are not crossed out will be the same as at the beginning, but some numbers could have been moved to the right or up. There are no numbers crossed out in the first row.

Now repeat the above rule with all (not crossed out) pairs $i$, $j$ as long as possible.
The whole procedure stops after finitely many steps, latest when all numbers which are not crossed out appear in the first row.
We claim that this produces a tableau $\hat{T}$ with some numbers crossed out (for a proof see below). By construction, each step does not change the equivalence class of column words (under the relations \eqref{PLfin1} and \eqref{PLfin2}.

\subsubsection*{Step 2: Creating $T'$ and $D(T)$} Now remove the boxes with crossed out numbers to obtain $T'$. The first row is weakly increasing, the columns are connected and are increasing sequences $i (i+1)(i+2)\ldots$ of consecutive numbers. Thanks to property (S) the local plactic relations make sure that the corresponding column word is in the same equivalence class as the column word of $T$.
 Finally, the tableau $D(T)$ is obtained from $T'$ by moving all the boxes to the left as far as possible. The specific shape of $T'$ makes sure that $D(T)$ is in fact a tableau (since the rows are of course still weakly increasing and the columns are not necessarily sequences anymore, but definitely strictly increasing). The column word of $D(T)$ is obviously equivalent to the column word of $T'$ via the relations \eqref{PLfin1}. This is the end of the algorithm.

\begin{example}
\label{algo}
{\rm We illustrate the algorithm with an example
\begin{eqnarray*}
\begin{array}{llll}
\Yvcentermath1{T=}&\Yvcentermath1{\young(111223346\ten,223334567,9)}&\rightsquigarrow&
\Yvcentermath1\young(1112233466\ten,2233345\cansix7,9)\\
\\
&\Yvcentermath1\young(1112233466\ten,223334\canfive57,9)&\rightsquigarrow&
\Yvcentermath1\young(11122333466\ten,22\canthree334\canfive\canfive57,9)\\
\\
\Yvcentermath1{\hat{T}=}&
\Yvcentermath1\young(111223334669\ten,22\canthree334\canfive\canfive57,\cannine)\\
\end{array}
\end{eqnarray*}
In the first step we applied \eqref{A3b}, \eqref{A4b}, \eqref{A1}, then \eqref{A2b}, then twice \eqref{A2a} followed by \eqref{A3b} and \eqref{A4a}, and finally ten times \eqref{A4c} followed by \eqref{A3a} and \eqref{A1}. Moreover,
\begin{eqnarray*}
\Yvcentermath1 T'&=&\Yvcentermath1\young(111,22)\young(22,33)\young(333,4)\young(4,5)\young(669\ten,7)\\
\Yvcentermath1 D(T)&=&\Yvcentermath1\young(111223334669\ten,2233457)
\end{eqnarray*}
The column word for $D(T)$ is $a_2a_1a_2a_1a_3a_1a_3a_2a_4a_2a_5a_3a_7a_3a_3a_4a_6a_6a_9a_{10}$. Under the local plactic relations, \eqref{PLfin1} and \eqref{PLfin2}, the corresponding element in $B$ is then $$a_1(a_2a_1)(a_2a_1)(a_3a_2)(a_3a_2)(a_4a_3)(a_3)(a_3)(a_5a_4)(a_6)(a_7a_6)(a_9)(a_{10}),$$
which is obviously equivalent to $D(T)$ under the Robinson-Schensted correspondence.
We demonstrate the reordering for the subsequence $a_2a_1a_2a_1a_3a_1a_3a_2a_4a_2\ldots$ of the column word: apply \eqref{PLfin1} to obtain $a_2a_1a_2a_1a_1a_3a_3a_2a_2a_4\ldots$ and then repeatedly \eqref{PLfin2} to arrive at $a_2a_1a_1a_2a_1a_3a_2a_3a_2a_4\ldots = a_1a_2a_1a_2a_1a_3a_2a_3a_2a_4\ldots$ Continuing in the same manner one derives the above basis word.}
\end{example}

\begin{proof} [Proof that $\hat{T}$ is a tableau]
We have to show that the rules from Step 1 produce a tableaux. We do this by carefully examining what happens in the different cases (1) -(6)(c).
So let $(i,j)$ be as in Step 1. The cases (1) and (4) are obvious, there is nothing to check, since the top entry in $C$ is smaller
than $j$. In case \eqref{A2b} we place $j$ below $j-1$, the last entry in
the column thanks to property $(S)$. The entry $a$ to the left of $j-1$ is
at most $j-1$, hence there must be an entry $b\leq j$ below $a$. The entry $%
a^{\prime}$ to the right of $j-1$ is at least $j-1$, hence any possible
entry $b^{\prime}$ below $a^{\prime}$ is at least $j$. Cases \eqref{A4a} and %
\eqref{A4b} are obvious. So placing $j$ always produces a new tableau, but
(disregarding the crossed out ones) possibly one entry removed.

\begin{enumerate}

\item[(i)]\label{A2aappl} Assume we removed $j$ via case \eqref{A2a} applied to
some column $C^{\prime}$ which then has to contain at least $j$ and $j-1$.
Applying then cases \eqref{A1}, \eqref{A2b}, \eqref{A3a} or \eqref{A4b}
cause no problems. Case \eqref{A4a} is obviously producing a tableau.

\item[(ii)] Assume we removed $j+1$ via case \eqref{A3a} from some column $%
C^{\prime}$. Again case \eqref{A1} is obvious. Then case \eqref{A2b} is not
applicable, since $C^{\prime}$ originally contained $j+1$, hence the entry
to the right of it is at least $j+1$. In case \eqref{A3a} we originally have
$j+1$ at the top of column $C^{\prime}$ and $j+2$ to the right of it. This
gets changed into $j$ and $j+1$ which produces a tableau. Case \eqref{A4a} is obviously causing no problems. In case \eqref{A4b} we originally had $j+1$ on the top of $C^{\prime}$
and all entries in the column to the right of it are greater than $j+2$. In
particular, we can replace the top entry by $j+1$.

\item[(iii)] Assume we applied case \eqref{A3b} for $(C^{\prime},j)$. Then applying %
\eqref{A1} is fine. For case \eqref{A2b} we argue as in \eqref{A2aappl}. For
cases \eqref{A3a} and \eqref{A4b} note that the top entry in $C^{\prime}$ is
at most $j$. Case \eqref{A4a} is clear.

\item[(iv)] Assume we applied case \eqref{A4b} for $(C^{\prime},j)$ with $j$
replacing $x$. Then applying \eqref{A1} is fine, because $x>j$. Case %
\eqref{A2b} is not applicable. Cases \eqref{A3a} and \eqref{A4a} are clear.
In case \eqref{A4b} we originally have $x$ on top of column $C^{\prime}$.
There might be some $a$ below it. Let $z$ be to the right of $x$ and
possibly some $b$ below $z$. Then $x$ gets replaced by $j$ and $z$ by $x$,
hence the first row is still weakly increasing. Since $j<x<a$ and $b>z\leq x$%
, the columns are still strictly increasing.

\item[(v)] Assume we applied case \eqref{A4c}. Then place $j$ at the bottom of
column $C^{\prime}$ and consider only the diagram containing the columns to
the right of $C^{\prime}$, including $C^{\prime}$. By induction on the
number of columns, we are done.
\end{enumerate}
\end{proof}

\bibliographystyle{plain}

\end{document}